\documentclass[final,onefignum,onetabnum,letterpaper]{siamonline190516}
\usepackage[scale=0.8]{geometry} 

\usepackage{soul}
\usepackage{lipsum}
\usepackage{amsfonts}
\usepackage{epstopdf}
\ifpdf
  \DeclareGraphicsExtensions{.eps,.pdf,.png,.jpg}
\else
  \DeclareGraphicsExtensions{.eps}
\fi

\usepackage{datetime}
\newdateformat{monthyeardate}{%
	\monthname[\THEMONTH] \THEDAY, \THEYEAR}

\usepackage{academicons}
\usepackage{xcolor}
\renewcommand{\orcid}[1]{\href{https://orcid.org/#1}{\textcolor[HTML]{A6CE39}{orcid.org/#1}}}

\usepackage{amsmath}
\allowdisplaybreaks
\usepackage{amssymb}
\usepackage{commath}
\usepackage{mathtools}
\usepackage{bbm}
\usepackage{bm}

\usepackage{color}
\usepackage{graphicx}
\usepackage[small]{caption}
\usepackage{subcaption}

\usepackage{relsize}
\usepackage{adjustbox}
\usepackage{algorithm}
\usepackage{algorithmicx}
\usepackage{algpseudocode}
\usepackage{booktabs}
\usepackage{tikz}
\usepackage{pifont}
\usetikzlibrary{bayesnet} 

\algblock{Indent}{EndIndent}
\algtext*{Indent}
\algtext*{EndIndent}

\usepackage{enumitem}
\setlist[enumerate]{leftmargin=.5in}
\setlist[itemize]{leftmargin=.5in}

\newsiamremark{remark}{Remark}
\newtheorem{example}{Example}
\newsiamremark{hypothesis}{Hypothesis}
\crefname{hypothesis}{Hypothesis}{Hypotheses}
\newsiamthm{claim}{Claim}

\headers{Residual prior for signal recovery}{Y.\ Xiao and A.\ Gelb}

\title{Joint Signal Recovery and Uncertainty Quantification via the Residual Prior Transform
\thanks{\monthyeardate\today  
}}
\author{ 
Yao Xiao\thanks{
Department of Mathematics, Dartmouth College, USA 
(\email{yao.xiao@dartmouth.com}, \orcid{0000-0002-0850-9624})
}
\and
Anne Gelb\thanks{
Department of Mathematics, Dartmouth College, USA 
(\email{anne.E.gelb@dartmouth.edu}, \orcid{0000-0002-9219-4572})
} 
}

\usepackage{amsopn}

\usepackage[normalem]{ulem}

\newcommand{\MAP}{\mathrm{MAP}}
\newcommand{\MMV}{\mathrm{MMV}}

\DeclareMathOperator{\diag}{diag}

\DeclareMathOperator*{\argmin}{arg\,min}

\DeclareMathOperator*{\kernel}{kernel}

\newcommand{\R}{\mathbb{R}} 
\newcommand{\C}{\mathbb{C}}

\usepackage{comment}

\ifpdf
\hypersetup{
	pdftitle={Xiao2024JSHBL},
 	pdfauthor={Yao Xiao}
}
\fi

\begin{document}

\maketitle

\begin{abstract}
	Conventional priors used for signal recovery are often limited by the assumption that the type of a signal's variability, such as piecewise constant or linear behavior, is known and fixed. This assumption is problematic for complex signals that exhibit different behaviors across the domain. The recently developed {\em residual transform operator} effectively reduces such variability-dependent error within the LASSO regression framework. Importantly, it does not require prior information regarding structure of the underlying signal. This paper reformulates the residual transform operator into a new prior within a hierarchical Bayesian framework.  In so doing, it unlocks two powerful new capabilities. First, it enables principled uncertainty quantification, providing robust credible intervals for the recovered signal, and second, it provides a natural framework for the joint recovery of signals from multimodal measurements by coherently fusing information from disparate data sources. Numerical experiments demonstrate that the residual prior yields high-fidelity signal and image recovery from multimodal data while providing robust uncertainty quantification.
\end{abstract}

\begin{keywords}
	Spatial sequence of measurements, 
  	joint sparsity,
	hierarchical Bayesian learning, 
	conditionally Gaussian priors, 
	(generalized) gamma hyperpriors
\end{keywords}

\begin{AMS}
	65F22, 
	62F15, 
	65K10, 
	68U10 
\end{AMS}

\begin{DOI}
    Not yet assigned
\end{DOI}

\section{Introduction} 
\label{sec:introduction}

Signal recovery from noisy and undersampled data is important across a wide variety of application areas. Sparsity-promoting techniques have demonstrated great success, with foundational ideas developed in compressive sensing \cite{candes2006robust,candes2006stable,chartrand2008iteratively,donoho2006compressed} and elegantly formalized within hierarchical frameworks like sparse Bayesian learning (SBL) \cite{tipping2001sparse,wipf2004sparse,wipf2007empirical}. Regardless of the specific framework, the objective is to find a point estimate that best balances data fidelity with prior knowledge regarding the underlying signal's structural information, generally the sparsity in certain transformed domains. However, a common feature across many of these methods is that the underlying assumption regarding the sparse-promoting prior operator is typically pre-determined and fixed. This reliance on a single, fixed model of variability, such as assuming purely piecewise constant or linear behavior, becomes problematic for more complex signals. Such assumptions are particularly inappropriate for functions that are piecewise polynomial, where smooth behavior can vary substantially from one interval to the next. Such a rigid assumption is also problematic in applications  such as scientific data compression, where efficient and robust representations of piecewise smooth, highly variable functions are paramount.

A novel {\em residual transform operator} was introduced in \cite{xiao2025new} to address the limitations of using fixed order transform operators in the context of Lasso regression. The foundational idea is that for any piecewise smooth signal, one can design two different operators, $\Phi_1$ and  $\Phi_2$, that respond similarly to the signal's underlying spatial variability, even when that variability is unknown. While either operator used alone may yield a poor choice of regularizer, the residual transform operator on the underlying (discretized) signal ${\bm f}$, given by $R{\bm f} = (\Phi_1 -\Phi_2){\bm f}$, is small by construction, especially in smooth regions. This provides a more robust method for promoting sparsity without making rigid assumptions about the signal's variable behavior. The residual transform operator was  successfully established as a powerful tool for obtaining a single point estimate in compressive sensing formulation \cite{xiao2025new}.

This investigation builds upon this concept by incorporating the regularizing residual transform operator into the SBL framework, namely, by using the residual transform operator to define the sparsity-promoting prior in the hierarchical Bayesian model. To differentiate our prior from the LASSO regression regularizer in \cite{xiao2025new}, we will refer to the corresponding sparsity-promoting transform operator as the {\em residual prior transform}. Our approach provides several significant advantages over its original deterministic formulation. First, our hierarchical Bayesian framework enables robust Uncertainty Quantification (UQ) by yielding a full conditional posterior distribution for the signal in addition to a single MAP point estimate. Credible intervals can then computed to assess the reconstruction's reliability \cite{mackay1992bayesian,glaubitz2025efficient}. Moreover, by utilizing uninformative hyperpriors, the SBL framework automatically infers the appropriate level of regularization directly from the data, which eliminates the need for the manual parameter tuning often required in compressive sensing. Finally, this work extends the application of the residual prior transform to complex multimodal problems, demonstrating its ability to leverage joint sparsity across a set of related spatial datasets. 


This investigation systematically evaluates the performance and robustness of the residual prior transform, demonstrating its advantages over standard first-order differencing prior transforms, particularly in capturing complex signal variability. A central goal of this work is to extend the residual transform operator beyond its original single measurement context to more challenging and realistic multimodal and two-dimensional scenarios. We will explore two critical paradigms: first, the joint recovery of multiple, distinct signals that are assumed to share a common sparsity structure, and second, the reconstruction of a single underlying signal from a diverse set of observation sources. Through our numerical experiments, which span from synthetic functions with mixed variable intervals and sparse characteristics to real world applications in Synthetic Aperture Radar (SAR) and traffic surveillance imaging, we will show that the residual prior transform consistently and more faithfully reconstructs regions of high variability where the standard first order differencing prior transform often fails. These results will establish that the residual prior transform is a more flexible and powerful choice for a wide range of realistic inverse problems encountered in practice.


\subsection*{Outline} 
The remainder of this paper is organized as follows. We begin in \Cref{sec:prelim} by reviewing the general hierarchical sparse Bayesian learning framework, including the update steps for both individual and joint recovery. This section also presents a local differencing operator and a global Fourier-based operator that will serve as the core ingredients for our proposed prior transform. The construction and error analysis of our new residual prior transform $R= \Phi_1-\Phi_2$ is formally developed in \Cref{sec:residualprior}. In \Cref{sec:numerics}, we demonstrate the effectiveness and robustness of this new prior through a series of numerical experiments on synthetic and real-world data. Finally, \Cref{sec:summary} provides concluding remarks and discusses potential directions for future work.

\subsection*{Notation} 

We use normal ($X$) and boldface capital letters ($\mathbf{X}$) to represent scalar and vector random variables, respectively. The notation $X\sim \pi$ indicates that a random variable $X$ is distributed according to the probability density $\pi$. For a collection of $L$ random variables $\mathbf{X}_1, \dots, \mathbf{X}_L$, we use the shorthand notation $\mathbf{X}_{1:L}=(\mathbf{X}_1, \dots, \mathbf{X}_L)$. This same notation can be applied to dummy variables, such as $\boldsymbol{x} \in \R^n$.

\section{Preliminaries} 
\label{sec:prelim} 

Consider a sequence of one-dimensional functions $\{f_l: [-\pi,\pi]\to\mathbb R\}_{l=1}^L$. Each $f_l$ corresponds to a distinct scenario, such as what might occur in a hyper-spectral data collection or in measuring physical quantities (e.g.~wind, temperature, or pressure) \cite{shaw2003spectral,grahn2007techniques,adao2017hyperspectral}, while collectively sharing a common set of abrupt changes (discontinuities). Importantly, the locations of these discontinuities are not known beforehand.

Our goal is to reconstruct {\em each} signal, $f_l(s)$, for $l = 1,\dots, L$, from its corresponding observations collected at a fixed set of $n$ equally-spaced grid points. These points are defined as $s_j = -\pi + j\Delta s$, where $\Delta s = \frac{2\pi}{n}$ for $j=0,\dots,n-1$. In this work, we model the noisy observations, denoted by $\mathbf{y}_l\in\R^{m_l}$, using the following forward process:   
\begin{equation}
    \label{eq: data_acquisition}
    \bm y_l = F_lf_l(\bm s) + \boldsymbol\epsilon_l = F_l\bm{x}_l + \boldsymbol\epsilon_l, \quad l=1,\dots,L.
\end{equation}
In this formulation, $\bm{x}_l \in \mathbb R^n$ represents a discretized version of the true signal $f_l$. This vector is created by sampling $f_l$ at the $n$ uniformly-spaced points within the interval $[- \pi,\pi]$. Our task is to recover this discretized vector $\bm{x}_l$, which serves as our reconstructed approximation of the continuous signal $f_l$. The vector $\bm y_l$ is collected from different sources, and we assume that the linear forward operators $F_l\in\R^{m_l\times n}$ are known.
Finally, the observations $\bm{y}_l$ are contaminated by additive, independent, and identically distributed (i.i.d.) Gaussian noise with zero mean and inverse variance $\alpha>0$ given by 
\begin{equation}
    \label{eq: iid_noise}\boldsymbol{\epsilon}_l\sim\mathcal N\left(0, \alpha^{-1}_l I_{m_l}\right),
\end{equation}
yielding the covariance matrix $\alpha^{-1}_l I_{m_l}$. This noise model extends directly to the 2D case by treating the 2D image as a vectorized signal and assuming the additive noise is i.i.d. across all pixels.

\subsection{The hierarchical Bayesian model}
\label{sub:prelim_model}

This work is motivated by the difficulties of recovering piecewise smooth signals with varying orders of smoothness. To address this limitation, we adapt a sparsity-promoting residual transform operator, originally developed for Lasso regression \cite{xiao2025new}, for use as a conditional Gaussian prior transform within the hierarchical sparse Bayesian learning framework. To highlight the effectiveness of the proposed residual prior transform, we benchmark its performance against the standard first-order differencing operator, essentially the total variation (TV), a widely used baseline for such problems. A brief review of the hierarchical sparse Bayesian learning model used in this investigation is provided below for self containment.

\subsubsection{Sparsity-promoting hierarchical learning}\label{sub:BCDalg}
We consider the generalized sparse Bayesian learning (GSBL) algorithm proposed in \cite{glaubitz2022generalized} which computes the posterior mean for each $\bm x_l$ in \eqref{eq: data_acquisition}. Note that each posterior mean is obtained {\em separately}, that is, no mutual information between measurements is used in the standard GSBL algorithm.

To demonstrate the basic structure and application of the GSBL algorithm, we adopt the following hierarchical sparse Bayesian learning model: For $l=1,\dots,L$, 
\begin{subequations}\label{eq:model1} 
\begin{equation}
	\mathbf{Y}_l | \mathbf{X}_l = \bm{x}_l 
		 \sim \mathcal{N}( F_l \bm{x}_l, \alpha^{-1}_lI_{m_l} ) \label{eq:model1likelihood}\end{equation}
	\begin{equation} \Phi_l \mathbf{X}_l | \boldsymbol{\Theta}_l = \boldsymbol{\theta}_l 
		 \sim \mathcal{N}(\mathbf{0},D_{\boldsymbol{\theta}_l}), \label{eq:model1prior}\end{equation}
	\begin{equation} (\Theta_l)_k 
		 \sim \Gamma((\beta_l)_k, (\vartheta_l)_k), \quad k=1,\dots,K_l, \label{eq:model1hypertheta}\end{equation}
\end{subequations} 
Here, \cref{eq:model1likelihood} describes the likelihood model, $\Phi_l \in \mathbb{R}^{K_l \times n}$ in \cref{eq:model1prior} encodes our prior belief that $\mathbf{X}$ is sparse under the transform $\Phi$ with a diagonal covariance matrix  $D_{\boldsymbol{\theta}_l} = \operatorname{diag}(\boldsymbol{\theta}_l) \in \mathbb{R}^{K_l\times K_l}$.  The hyperprior $\boldsymbol{\theta}_l = [(\theta_l)_1, \dots, (\theta_l)_{K_l}]$ in \cref{eq:model1hypertheta} is modeled by the gamma distribution $\Gamma$. For each $l=1,\dots,L$, setting $\{(\beta_l)_k\}_{k=1}^{K_l} = 1$ and $\{(\vartheta_l)_k\}_{k=1}^{K_l} \approx 0$ removes the moderating influence of the hyperprior parameter vector, which enables $\boldsymbol\theta_l$ to fluctuate wildly according to the data. 
The specific values used for the hyperprior parameter vector, $\{(\vartheta_l)_k\}_{k=1}^{K_l}$, are provided in \Cref{sec:numerics}.
It is important to emphasize that \cref{eq:model1} describes the model used to {\em individually} recover the map estimate for {\em each} data acquisition, $l = 1,\dots, L$, that is,  no information is shared. 

To reflect the practical challenge where signal smoothness is often unknown and variable, we deliberately apply a single, consistent prior operator $\Phi_l=\Phi$ across all $L$ measurements.   This approach underscores the power of the proposed prior transform, which is designed to be effective without being fine-tuned to the specific variability within each measurement. This key advantage removes the burden of selecting and tuning different operators for each data set, which is a classic hurdle for conventional prior transforms. Consequently, the residual prior transform is exceptionally well-suited for multimodal joint recovery algorithms, as it can be incorporated seamlessly without the need for modality-specific parameter adjustments.

Corresponding  to \cref{eq:model1} are the likelihood, prior, and hyperprior densities,
\begin{subequations}\label{eq:model1_PDFs} 
\begin{equation}
	\pi_{\mathbf{Y}_l|\mathbf{X}_l} ( \bm{y}_l | \bm{x}_{l}) 
		 \propto \exp\left( - \frac{\alpha_l}{2} \| F_l \bm{x}_l - \bm{y}_l \|_2^2 \right), \label{eq:model1pdflikelihood}\end{equation} 
\begin{equation}	\pi_{\mathbf{X}_l|\boldsymbol\Theta_l} ( \bm{x}_l | \boldsymbol{\theta}_l ) 
		 \propto \det( D_{\boldsymbol{\theta}_l} )^{1/2} \exp\left( - \frac{1}{2} \| D_{\boldsymbol{\theta}_l}^{1/2} \Phi_l \bm{x}_l \|_2^2 \right), \label{eq:model1pdfprior}\end{equation}
\begin{equation}
\pi_{\boldsymbol{\Theta}_l} (\boldsymbol{\theta}_l)
=\prod_{k=1}^{K_l}\Gamma((\theta_l)_k|\beta_l, \vartheta_l)
		 \propto  \det(D_{\theta_l})^{\beta_l - 1} 
         \exp\left( - \sum_{k=1}^{K_l} \frac{(\theta_l)_k}{\vartheta_l}  \right).
\label{eq:model1pdfhypertheta}
\end{equation}
\end{subequations}

\begin{remark}
To relax the i.i.d. noise assumption in \eqref{eq: iid_noise}, the noise can be modeled as $\boldsymbol{\epsilon}_l\sim\mathcal{N}(0,A^{-1}_l)$, where $A_l=\diag(\boldsymbol\alpha_l)=\diag([(\alpha_l)_1,\dots,(\alpha_l)_{m_l}])$ is a diagonal, positive-definite inverse noise covariance matrix. 
    This yields a more generalized likelihood function than  \eqref{eq:model1pdflikelihood}: 
\[ \pi_{\mathbf Y_l|\mathbf X_l}( \bm{y}_l | \bm{x}_l) 
		 \propto \exp\left( - \frac{1}{2} ( F_l \bm{x}_l - \bm{y}_l)^TA_l  ( F_l \bm{x}_l - \bm{y}_l)\right) \]
Note that this model simplifies back to the i.i.d. case described in \eqref{eq: iid_noise} with its corresponding likelihood function \eqref{eq:model1pdflikelihood} if all inverse variances are equal ($(\alpha_l)_1=\dots=(\alpha_l)_{m_l}=\alpha_l$). 
\end{remark}

The relationship between the unknown discretized signal $\bm x$, its hyperparameter vector  $\boldsymbol{\theta}$, and the observation $\bm y$ is governed by Bayes' Theorem. In general formulation, this theorem provides the joint posterior density function as a product of the likelihood, prior, and hyperprior as
\begin{equation}\label{eq:model1_Bayes} 
	\pi_{\mathbf{X},\boldsymbol{\Theta}|\mathbf{Y}}( \bm{x}, \boldsymbol{\theta} | \bm{y} ) 
		\propto \pi_{\mathbf{Y}|\mathbf{X}}( \bm{y} | \bm{x} ) 
			\pi_{\mathbf{X}|\boldsymbol{\Theta}}( \bm{x} | \boldsymbol{\theta} ) \pi_{\boldsymbol{\Theta}}(\boldsymbol{\theta}).
\end{equation}

Direct substitution of each equation of \cref{eq:model1_PDFs} into \cref{eq:model1_Bayes} yields
\begin{equation}\label{eq:model1_posterior} 
	\pi_{\mathbf{X}_l, \boldsymbol{\Theta}_l | \mathbf{Y}_l}( \bm{x}_l, \boldsymbol{\theta}_l | \bm{y}_l ) 
		\propto 
        \det( D_{\boldsymbol{\theta}_l} )^{ \beta_l - 1/2}\exp\left( - \frac{ \alpha_l}{2} \| F_l \bm{x}_l - \bm{y}_l \|_2^2 - \frac{1}{2} \| D_{\boldsymbol{\theta}_l}^{1/2} \Phi_l \bm{x}_l \|_2^2 - \sum_{k=1}^{K_l} \frac{(\theta_l)_k}{\vartheta_l}  \right) .
\end{equation}

\subsubsection{MAP estimation}
\label{sub:prelim_model1_MAP}

We now address Bayesian inference for the hierarchical sparse Bayesian learning model in \cref{eq:model1_Bayes}. 
To this end, a common strategy is to solve for the \emph{Maximum a Posteriori (MAP) estimate} $(\bm{x}^{\MAP}_l,\boldsymbol{\theta}^{\MAP}_l)$ for given measurements $\bm{y}_l$, which is the maximizer of the joint posterior density \cref{eq:model1_posterior}. 
Equivalently, the MAP estimate is the minimizer of the negative logarithm of the posterior, i.e.,  
\begin{equation}\label{eq:model1_MAP_estimate}
	(\bm{x}^{\MAP}_l,\boldsymbol{\theta}_l^{\MAP}) 
		= \argmin_{ \bm{x}_l, \boldsymbol{\theta}_l } \left\{ \mathcal{G}( \bm{x}_l, \boldsymbol{\theta}_l ) \right\}, 
\end{equation} 
where the objective function $\mathcal{G}$ is 
\begin{equation}\label{eq:model1_obj_fun} 
	\mathcal{G}( \bm{x}_l, \boldsymbol{\theta}_l ) 
		= - \log \pi_{ \mathbf{X}_l, \boldsymbol{\Theta}_l  | \mathbf{Y}_l }( \bm{x}_l, \boldsymbol{\theta}_l| \bm{y}_l ).
\end{equation} 
Substituting \cref{eq:model1_posterior} into \cref{eq:model1_obj_fun} leads to
\begin{equation}\label{eq:G}
\mathcal{G}( \bm{x}_l, \boldsymbol{\theta}_l ) = \frac{\alpha_l}{2}\| F_l \bm{x}_l - \bm{y}_l \|_2^2 + \frac{1}{2}\bigl\| D_{\boldsymbol{\theta}_l}^{1/2} \Phi_l \bm{x}_l\bigr\|_2^2+ \sum_{k=1}^{K_l} \frac{(\theta_l)_k}{\vartheta_l}+ \sum_{k=1}^{K_l} \Bigl(\beta_l - 1/2\Bigr)\log\!\bigl((\theta_l)_k\bigr)
\end{equation}
up to constants that neither depend on $\bm{x}_l$ nor $\boldsymbol{\theta}_l$. 
As mentioned in the \Cref{sec:introduction}, a prevalent algorithm used to approximate the minimizer of $\mathcal{G}$, and equivalently the MAP estimate $(\bm{x}^{\MAP}_l,\boldsymbol{\theta}^{\MAP}_l )$ in \cref{eq:model1_MAP_estimate}, is the GSBL framework \cite{lanza2020residual,tipping2001sparse,murphy2007conjugate}.
The GSBL algorithm uses a block-coordinate descent method \cite{wright2015coordinate,beck2017first} 
to compute the minimizer of the objective function $\mathcal{G}$ by alternating the  (i) minimization of $\mathcal{G}$ w.r.t.\ $\bm{x}_l$ for fixed $\boldsymbol{\theta}_l$, (ii) the minimization $\mathcal{G}$ w.r.t.\ $\boldsymbol{\theta}_l$ for fixed $\bm{x}_l$. 
Specifically, for an initial guess of hyperparameter vector $\boldsymbol{\theta}_l$ , the GSBL algorithm proceeds through a sequence of updates of the form 
\begin{subequations}\label{eq:GSBL}
\begin{equation}\label{eq:x_updateGSBL}
	\bm{x}^{\MAP}_l = \argmin_{\bm{x}_l} \left\{ \mathcal{G}(\bm{x}_l,\boldsymbol{\theta}_l ) \right\}, \end{equation}
    \begin{equation}\label{eq:update_beta1}
	\boldsymbol{\theta}^{\MAP}_l = \argmin_{\boldsymbol{\theta}_l} \left\{ \mathcal{G}(\bm{x}_l,\boldsymbol{\theta}_l ) \right\}, \end{equation}
\end{subequations}
until a convergence criterion is met.  \Cref{alg:GSBL} provides the basic GSBL method for estimating the underlying signals ${\bm x}_l$ from the model described in \cref{eq:GSBL}.

\begin{algorithm}[h!]
    \caption{GSBL algorithm of estimating $\bm{x}^{\MAP}$ }
    \label{alg:GSBL}
        \hspace*{\algorithmicindent} \textbf{Input:} Forward operator $F_{1:L}$, measurements $\bm y_{1:L}$, hyperprior parameters $(\beta_{1:L}, \vartheta_{1:L})$ in \eqref{eq:model1pdfhypertheta}, the prior operator $\Phi_{1:L}$ in \eqref{eq:model1}, and the noise precision $\alpha_{1:L}$ in \eqref{eq: iid_noise}.\\
        \hspace*{\algorithmicindent} \textbf{Output:} Signal estimate $\boldsymbol{x}^{\MAP}_{1:L}$.
    \begin{algorithmic}[1]
    \State{For $l=1,\dots,L$}
    \Indent
        \State{Initialize $\boldsymbol\theta_{1:L} = \boldsymbol{1}$,  and $\boldsymbol{x}_{1:L}=\boldsymbol{0}$.}
        \Repeat
            \State{Update $\boldsymbol{x}^{\MAP}_{1:L}$ by \eqref{eq:x_update}}
            \State{Update $\boldsymbol\theta^{\MAP}_{1:L}$ by \eqref{eq: update_theta}}
        \Until{convergence or maximum number of iterations is reached}. 
    \EndIndent
    \end{algorithmic}
\end{algorithm}

\begin{remark}\label{rem:hyperparameters}
Since the prior operator $\Phi$ is applied uniformly to all $l$ measurements, we have $K = K_1=\cdots=K_L$.
The individual recovery model (\cref{alg:GSBL}) is highly flexible, allowing for distinct hyperprior parameter vectors,  $\boldsymbol{\beta}_l = \{(\beta)_k\}_{k= 1}^K$ and $\boldsymbol{\vartheta}_l = \{(\vartheta_l)_k\}_{k= 1}^K$ in \cref{eq:model1hypertheta}, for each measurement $\bm y_l$. This structure allows regularization to be tailored independently for each measurement and even for each component of the prior operator. Consequently,  \cref{alg:GSBL} computes a separate hyperprior  vector $\boldsymbol{\theta}_l$ independently for each dataset, guided by uniform hyperprior parameter inputs.    In contrast, the joint recovery obtained by \cref{alg:MMV_GSBL} is inherently designed to enforce a common structure. It achieves this by using a single, shared pair of  hyperprior parameter vectors, $\boldsymbol{\beta}=\{\beta_k\}_{k=1}^K$ and $\boldsymbol{\vartheta}=\{\vartheta_k\}_{k=1}^K$ to compute one shared hyperprior vector $\boldsymbol{\theta}$ that couples all $L$ measurements. For the experiments in this paper, we simplify the hyperprior parameter vectors in two ways: (1) we constrain these vectors to have constant entries, i.e. $(\beta_l)_k=\beta_l$ and $(\vartheta_l)_k = \vartheta_l$ for $l=1,\dots,L$; and (2) we set these scalars to be identical across all $L$ measurements, i.e.~$(\beta_l)_k = \beta$ and $(\vartheta_l)_k = \vartheta$, for all $k$ and $l$, respectively.\footnote{The second simplification applies only to \Cref{alg:GSBL}, since by construction joint sparsity enforces this constraint in \Cref{alg:MMV_GSBL}. The first  simplification applies to both algorithms.} 
More details will be provided in \cref{sec:jointsparse}.

\end{remark}

We now detail the specific update steps of \Cref{alg:GSBL}, which are used to compute the MAP estimate, $(\bm{x}^{\MAP}_l,\boldsymbol{\theta}_l^{\MAP})$, as defined in \eqref{eq:model1_MAP_estimate}.

\subsubsection*{Updating the parameter vectors $\bm{x}_{1:L}$} 
Updating each sequential $\bm{x}_l$, $l = 1,\dots, L$, given corresponding $\boldsymbol{\theta}_l$ in \cref{eq:x_updateGSBL} (step 4 in \Cref{alg:GSBL}) reduces to solving the quadratic optimization problem 
\begin{equation}\label{eq:x_update} 
	\bm{x}^{\MAP}_{l} 
		= \argmin_{\bm{x}_l} \left\{ \alpha_l \| F_l \bm{x}_l - \bm{y}_{l} \|_2^2 + \| D_{\boldsymbol{\theta}_{l}}^{1/2} \Phi_l \bm{x}_l \|_2^2 \right\},
\end{equation}
where $D_{\boldsymbol{\theta}_{l}} = \diag(\boldsymbol{\theta}_{l})$. 

There is a unique solution to \cref{eq:x_update} for each $l=1,\dots,L$ as long as the \emph{common kernel condition} 
\begin{equation}\label{eq:common_kernel}
    \kernel(F_l) \cap \kernel(\Phi_l) = \{ \mathbf{0} \}
\end{equation} 
is satisfied. This assumption is widely adopted in the context of regularized inverse problems \cite{tikhonov1995numerical, kaipio2006statistical}.  Here $\kernel(G) := \{ \, {\bm t} \in \R^n \mid G {\bm t} = \mathbf{0} \, \}$  is the kernel of an operator $G: \R^n \to \R^m$, i.e., the set of vectors that are mapped to zero by $G$.
Due to its exclusive use of the $2-$norm, methods such as the preconditioned conjugate gradient (PCG) method \cite{saad2003iterative} (potentially combined with an early stopping based on Morozov's discrepancy principle \cite{calvetti2015hierarchical,calvetti2018bayes,calvetti2020sparse}) and the gradient descent approach \cite{glaubitz2022generalized} can all be used to efficiently solve \cref{eq:x_update}.  
As there is no general advantage of one method over  the other, which method to use should be made based on the specific problem (and the structure of $F_l$ and $\Phi_l$) at hand.  

\subsubsection*{Updating the hyperparameter vectors $\boldsymbol{\theta}_{1:L}$} 
As prescribed in \cref{eq:model1hypertheta}, for each $l = 1,\dots, L$, the hyperparameter vector $\boldsymbol{\theta}_l$ is modeled using the  Gamma distribution ($\Gamma(\beta_l,\vartheta_l)$). 
To update the hyperparameter vector $\boldsymbol{\theta}_l$ in \cref{eq:update_beta1} (step 5 in \Cref{alg:GSBL}) for fixed  $\bm{x}_l$, 
we substitute \cref{eq:G} into \cref{eq:update_beta1} and ignore all terms that do not depend on $\boldsymbol{\theta}_l$, yielding 
\begin{equation}\label{eq:update_beta2} 
	(\theta_{l})_k^{\MAP} = \argmin_{(\theta_{l})_k} \left\{  \tfrac{(\theta_{l})_k[\Phi_l {\bm x}_l]_k^2}{2} + \tfrac{(\theta_l)_k}{\vartheta_l}+  \eta\log\!\bigl((\theta_l)_k\bigr) \right\}, \quad k = 1,\dots, K,
\end{equation} 
where $\eta = \beta_l - 1/2$ and $[\Phi_l {\bm x}_l]_k$ denote the $k$-th entry of the vector $\Phi_l {\bm x}_l \in \R^K$. 
Differentiating the objective function in \cref{eq:update_beta2} w.r.t.\ $(\theta_l)_k$ and setting this derivative to zero gives
\begin{equation}\label{eq:update_beta3} 
	0 =  [\Phi_l {\bm x}_l]_k^2/2 + \vartheta_l^{-1} - (\theta_l)_k^{-1} \eta.
\end{equation} 
The minimizer of $\mathcal{G}$ w.r.t. $\boldsymbol{\theta}_l$ for fixed $\bm x_l$ is
\begin{equation}
    \label{eq: update_theta}
    (\theta_{l})_k^{\MAP} = \tfrac{\eta}{[\Phi_l {\bm x}_l]_k^2/2 + \vartheta_l^{-1}}.
\end{equation}

\subsubsection{Exploiting joint sparsity from multiple measurements}
\label{sec:jointsparse}
In applications such as synthetic aperture radar (SAR), multiple measurements of the same static image are acquired via different polarization configurations \cite{jakowatz2012spotlight}. Analogously, parallel magnetic resonance imaging (pMRI) \cite{chun2017compressed,chun2015efficient} collects measurements of a static image from multiple viewing angles that are {\em jointly sparse}. Finally, the state variable solutions for numerical  systems of conservation laws will often have shocks in the same physical location \cite{LeVeque02}.  In each case, the acquired data may be referred to as {\em multiple measurement vectors} (MMVs).  Compressive sensing algorithms \cite{adcock2019joint,gelb2019reducing,sanders2017composite} and sparse Bayesian learning methods \cite{glaubitz2024leveraging,green2023leveraging,green2024complexmmv,xiao2023sequential,zhang2022empirical} have been effectively developed to exploit mutual information and the assumption of joint sparsity in the underlying signals and images. In particular, the algorithm developed in \cite{glaubitz2024leveraging} codifies the presumed joint sparsity of piecewise constant images in the construction of the hyperprior parameter vector ${\boldsymbol{\vartheta}}$ in \cref{eq:model1hypertheta}, which is then to be seamlessly integrated into \cref{alg:GSBL}, resulting in \Cref{alg:MMV_GSBL} \cite{glaubitz2024leveraging} (for general $\Phi$). Importantly, no additional information is needed regarding  hyperparameter or variance estimations. The approach is described below. 

We begin by specifying  the likelihood, prior, and hyperprior densities hierarchical sparse Bayesian learning framework as
\begin{subequations}\label{eq:model2_PDFs} 
\begin{equation}
	\pi_{\mathbf{Y}_{1:L}|\mathbf{X}_{1:L}} ( \bm{y}_{1:L} | \bm{x}_{1:L}) 
		 \propto \exp\left( - \frac{1}{2} \sum_{l=1}^L\alpha_l\| F_l \bm{x}_l - \bm{y}_l \|_2^2 \right), \label{eq:model2pdflikelihood}\end{equation} 
\begin{equation}	\pi_{\mathbf{X}_{1:L}|\boldsymbol\Theta} ( \bm{x}_{1:L} | \boldsymbol{\theta}_\MMV ) 
		 \propto \det( D_{\boldsymbol{\theta}_\MMV} )^{L/2} \exp\left( - \frac{1}{2} \| D_{\boldsymbol{\theta}_\MMV}^{1/2} \Phi_l \bm{x}_l \|_2^2 \right), \label{eq:model2pdfprior}\end{equation}
\begin{equation}	\pi_{\boldsymbol{\Theta}} (\boldsymbol{\theta}_\MMV)
=\prod_{k=1}^K\Gamma((\theta_\MMV)_k|\beta,\vartheta)
		 \propto  \det(D_{\theta_\MMV})^{\beta - 1} \exp\left( - \sum_{k=1}^K (\theta_\MMV)_k/\vartheta  \right).\label{eq:model2pdfhypertheta}\end{equation}
\end{subequations} 
The joint sparsity is characterized by the hyperprior vector $\boldsymbol\theta_\MMV$.  In particular, in contrast to updating hyperparameter vectors $\boldsymbol\theta_{1:L}$ {\em separately}, as is done in \Cref{alg:GSBL}, here each sequential MAP estimate for $\bm{x}_{1:L}$ depends upon the MAP estimate of $\boldsymbol\theta_\MMV$. 

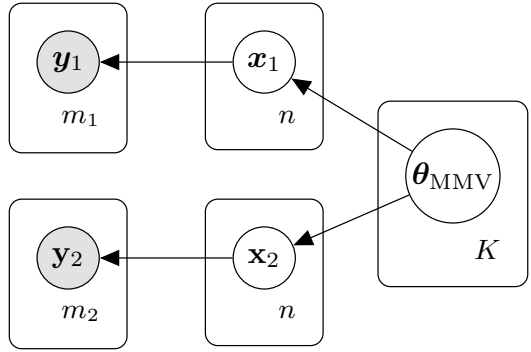
\begin{figure}[h!]
\centering
\resizebox{0.4\textwidth}{!}{%
\begin{tikzpicture}
    \node[obs] (y1) {$\bm{y}_1$}; %
\node[latent, right=1.5 of y1] (x1) {$\bm{x}_1$} ; %
\node[const, right=1.75 of x1] (theta_aux) {}; %
\node[latent, below=.725 of theta_aux] (theta) {$\boldsymbol{\theta}_\MMV$} ; %
\node[obs, below=1.5 of y1] (y2) {$\mathbf{y}_2$}; %
\node[latent, right=1.5 of y2] (x2) {$\mathbf{x}_2$} ; %
\edge {x1} {y1}; %
\edge {x2} {y2}; %
\edge {theta} {x1}; %
\edge {theta} {x2}; %
\plate[inner sep=0.3cm, xshift=0cm, yshift=0cm] {plate_y1} {(y1)} {$m_1$}; %
\plate[inner sep=0.3cm, xshift=0cm, yshift=0cm] {plate_x} {(x1)} {$n$}; %
\plate[inner sep=0.3cm, xshift=0cm, yshift=0cm] {plate_theta} {(theta)} {$K$}; %
\plate[inner sep=0.3cm, xshift=0cm, yshift=0cm] {plate_y2} {(y2)} {$m_2$}; %
\plate[inner sep=0.3cm, xshift=0cm, yshift=0cm] {plate_x} {(x2)} {$n$}; %
\end{tikzpicture} 
}%
\caption{
    Graphical representation of \Cref{alg:MMV_GSBL} for two ($L=2$) measurement and parameters vectors, $\bm{y}_1, \bm{y}_2$ and $\bm{x}_1, \bm{x}_2$, respectively. 
  Shaded (plain) circles indicate observed (hidden) random variables, respectively. 
The arrows represent the conditional dependencies among the random variables. 
The likelihood \cref{eq:model2pdflikelihood} bridges the parameter vectors $\bm{x}_1,\bm{x}_2$ and the measurement vectors $\bm{y}_1,\bm{y}_2$; 
The joint-sparsity-promoting prior \cref{eq:model2pdfhypertheta} links the shared hyperparameter vector $\boldsymbol{\theta}_\MMV$ to $\bm{x}_1, \bm{x}_2$. 
Using shared Gamma distributed hyperparameter vector $\boldsymbol{\theta}_\MMV$ (instead of separate ones for $\bm{x}_1, \bm{x}_2$) results in $\Phi \bm{x}_1$ and $\Phi \bm{x}_2$ sharing the support. 
}
\label{fig:graphical_model}
\end{figure}

We employ \cref{eq:model2_PDFs} to reformulate the objective function \cref{eq:G} as \cite{glaubitz2024leveraging}
\begin{align}
    \label{eq:G_group}
    \begin{split}
    \mathcal{G}(\bm x_{1:L},\boldsymbol{\theta}_\MMV)=&\frac{1}{2}\left( \sum_{l=1}^L\alpha_l \norm{F_l \bm{x}_l - \bm{y}_l}_2^2 + \norm{D_{\boldsymbol{\theta}_\MMV}^{1/2}\Phi_l \bm{x}_l}_2^2\right) \\
    &+ \sum_{k=1}^K \frac{(\theta_\MMV)_k}{\vartheta} +  \sum_{k=1}^K \Bigl(\tfrac{2\beta - 2 + L}{2}\Bigr)\log\!\bigl((\theta_\MMV)_k\bigr).
    \end{split}
\end{align}
The MAP estimate is found by minimizing the corresponding joint objective function in \eqref{eq:G_group}, specifically given by 
\begin{align}
    \label{eq:GSBL_group}
    \begin{split}
    \bm{x}^{\MAP}_{1:L} =& \argmin_{\bm{x}_{1:L}} \left\{ \mathcal{G}(\bm{x}_{1:L},\boldsymbol{\theta}_\MMV ) \right\},\\
    \boldsymbol{\theta}^{\MAP}_\MMV =& \argmin_{\boldsymbol{\theta}_\MMV} \left\{ \mathcal{G}(\bm{x}_{1:L},\boldsymbol{\theta}_\MMV ) \right\},
    \end{split}
\end{align}
until a convergence criterion is met.  Notably, the update rules governing $\bm x_l$ for $l=1,\dots,L$ remain unchanged, following \eqref{eq:x_update}. Finally, as already mentioned, although the model supports vector-valued hyperprior parameters, for simplicity we use constant scalar values $\beta_k=\beta$ and $\vartheta_k=\vartheta$ in our numerical experiments.  \Cref{alg:MMV_GSBL} provides the basic steps for the GSBL procedure that exploits the joint sparsity across the $L$ measurements, while \cref{fig:graphical_model} (Figure 2 in  \cite{glaubitz2024leveraging}) provides a graphical depiction of the process.

\begin{algorithm}[h!]
    \caption{MMV-GSBL algorithm of estimating $\{\boldsymbol {x}^{\MAP}\}_{l = 1}^L$}
    \label{alg:MMV_GSBL}
        \hspace*{\algorithmicindent} \textbf{Input:} Measurements $\boldsymbol y_{1:L}$, hyperprior parameters $(\beta, \vartheta)$, the prior operator $\Phi_{1:L}$ in \eqref{eq:model1}, and the noise precision $\{\alpha_l\}_{l=1}^L$ in \eqref{eq: iid_noise}.\\
        \hspace*{\algorithmicindent} \textbf{Output:} Signal estimate $\boldsymbol{x}^{\MAP}_{1:L}$.
    \begin{algorithmic}[1]
    \State{Initialize $\boldsymbol\theta_{\MMV} = \boldsymbol{1}$,  and $\boldsymbol{x}_{1:L}=\boldsymbol{0}$.}
    \Repeat
        \State{Update $\boldsymbol{x}^{\MAP}_{1:L}$ by \eqref{eq:x_update}}
        \State{Update $\boldsymbol\theta^{\MAP}_\MMV$ by \eqref{eq: update_theta_group}}
    \Until{convergence or maximum number of iterations is reached }. 
    \end{algorithmic}
\end{algorithm}

\subsubsection*{Updating the {\em group} hyperprior $\boldsymbol{\theta}_\MMV$}
For each $l = 1,\dots,L$, the standard GSBL algorithm (\cref{alg:GSBL}) prescribes as one of its inputs a  constant vector  ${\boldsymbol{\vartheta}}_l$ as the hyperprior parameter vectors needed for computing the hyperprior vectors $\boldsymbol\theta_l$ in \cref{eq:model1hypertheta}. While using the gamma distribution  helps to promote the sparsity of $\Phi_l{\bm x_l}$, the {\em support locations} in the sparse domain are typically unknown.  The idea in \cite{glaubitz2024leveraging} is to leverage the joint sparsity which can be extracted from the multiple measurements to further refine how this is done. 
Substituting \eqref{eq:G_group} into $\eqref{eq:GSBL_group}$ and removing terms irrelevant to $\boldsymbol\theta_\MMV$ is equivalent to computing
\begin{equation}\label{eq:theta_group}
    \left(\theta^{\MAP}_{\MMV}\right)_k = \argmin_{(\theta_\MMV)_{k}} \left\{ \tfrac{(\theta_\MMV)_{k} \sum_{l=1}^L[\Phi_l {\bm x}_l]_k^2}{2} + \tfrac{(\theta_\MMV)_{k}}{\vartheta} + \eta \log( (\theta_\MMV)_{k} ) \right\}, \quad k = 1,\dots, K,
\end{equation}
where $\eta=\beta-1+L/2$.
The joint sparsity is incorporated into \cref{eq:theta_group} via hyperprior parameters $\beta$ and $\vartheta$ shared across $l=1,\dots,L$, distinguishing it from \eqref{eq:update_beta2}.
Equating the derivative of the objective function \eqref{eq:theta_group} with respect to $(\theta_\MMV)_k$	to zero results in a closed-form solution, analogous to \cref{eq: update_theta},
\begin{equation}
    \label{eq: update_theta_group}
    \left(\theta^{\MAP}_{\MMV}\right)_k = \tfrac{\beta-1+L/2}{\sum_{l=1}^L[\Phi_l {\bm x}_l]_k^2/2 + \vartheta^{-1}}.
\end{equation}

\subsubsection{Uncertainty Quantification}
While we compute the MAP estimate as a point solution for our reconstructed signal given MMV data $\bm y_{1:L}$, the hierarchical sparse Bayesian learning structure still allows us to partially quantify the uncertainty in the recovered parameters by analyzing the posterior distribution $\pi_{\mathbf{X}_{1:L},\boldsymbol\Theta|\mathbf{Y}_{1:L}=\bm y_{1:L}}$. Specifically, for fixed hyperparameter vector $\boldsymbol{\theta}_\MMV$, Bayes' theorem provides the fully conditional posterior for the parameter vectors 
$\mathbf{X}_{1:L}$ as
\begin{equation}
    \label{eq: bayes_UQ}
    \pi_{\mathbf{X}_{1:L}|\boldsymbol\Theta=\boldsymbol{\theta}_\MMV,\mathbf{Y}_{1:L} =\bm y_{1:L}}(\bm x_{1:L})
    \propto
    \pi_{\mathbf{Y}_{1:L}|\mathbf{X}_{1:L}}(\bm y_{1:L}|\bm x_{1:L})
    \pi_{\mathbf{X}_{1:L}|\boldsymbol{\Theta}}(\bm x_{1:L}|\boldsymbol{\theta}_\MMV).
\end{equation}
By substituting the likelihood density $\pi_{\mathbf{Y}_{1:L}|\mathbf{X}_{1:L}}(\bm y_{1:L}|\bm x_{1:L})$ from \eqref{eq:model2pdflikelihood} and the conditional prior density $\pi_{\mathbf{X}_{1:L}|\boldsymbol{\Theta}}(\bm x_{1:L}|\boldsymbol{\theta})$ from \eqref{eq:model2pdfprior} into \eqref{eq: bayes_UQ}, we obtain
\begin{equation}
    \label{eq: pos_UQ}
    \pi_{\mathbf{X}_{1:L}|\boldsymbol\Theta=\boldsymbol{\theta}_\MMV,\mathbf{Y}_{1:L} =\bm y_{1:L}}(\bm x_{1:L})\propto
    \frac{1}{2}\left( \sum_{l=1}^L\alpha_l \norm{F_l \bm{x}_l - \bm{y}_l}_2^2 + \norm{D_{\boldsymbol{\theta}_\MMV}^{1/2}\Phi_l \bm{x}_l}_2^2\right).
\end{equation}
For convenient notation, we define the posterior mean and covariance matrix as
\begin{equation}
\label{eq: pos_mean_var}
    \boldsymbol{\mu}_l=\alpha\Sigma_lF^T\bm y, 
    \quad
    \Sigma_l=(\alpha F^TF+\Phi_l^TD_{\theta_\MMV}^{-1} \Phi_l)^{-1}.
\end{equation}
With these definitions, the MAP estimate for each posterior distribution of $\mathbf {X}_l$ is proportional to a normal distribution, given by
\begin{equation*}
\pi_{\mathbf{X}_{l}|\boldsymbol\Theta=\boldsymbol{\theta},\mathbf{Y}_{l} =\bm y_{l}}(\bm x_{l})\propto\mathcal{N}\left(\boldsymbol{\mu}_l|\Sigma_l\right).
\end{equation*}
When the common kernel condition in equation \eqref{eq:common_kernel} holds, the covariance matrix $\Sigma_l$ is guaranteed to be symmetric and positive definite (SPD). This allows us to quantify the uncertainty in  $\mathbf{X}_l$ by sampling from the resulting normal distribution to determine the sample mean and credible intervals. For more details on sampling from high-dimensional Gaussian distributions, see \cite{mackay1992bayesian,glaubitz2025efficient}.

\begin{remark}
    The uncertainty analysis in this work is conditional.  It quantifies the variability in the solution $\bm x_{1:L}$ but relies on a point estimate for the hyperparameter vector $\boldsymbol{\theta}_\MMV$. As such, our approach does not capture the uncertainty inherent in the estimation of $\boldsymbol{\theta}_\MMV$. A more comprehensive uncertainty analysis would require sampling from the full joint posterior distribution, $p(\bm x_{1:L}, \boldsymbol{\theta}_\MMV|\bm y_{1:L})$. This is typically achieved using a sampling-based approach, such as a Markov chain Monte Carlo (MCMC) method, which would fully characterize all sources of uncertainty in the model \cite{owen2013monte}.
\end{remark}

\subsection{Edge detection}
\label{sec:edgeDetection}

Before integrating the {\em residual prior transform} into either \Cref{alg:GSBL} or \Cref{alg:MMV_GSBL}, we first review the necessary ingredients for its construction. To this end, we note that sparsity-promoting operators can be viewed as {\em edge detectors} for piecewise smooth functions. An effective edge detector should be able to identify true discontinuities from the smooth variations in the underlying signals. Such effectiveness hinges on how well the operator is matched to the variability of signal. Hence,  conventional sparsity-promoting operators typically require a priori knowledge of smoothness order. The residual prior transform overcomes this challenge by exploring the principle that operators designed to match the same type of variability will produce similar responses.   The core idea is to take two distinct edge detectors that are both matched to the same class of smooth functions. By design, the {\em residual} between the outputs of these two operators will be minimal in smooth regions while magnifying their disagreement at discontinuities. This creates a new operator that is inherently sparse in a known way, without requiring precise a priori knowledge of the signal's variability.

While a wide variety of operator pairings could be used to construct a residual prior transform, this study adopts the specific operators from \cite{xiao2025new}, namely a local differencing operator and a global Fourier-based concentration factor edge detector. This deliberate pairing is motivated by each having well-established convergence properties \cite{gelb2006adaptive}, since due to linearity, the convergence and robustness properties of the  residual transform operator can also be established. For clarity, a brief review of each operator follows.

\subsubsection{Jump function approximation}
\label{sec:Edgesgeneral}
Let ${\bm f}\in\mathbb{R}^n$ be the discrete representation of a piecewise smooth, $2\pi$-periodic function, $f(s)$. This vector is created by sampling $f(s)$ at $n$ uniformly-spaced grid points, $s_j = -\pi + (j-1)\Delta s = -\pi + \frac{2\pi (j-1)}{n}$, for $j = 1,\dots,n$. In other words, $f$ is simply the vector of sampled values:
\begin{equation}
\label{eq:discretization_f}
{\bm f} = \{ f(s_j) \}_{j=1}^n. 
\end{equation} 
The {\em jump function}, $[f]$, is defined as the difference between the right- and left-hand limits of a function $f$ at a specific point $s$:
\begin{equation} \label{eq:jumpfunction}
[f](s):=f(s^+)-f(s^-),
\end{equation}
This quantifies the size of a discontinuity, or ``jump'', in the function $f$ at that location.  
If a function $f$ contains $\mathcal{W}$ simple discontinuities at locations $\{\xi_w\}_{w=1}^{\mathcal{W}} \subset [-\pi,\pi)$, its jump function \eqref{eq:jumpfunction} can be equivalently represented by:
\begin{equation}
\label{eq:jumpfunction2}
[f](s) =\sum_{w=1}^{\mathcal W} f(\xi_w)I_{\xi_w}(s), 
\end{equation}
where the indicator function, $I_{\xi_w}(s)$, is  defined as
\[I_{\xi_w}(s) = \begin{cases}
1, & \text{if }  s=\xi_w \\
0, & \text{otherwise}
\end{cases}.\]
The goal of an edge detector is to approximate $[f](s_j)$ within each cell $[s_j,s_{j+1})$, $j = 1,\dots,n$, with $f(s_{n+1}) = f(s_1)$ due to periodicity.  The approximation of the jump function at each cell location, $[f](s)$, is what populates the entries of the {\em edge vector}. We denote the ground truth edge vector as ${\bm g} \in \mathbb{R}^n$, with 
\begin{equation}
    \label{eq:groundtruthedge}
    {\bm{g}}_j = [f](s_j), \quad j = 1,\dots,n.
\end{equation}

\subsubsection{A local approach to edge detection}
\label{sec: local_approach}
We employ a local differencing formulation to detect edges within the cell between two adjacent grid points, $s_j$ and $s_{j+1}$.  This method leverages the local behavior of the signal to identify sharp changes. If an edge, or jump discontinuity, is located at $\xi$ within the interval $s_j \leq \xi \leq s_{j+1}$ as defined in \eqref{eq:jumpfunction2}, we can describe this relationship using an asymptotic statement under certain smoothness conditions:
\begin{equation}
\label{eq:asym_stmt}
\Delta f_{j+\frac{1}{2}} := f_{j+1}-f_j=
\begin{cases}
[f](\xi) + \mathcal{O}(\Delta s), & \text{if } j=j_\xi : \xi\in[s_j, s_{j+1}) \\
\mathcal{O}(\Delta s), & \text{otherwise}
\end{cases},
\end{equation}
where $f_j$ is the discretized representation defined in \cref{eq:discretization_f} for $j = 1,\dots,n$ with $\Delta s = \frac{2\pi}{n}$. Note that \cref{eq:asym_stmt} provides the theoretical foundation for using TV or equivalent first-order differencing operators. Moreover, \cref{eq:asym_stmt} relies on the assumption that $\abs{[f](\xi_w)}\gg\mathcal{O}(\Delta s)$, $w = 1,\dots,\mathcal{W}$, which allows us to effectively distinguish true edges from the smooth variations in the function, preventing false positives.

First-order local differencing operators are \sout{is} optimally suited only for piecewise constant functions, and its performance degrades as the function exhibits greater variability.  Higher-order differencing operators can effectively detect edges in more complex signals, such as piecewise polynomial functions \cite{archibald2005polynomial}. The $2p+1$-order differencing operator is a direct extension of \eqref{eq:asym_stmt} and is given by the scaled version of the Newton divided difference formula
\begin{equation}
\label{eq:general_differencing}
\Delta^{2p+1} f_{j+\frac{1}{2}} 
= \sum_{l=0}^p (-1)^{l}\binom{2p+1}{p-l
} (f_{j+1+l} - f_{j-l}), \quad j = 1,\dots, n.
\end{equation}
Since $f$ is periodic, the formula simply replaces  the necessary components $f_{n+1+p}$ and $f_{1-p}$ with their equivalent values $f_{1+p}$ and $f_{n-p+1}$ respectively.
Its convergence to the jump function is characterized by the asymptotic statement  \cite{gelb2000detection,gelb2002spectral,gelb2006adaptive} 
\begin{equation}
\label{eq:differencing_jump}
\Delta^{2p+1} f_{j+l+\frac{1}{2}}=
\begin{cases}
(-1)^{l} q_{l,p}[f](\xi) + \mathcal{O}(\Delta s), & \abs{l}\leq p \\
\mathcal{O}(\Delta s)^{2p+1}, & \abs{l} > p
\end{cases},
\end{equation}
where
\[
q_{l,p} = (-1)^{l}\sum_{k=0}^{p-\abs{l}}\binom{2p+1}{k}(-1)^k = \binom{2p}{p+\abs{l}}.
\]
Specifically, \eqref{eq:differencing_jump} depicts how \eqref{eq:general_differencing} behaves both within and outside of the $2p$-cell neighborhood of each jump discontinuity. Consequently, the $2p+1$-order local edge detection method for discretized vector $\bm f$ over $n$ uniformly-spaced
grid points is defined as \cite{archibald2005polynomial,gelb2006adaptive}
\begin{equation}
\label{eq:local_edge_detector}
T_n^{\Delta_{2p+1}}{f}_j = \frac{1}{q_{0,p}}\Delta^{2p+1} f_{j+\frac{1}{2}},  \quad j = 1,\dots, n,
\end{equation}
where 
\begin{equation}\label{eq:q0P}
 q_{0,p}=\binom{2p}{p}.
\end{equation}

For the purpose of constructing the residual prior transform, we express the local edge detection method from  \eqref{eq:local_edge_detector} in a matrix-vector form, $T_{n}^p {\bm f} \approx {\bm g}$, with ${\bm g}$ defined in \cref{eq:groundtruthedge}. This formulation is defined by the matrix entries $T_{n}^p(j,j')$ \cite{archibald2016image}
\begin{equation}
    \label{eq:localedgeMV}
   \underset{j,\, j' = 1, \dots, n}
    {T_n^p(j, j')} = \begin{cases}
        \tfrac{(-1)^{j'+1}}{q_{0,p}}\binom{2p+1}{p-j'}, & j-p\leq j'\leq j\\
        \tfrac{(-1)^{j'}}{q_{0,p}}\binom{2p+1}{p-j'}, & j+1\leq j'\leq j+p+1\\
        0, & \text{otherwise}
    \end{cases},
\end{equation} 
yielding $(T_n^p\bm f)_j=T_n^{\Delta_{2p+1}}f_j$. The convergence  of $T^p_n{\bm f}$ to ${\bm g}$ is ensured by \eqref{eq:differencing_jump}. In particular, $p = 0$ reduces the operator to standard first-order differencing.
Even-order operators for the formulation in \eqref{eq:local_edge_detector} are omitted from this discussion due to their requirement for non-symmetric stencils. Nevertheless, such differencing operators can be generated in a similar fashion and readily incorporated into these constructions.

\subsubsection{A global approach to edge detection}\label{sec:globaledge}
A {\em globally}-based edge detector for the discretized vector $\bm f$ from \eqref{eq:discretization_f} is implemented using the discrete Fourier transform (DFT). This method identifies edges by analyzing the overall frequency content of the signal rather than by examining local neighborhoods.
The DFT matrix, denoted by $F$, is defined by its entries 
\begin{equation}
    \label{eq: DFT}
    \hat{F}_{j,k}=\tfrac{1}{n}e^{-iks_j}, \quad j = 1,\dots,n, \quad k = -\frac{n}{2}, \dots, \frac{n}{2}-1.
\end{equation}
The discrete (pseudo-spectral) Fourier coefficient vector $\hat{\bm f}\in\C^n$,  is then computed using this matrix, where each component is calculated as
\begin{equation}\label{eq: fourcoefdiscrete}
    \hat{f}_k = \frac{1}{n} \sum_{j=1}^n f(s_j) e^{-i k s_j}.
\end{equation}
This entire transformation can be expressed compactly in a matrix-vector product
\begin{equation}
    \label{eq:discretefourMV}
\hat{\bm f} = \hat{F}{\bm f}. 
\end{equation}
We implement the {\em concentration factor edge detection} method \cite{gelb1999detection,gelb2006adaptive} to  detect edges from the discrete Fourier data, given by 
\begin{equation}
\label{eq:consum}
S^{\sigma}_n f(s)
:=  \pi i \sum_{k=-\frac{n}{2}}^{\frac{n}{2}} \operatorname{sgn}(k) \, \sigma\left(\frac{|k|\Delta s}{\pi}\right)sinc\left(\frac{|k|\Delta s}{2}\right)\hat{f}_k \, e^{i k s} \approx [f](s),
\end{equation}
where  $[f](s)$ is defined in \cref{eq:jumpfunction2}. This technique relies on {\em an admissible} concentration factor function $\sigma(\cdot)$ satisfying the requirements

\begin{enumerate}
    \item  $\displaystyle K^\sigma_n(s) = \sum_{k=1}^{\frac{n}{2}} \sigma{\left(\frac{2k}{n}\right)} \sin(ks)$ must be an odd function.

    \item $\frac{\sigma(\eta)}{\eta}\in C^2(0,1)$, ensuring the continuity of both the first and second derivatives.

    \item $\int^1_\epsilon \frac{\sigma(\eta)}{\eta} \, d\eta \to 1$, $\epsilon = \epsilon(n) > 0$ being small, which guarantees the appropriate normalization in the limit.
\end{enumerate}
The concentration factor $\sigma(\cdot)$ acts as a scaled band-pass filter with tuning parameters determining crucial properties such as the order of convergence, along with the degree of localization and smoothing.  By modulating the Fourier coefficients, \cref{eq:consum} converges to the jump function $[f](s)$ in \cref{eq:jumpfunction} for any $s \in [-\pi,\pi)$.  While the {\em discrete} Fourier coefficients in \eqref{eq: fourcoefdiscrete} are used to generate \eqref{eq:consum}, the admissible concentration factor, $\sigma(\cdot)$, was originally derived using {\em continuous} Fourier coefficients, i.e. $\hat f_k=\tfrac{1}{2}\int^\pi_{-\pi}f(s)e^{-iks}ds$, through integration by parts. In a similar vein, the discrete case requires summation by parts which introduces the additional sinc factor in \cref{eq:consum}, a detail discussed further in \cite{gelb1999detection}. 

We apply \eqref{eq:consum} at $s=s_{j+\zeta}$  for $j=1,\dots,n$ and $\zeta\in[0,1)$ to identify a jump discontinuity within the grid cells $[s_j,s_{j+1}]$.  By substituting the matrix-vector formulation from \eqref{eq:discretefourMV} into \eqref{eq:consum}, we can express this approximation as a matrix-vector product, $ S_{n,\zeta}^\sigma {\bm f}\approx\bm g$, the true edge vector in \cref{eq:groundtruthedge}. The entries of the matrix $S_{n,\zeta}^\sigma$ are then given by 
\begin{equation}\label{eq:matrix_S_mu}
    S_{n,\zeta}^{\sigma}(j,j')=\sum_{k=1}^{\tfrac{n}{2}}\frac{\sigma\left(\frac{k\Delta s}{\pi}\right)}{k}\left(\cos{(k\Delta s\varphi^+_{j,j'}}-\cos{(k\Delta s\varphi^-_{j,j'})}\right),\quad j,j' = 1,\dots, n,
\end{equation}
where $\varphi^+_{j,j'} = j -j' +(\tfrac{1}{2}+ \zeta)$ and $\varphi^-_{j,j'} = j-j' -(\tfrac{1}{2} - \zeta)$. 

A variety of admissible concentration factors $\sigma(\cdot)$ have been proposed based on their associated localization and convergence characteristics. For the {\em residual prior transform} detailed in \Cref{sec:residualprior}, we select 
\begin{equation}
\label{eq: trig_conc}
    \sigma(\eta) := \sigma_{2p+1}(\eta) = \frac{2^{2p}\eta \sin^{2p}(\frac{\pi}{2}\eta)}{q_{0,p}}.
\end{equation}
since when evaluated at $s_{j+\tfrac{1}{2}}$, {\em exactly} matches the differencing behavior of the local edge detection method given by \eqref{eq:localedgeMV}.  That is, $S^{\sigma_{2p+1}}_nf(s_{j+\frac{1}{2}}) = T_n^{\Delta_{2p+1}}f_j$, $j = 1,\dots,n$,  with  $q_{0,p}$ defined in \cref{eq:q0P}. Hence, substituting \eqref{eq: trig_conc} into \eqref{eq:matrix_S_mu} 
provides
\begin{equation}
    \label{eq:matrixconcsum}
S^{\sigma_{2p+1}}_{n,{\zeta}}(j, j') = \frac{2^{2p+1}}{n q_{0,p}}\sum_{k = 1}^{\frac{n}{2}} \sin^{2p}{\left(\frac{k\Delta s}{2}\right)}\left(\cos{(k\Delta s\varphi^+)}- \cos{(k\Delta s\varphi^-)}\right),
 \end{equation}
for $j,j'= 1,\dots, n$,
yields $(S^{\sigma_{2p+1}}_{n,{\zeta}}\bm f)_j=S^{\sigma_{2p+1}}_{n}f(s_{j+\zeta})$ and in general
\begin{equation}
    T_{n}^p{\bm f} = S^{\sigma_{2p+1}}_{n,\frac{1}{2}}{\bm f}.
\label{eq:equivalencies}
\end{equation}
Both $T_n^p{\bm f}$  and $S^{\sigma_{2p+1}}_{n,\frac{1}{2}}{\bm f}$  aim to approximate the edge vector $\bm g$ for a piecewise smooth function $f$ defined on  $[-\pi, \pi)$. The precision of this approximation is contingent upon three factors: (1) the underlying variability of $f(s)$; (2) the refinement of the discretization step $\Delta s$; and (3) the appropriateness of the selected order  $2p+1$. More detailed discussions of the equivalency between the local differencing and pseudo-spectral band pass filtering approaches, as given by \eqref{eq:equivalencies}, are available in \cite{gelb2006adaptive, xiao2025new}.
 
\section{The residual prior transform}
\label{sec:residualprior}

With the necessary components for its construction established, we now formalize the definition of the residual prior transform and provide the error analysis of its constituent parts. 
To demonstrate the validity of the MAP estimate in equation \eqref{eq:x_MAP}, we begin with the ordinary least squares estimate of $\bm x$ from \eqref{eq: data_acquisition} given by
\begin{eqnarray}\label{eq:x_MAP}
    \bm x^{EST} = (F^TF)^\dagger F^T\bm y= (F^TF)^\dagger F^T(F\bm x+\boldsymbol\epsilon)=\bm f+ (F^TF)^\dagger F^T\boldsymbol{\epsilon} = {\bm f} + \boldsymbol{\epsilon}^{EST},
\end{eqnarray}
with ${\bm f}$ being the true signal discretized by \cref{eq:discretization_f} and $\boldsymbol{\epsilon}^{EST}$ being the corresponding Gaussian distributed error. The symbol $\dagger$ denotes the Moore-Penrose pseudo-inverse of a matrix. For notational simplicity we have omitted the subscript $l$. The error mean follows as
\begin{eqnarray}\label{eq:expected_eps}
    E\left[\boldsymbol{\epsilon}^{EST}\right] = E\left[(F^TF)^{\dagger}F^T\boldsymbol{\epsilon}\right] = (F^TF)^{\dagger}F^TE[\boldsymbol{\epsilon}]=0.
\end{eqnarray}
Moreover, the variance of $\boldsymbol{\epsilon}^{EST}$ can be computed directly from $Cov(\boldsymbol{\epsilon})=\alpha^{-1}I$ as
\begin{align}
Cov\left[\boldsymbol{\epsilon}^{EST}\right]&= Cov\left[(F^TF)^{\dagger}F^T\boldsymbol{\epsilon}\right]=(F^TF)^{\dagger}F^T Cov(\boldsymbol{\epsilon}) \left((F^TF)^{\dagger}F^T\right)^T\nonumber\\
    &=\alpha^{-1}(F^TF)^{\dagger}F^T F\left((F^TF)^{\dagger}\right)^T = \alpha^{-1}(F^TF)^{\dagger}.\label{eq:var_eps}
\end{align}
Finally, we observe from \cref{eq:equivalencies}  that
\[|T_n^p{\bm x}^{EST} - S_{n,\zeta}^{\sigma_{2p+1}}{\bm x}^{EST}| \approx |T_n^p{\boldsymbol{\epsilon}}^{EST} - S_{n,\zeta}^{\sigma_{2p+1}}{\boldsymbol{\epsilon}}^{EST}| = |(T_n^p - S_{n,\zeta}^{\sigma_{2p+1}}){\boldsymbol{\epsilon}}^{EST}|.\]
The above analysis motivates the adaptation of the residual transform operator from \cite{xiao2025new} to use in what have coined the residual prior transform in our hierarchical sparse Bayesian learning framework.  Specifically, 

\begin{definition} The residual prior transform is defined as 
    \label{def:residtransform}

\begin{equation}
    \label{eq:residualprior}
    R_{n,\zeta}^p := T_{n}^p-S_{n,\zeta}^{\sigma_{2p+1}},
\end{equation}
where $n$ is the vector length for discretized signal ${\bm f}$, $2p+1$ is the order of the convergence in smooth regions of the corresponding jump function $[f](s)$ in \cref{eq:jumpfunction}, and $\zeta \in [0,1)$ corresponds to the the point $s_{j+\zeta}$ within each grid cell $[s_j,s_{j+1}]$ on which the jump function is recovered using \cref{eq:consum}.
\end{definition}

As already noted and will be described in more detail in \cref{sec:noiseresidual}, for a general piecewise smooth function $f(s)$ on $[-\pi,\pi]$, the residual prior transform yields a more accurate description of the underlying discretized signal ${\bm f}$ in \cref{eq:discretization_f}. This is because the variability in the smooth regions is better matched by the residual prior transform. 
By contrast, effective use of both the standard local differencing  $T_n^p$ and the band pass filtering  $S_{n,\zeta}^{\sigma_{2p+1}}$ transforms rely on  prior knowledge of {\em fixed} variability for piecewise smooth functions.  {\em By design} they fail to capture the inherent variable smoothness in real world signals, which inevitably may lead to a significant mischaracterization of the underlying signal. Moreover, because the genuine variability is rarely known a priori,  conventional priors are subject to the rigid assumption of a fixed order of $2p+1$ (often $p = 0$ corresponding to a sparse gradient domain), which further compromises their effectiveness in capturing the complex structural information in the underlying scene. The residual prior transform offers a promising alternative since it is robust to highly variable behaviors in smooth regions of signals. The following error analysis, corroborated by the numerical experiments in \Cref{sec:numerics}, rigorously establishes that $R_{n,\zeta}^p$ achieves a higher reconstruction fidelity than $T_n^p$.
\begin{remark}\label{rem: onlyp0}
    All of our numerical experiments use $p=0$, as higher-order operators are redundant within our residual framework. The core idea is that the residual prior transform, $\Phi=\Phi_1-\Phi_2$, is designed such that the responses of $\Phi_1$ and $\Phi_2$ are similar over smooth regions of the signal. Consequently, any information captured by higher-order components, which primarily describe these smooth regions, would be largely canceled in the subtraction, having a negligible effect on the final result.\footnote{Numerical experiments in \cite{xiao2025new} use the residual transform operator $R_{n,\zeta}^p$ in the context of LASSO regression for $p = 1$.}
\end{remark}

\subsection{MAP estimate error analysis}\label{sec:noiseresidual} 
We now describe how the choice of each sparsity-promoting prior transform affects the performance of the MAP recovery in \cref{eq:model1_MAP_estimate}. This analysis is grounded in the established equivalency between the two transform matrices, as shown in \eqref{eq:equivalencies}.

The  expectation  $E(T_n^p{\bm x}^{EST})$  is calculated componentwise as
\begin{eqnarray*}
    E\left[(T_n^p{\bm x}^{EST})_j\right]
    &=&E\left\{ \tfrac{1}{q_{0,p}}\sum_{l=0}^p (-1)^{l}\binom{2p+1}{p-l
} (f_{j+1+l} + \epsilon^{EST}_{j+1+l} - (f_{j-l}+\epsilon^{EST}_{j+l}))\right\}\\
&=&E\left[(T_n^p {\bm f})_j\right]+E\left[(T_n^p\boldsymbol{\epsilon}^{EST})_j\right]
= (T_n^p{\bm f})_j,
\end{eqnarray*}
that is, $T_n^p$ as defined by \cref{eq:localedgeMV} is an unbiased estimator. By contrast,  the  $Var(T_n^p{\bm x}^{EST})$  depends on the noise precision $\alpha$ defined in \cref{eq: iid_noise}, the order of the difference operator $2p+1$, and the stencil location.  Specifically, incorporating $Cov(\boldsymbol{\epsilon}^{EST})$ we have 
\begin{eqnarray*}
Var\left[(T_n^p\bm{x}^{EST})_j\right]&=&E\left\{\left( \tfrac{1}{q_{0,p}}\sum_{l=0}^p (-1)^{l}\binom{2p+1}{p-l
} (f_{j+1+l} + \epsilon^{EST}_{j+1+l} - (f_{j-l}+\epsilon^{EST}_{j-l}))-E\left[(T_n^p\bm{x}^{EST})_j\right]\right)^2\right\}\\
&=&E\left\{\left( \tfrac{1}{q_{0,p}}\sum_{l=0}^p (-1)^{l}\binom{2p+1}{p-l} \left( \epsilon^{EST}_{j+1+l}  -\epsilon^{EST}_{j-l}\right)\right)^2\right\}.
\end{eqnarray*}
Based on $E[(\cdot)^2]=Var[\cdot]+(E[\cdot])^2$, \eqref{eq:expected_eps} allows for further simplification, so that
\begin{align}
Var\left[(T_n^p\bm{x}^{EST})_j\right]
    &= Var\left[ \tfrac{1}{q_{0,p}}\sum_{l=0}^p (-1)^{l}\binom{2p+1}{p-l} \left( \epsilon^{EST}_{j+1+l}  -\epsilon^{EST}_{j-l}\right)\right]\nonumber\\
    &= Var\left[(T_n^p \boldsymbol\epsilon^{EST})_j\right]=\alpha^{-1}\left(T_n^p\left(F^TF\right)^\dagger{T_n^p}^T\right)_{j,j}.\label{eq:var_Tx}   
\end{align}
with the last equality following from \cref{eq:var_eps}.

We note again that $p = 0$ is used in all of our  numerical experiments (see \cref{rem: onlyp0}).  In this case we simply have
\[Var\left[(T_n^0\bm{x}^{EST})_j\right]=Var\left[\epsilon^{EST}_{j+1}-\epsilon^{EST}_{j}\right]
=\alpha^{-1}\begin{bmatrix}
    -1\\ 1
\end{bmatrix}^T(F^T F)^{\dagger}_{j:j+1,j:j+1}\begin{bmatrix}
    -1\\ 1
\end{bmatrix}, \]
where $j:j+1$ in the subscript refers to indices from $j$ to $j+1$ inclusively, and $(F^T F)^{\dagger}_{j:j+1, j:j+1}$ is the sub-block matrix of $(F^T F)^{\dagger}$.

Analogously,  $S_{n,1/2}^{\sigma_{2p+1}}$ is an unbiased estimator since for $\hat{f}_k$ given in \cref{eq: fourcoefdiscrete} and $\hat{\epsilon}_k$ similarly defined as $\hat{\epsilon}_k = \frac{1}{n}\sum_{j = 1}^n \epsilon^{EST}_j e^{iks_j} =(\hat{F}\boldsymbol\epsilon^{EST})_k$, $k = -\frac{n}{2},\cdots,\frac{n}{2}-1$, we have (also see \cite{archibald2002reducing})
\begin{eqnarray*}
E\left[(S_{n,1/2}^{\sigma_{2p+1}}\bm{x}^{EST})_j\right]
    &=&E\left\{ \pi i \sum_{k=-\frac{n}{2}}^{\frac{n}{2}} \operatorname{sgn}(k) \, \sigma\left(\frac{|k|\Delta s}{\pi}\right)sinc\left(\frac{|k|\Delta s}{2}\right)\left(\hat {f}_k+\hat\epsilon_k\right) \, e^{i k s_j} e^{i \tfrac{k\Delta s}{2}}\right\}\\
    &=&E\left[(S_{n,1/2}^{\sigma_{2p+1}}{\bm f})_j\right]+E\left[(S_{n,1/2}^{\sigma_{2p+1}} \boldsymbol\epsilon^{EST})_j\right]\\
    &=&(S_{n,1/2}^{\sigma_{2p+1}} {\bm f})_j. 
\end{eqnarray*}

As in the case for the local differencing estimator, $Var(S_{n,1/2}^{\sigma_{2p+1}}{\bm x}^{EST})$ depends on the noise variance and order $2p+1$, but unlike the case for $T_n^p$ the effects are {\em global}. 
The estimator transforms the noise into the Fourier domain, yielding $\hat{\boldsymbol{\epsilon}}$ with the covariance matrix given by
\begin{eqnarray}
    \label{eq:cov_hat_eps}
    Cov(\hat{\boldsymbol{\epsilon}})=Cov(\hat{F}\boldsymbol{\epsilon}^{EST})=\hat{F}Cov(\boldsymbol{\epsilon})\hat{F}^H=\alpha^{-1}\hat{F}(F^T F)^{\dagger}\hat{F}^H,
\end{eqnarray}
where $\hat{F}^H$ is the conjugate transpose of the discrete Fourier transform matrix $\hat{F}$ defined in \cref{eq:discretefourMV}.
The global impact is thus evident from direct calculation:
\begin{eqnarray*}
Var\left[(S_{n,1/2}^{\sigma_{2p+1}} \bm{x}^{EST})_j\right]
&=&E \left[\pi i \sum_{k=-\frac{n}{2}}^{\frac{n}{2}} c_k\left(\hat f_k+\hat\epsilon_k\right) -E\left[(S_{n,1/2}^{\sigma_{2p+1}}\bm x^{EST})_j\right] \right]\\
  &\cdot & \left[\pi i \sum_{l=-\frac{n}{2}}^{\frac{n}{2}} c_l\left(\hat f_l+\hat\epsilon_l\right)-E\left[(S_{n,1/2}^{\sigma_{2p+1}}\bm x^{EST})_j\right]\right]^\ast\\
&=&E \left[\left(\pi i \sum_{k=-\frac{n}{2}}^{\frac{n}{2}} c_k\hat\epsilon_k\right)
    \left( \pi i \sum_{l=-\frac{n}{2}}^{\frac{n}{2}}c_l\hat\epsilon_l\right)^\ast\right]
    \\
& = & \pi^2\sum_{k=-\frac{n}{2}}^{\frac{n}{2}}\sum_{l=-\frac{n}{2}}^{\frac{n}{2}}c_k c_l^\ast E\left[\hat\epsilon_k\hat\epsilon_l^\ast\right]
    =\tfrac{\pi^2}{\alpha} \sum_{k=-\frac{n}{2}}^{\frac{n}{2}}\sum_{l=-\frac{n}{2}}^{\frac{n}{2}}c_k c_l^\ast (\hat{F} (F^T F)^{\dagger} \hat{F}^H)_{kl}
\end{eqnarray*}
where $\ast$ denotes the complex conjugate
and $c_k=\operatorname{sgn}(k) \, \sigma\left(\frac{|k|\Delta s}{\pi}\right)sinc\left(\frac{|k|\Delta s}{2}\right)e^{i k s}e^{i \tfrac{k\Delta s}{2}}$. The final equality follows from \eqref{eq:cov_hat_eps}. By defining column vector  $\bm c \in \mathbb{C}^{n}$ with entries $\{c_k\}_{k = -\frac{n}{2}}^{\frac{n}{2}-1}$, we can now formulate a quadratic form similar to \eqref{eq:var_Tx} as
\begin{eqnarray*}
    Var\left[(S_{n,1/2}^{\sigma_{2p+1}} \bm{x}^{EST})_j\right]=\tfrac{\pi^2}{\alpha}\bm c^T(\hat{F} (F^T F)^{\dagger} \hat{F}^H)\bm c^\ast=\alpha^{-1}\left(S_{n,1/2}^{\sigma_{2p+1}}(F^TF)^\dagger{S_{n,1/2}^{\sigma_{2p+1}}}^T\right)_{j,j},
\end{eqnarray*}
which is mainly determined by noise precision $\alpha$ and concentration factor $\sigma$.

We are now able to approximate the impact of noise on the residual prior transform.  Starting with the component calculation of $E(R_{n,1/2}^p{\bm x}^{EST})$, from linearity we have that
\begin{eqnarray*}
    E\left[(R_{n,1/2}^p{\bm x}^{EST})_j\right]&=&E\left[\left(\left(T_n^p-S_{n,1/2}^{\sigma_{2p+1}}\right){\bm x}^{EST}\right)_j\right]=  E\left[(T_n^p\bm{x}^{EST})_j\right]-  E\left[(S_{n,1/2}^{\sigma_{2p+1}}\bm{x}^{EST})_j\right] \\
    &=& (T_n^p{\bm f})_j - (S_{n,1/2}^{\sigma_{2p+1}}{\bm f})_j = 0.
\end{eqnarray*}
To calculate $Var\left[(R_{n,1/2}^p{\bm x}^{EST})_j\right]$, we first observe that
\begin{eqnarray*}
&&Cov\left((T_n^p{\bm x}^{EST})_j,(S_{n,1/2}^{\sigma_{2p+1}}{\bm x}^{EST})_j\right)\\
&=&E\left[\left((T_n^p{\bm x}^{EST})_j-E\left[(T_n^p{\bm x}^{EST})_j\right]\right)
\left((S_{n,1/2}^{\sigma_{2p+1}}{\bm x}^{EST})_j-E\left[(S_{n,1/2}^{\sigma_{2p+1}}{\bm x}^{EST})_j\right]\right)
\right]\\
&=&E\left[ (T^n_p\boldsymbol{\epsilon}^{EST})_j(S_{n,1/2}^{\sigma_{2p+1}}\boldsymbol{\epsilon}^{EST})_j \right]
=\left(T^p_n Cov(\boldsymbol{\epsilon}^{EST}){S_{n,1/2}^{\sigma_{2p+1}}}^T\right)_{j,j}\\
&=&\alpha^{-1}\left(T^p_n\left(F^TF\right)^\dagger {S_{n,1/2}^{\sigma_{2p+1}}}^T\right)_{j,j}, \end{eqnarray*}
where the last equality again follows from \cref{eq:var_eps}, just as in \cref{eq:var_Tx}.
The above covariance is central to understanding the behavior of the residual prior transform. By design, the estimators $T^p_n$ and $S_{n,1/2}^{\sigma_{2p+1}}$ provide close approximations of each other in the signal's smooth regions, meaning their outputs are highly positively correlated. This large, positive covariance term directly offsets the sum of the individual variances, leading to a significantly suppressed variance of the residual prior transform in smooth regions,
\begin{eqnarray*}
    &&Var\left[(R_{n,1/2}^p\bm x^{EST})_j\right]=Var\left[\left(\left(T_n^p-S_{n,1/2}^{\sigma_{2p+1}}\right){\bm x}^{EST}\right)_j\right]\\
    &=&Var\left[ \left(T_n^p {\boldsymbol \epsilon}^{EST}\right)_j\right]+Var\left[\left(S_{n,1/2}^{\sigma_{2p+1}} {\boldsymbol \epsilon}^{EST}\right)_j\right]-2Cov\left((T_n^p{\boldsymbol \epsilon}^{EST})_j,(S_{n,1/2}^{\sigma_{2p+1}}{\boldsymbol \epsilon}^{EST})_j\right).
\end{eqnarray*}
In contrast, around a jump discontinuity, the agreement between the two estimators $T^p_n$ and $S_{n,1/2}^{\sigma_{2p+1}}$ is disrupted. This loss of correlation prevents the variance cancellation seen in smooth regions, making the variance of $R_{n,1/2}^p$ pronounced specifically at the jump. This resulting ``adaptive variance'', which is suppressed in smooth regions and amplified at discontinuities, is fundamentally advantageous for signal recovery. By effectively leveraging the complementary information from both estimators, $R_{n,1/2}^p$ becomes a more sensitive and robust detector than either estimator used alone. 
This property directly enhances the ``separation of scales'', the critical principle, that allows us to distinguish genuine jump discontinuities from smooth variations. In particular, $R_{n,1/2}^p\bm x$ is subject to a much larger variance at the location of a discontinuity compared to its minimal variance in smooth regions. This creates a sharper contrast, making the locations of jumps easier to identify and leading to a more reliable signal recovery.

\section{Numerical results} 
\label{sec:numerics} 

We now evaluate the performance of the residual prior transform within the hierarchical sparse Bayesian learning framework. The results are computed using the generalized sparse Bayesian learning algorithm \cref{alg:GSBL} for individual reconstruction and the MMV version \cref{alg:MMV_GSBL} for joint recovery. We examine  the efficacy of our residual prior transform for both 1D and 2D recoveries and also explore the benefits of leveraging joint information across multiple measurement vectors (MMVs).

To provide a comprehensive analysis, our results are organized in two parts. The first focuses on 1D signals, where we also visualize the uncertainty quantification of our reconstructions. The second extends our investigation to 2D images, showcasing the effectiveness of the residual prior transform in higher dimensions. Moreover, we consider representative multimodal problems. For display purposes, we apply test suites to $L=3$ measurements, though the framework readily accommodates a larger number of modalities. Our experiments cover both the recovery of a single function from multimodal data and the joint recovery of multiple functions sharing common structural features. Each of the three measurements is simulated using a distinct forward operator to model fundamental signal processing tasks. Particularly, we exploit tests of denoising,  deblurring, and undersampling, for which we investigate both spatial undersampling and the recovery from incomplete Fourier data.

The model is always given by \cref{eq: data_acquisition}, with the linear forward transform  $F_l$, $l = 1,\dots,L$, dependent on the acquisition modality.  For denoising, we use the identity matrix  $F_l =I_n$. For deblurring, we use the matrix $F_l = \tilde F \in \mathbb{R}^{n \times n}$ to apply a Gaussian blur, defined by its point spread function (psf) in
\begin{equation}
 \label{eq:bluroperator}   
\tilde{F}_{ij}= \frac{1}{2\pi\gamma^2}\exp{\left(-\frac{i^2+j^2}{2\gamma^2}\right)}. \end{equation}
This process smooths out high frequency details, making the underlying structure less distinct and introducing ill-posedness. The extent of this smoothing is controlled by the psf parameter $\gamma$. We also consider a spatial undersampling case, as defined by $F_l =HI_n\in\R^{m\times n}$, where $H$ randomly removes rows at a ratio $r$ such that  $m=\operatorname{round}((1-r)n)$, leading to a loss of information. Finally, for incomplete Fourier data, the forward operator is analogously defined as $F_l=H\hat F$ with $\hat F$ in \eqref{eq: DFT}.

The noise precision $\alpha$ in \cref{eq: iid_noise} is determined by the signal to noise ratio (SNR)
\[\text{SNR}_\text{dB} = 10\log_{10}\left(\tfrac{P_{signal}}{P_{noise}}\right) = 10\log_{10}\left(\frac{\alpha\norm{\bm f}^2}{n}\right).\]
For notational clarity, we omit the subscript $l$, as the following definition holds for each of the $L$ measurements.
We evaluate the performance of each recovery using the pointwise error
\begin{equation}
    \label{eq: abs_err}
   E_{j}^\text{abs}=\abs{\bm x(j)-f(s_j)}, \quad j = 1,\dots,n,
\end{equation}
where $\bm x(j)$ is the $jth$ component of the 
recovery obtained from either \cref{alg:GSBL} or \cref{alg:MMV_GSBL}, and $f(s_j)$ is the $ith$ component of  the discretized true solution $f$ in \cref{eq:discretization_f}.

In all figures, solutions with the superscript $loc$ are obtained using the local differencing prior transform  $\Phi= T_n^p$ in \eqref{eq:x_update}, which is applied within either \Cref{alg:GSBL} ($\bm x^{loc}_\text{GSBL}$) for separate recovery or \cref{alg:MMV_GSBL} ($\bm x^{loc}_\text{MMV-GSBL}$) for joint recovery. 
Similarly, solutions using the residual prior transform   $\Phi= R^p_{n,\zeta}$, as given by \cref{def:residtransform}, are designated with the superscript $R$. This gives us $\bm x^{R}_\text{GSBL}$ for separate recovery and $\bm x^{R}_\text{MMV-GSBL}$ for joint recovery, respectively. The discretization resolution $n$ was chosen based on the complexity of the data. A resolution of $n=128$ was used for both 1D and 2D synthetic data, whereas the real-world tests were processed at a higher resolution of $n=400$ to better capture their intricate details. Whenever the residual prior transform $\Phi=R_{n,\zeta}^p$ was employed, the shift parameter was fixed at $\zeta=1/4$.
Other parameter choices produce consistent outcomes.

For individual recovery with \cref{alg:GSBL}, we constrained the distinct hyperparameters $\vartheta_l$ to a constant value, thereby mirroring the inherently shared hyperparameter $\vartheta$ of the joint recovery algorithm \cref{alg:MMV_GSBL}. Particularly, we set $\vartheta=10^{-4}$ for 1D signal recovery problems. For 2D image recovery problems, $\vartheta$ was adapted to the specific prior transform. We set $\vartheta = 10^{-2}$ for $\Phi=T_n^p$, a value consistent with common practice in related hierarchical sparse Bayesian learning algorithms \cite{tipping2001sparse,babacan2010sparse,bardsley2012mcmc}. However, we found that $\vartheta=10^{-3}$ yielded better results for $\Phi=R_{n,\zeta}^p$. This aligns with the Lasso regression recovery estimates in \cite{xiao2025new}, where the residual transform operator requires a stronger regularization penalty to achieve optimal performance.

\subsection{Joint Recovery of Multimodal One-dimensional Signals}
\label{sec:1dsignal}

\begin{figure}[h!]
    \centering
    \begin{subfigure}[b]{.30\textwidth}
        \includegraphics[width=\textwidth]{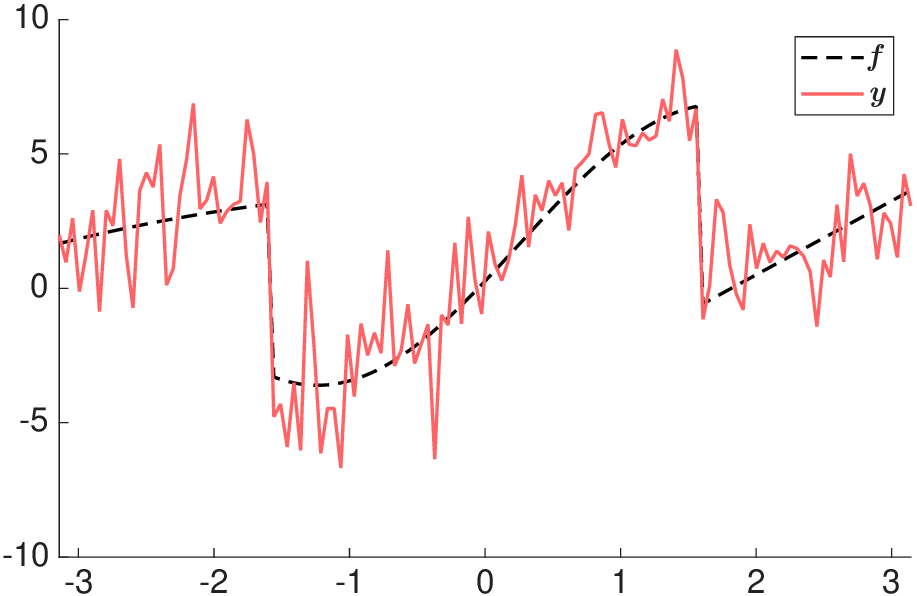}
    \end{subfigure}
    \begin{subfigure}[b]{.30\textwidth}
        \includegraphics[width=\textwidth]{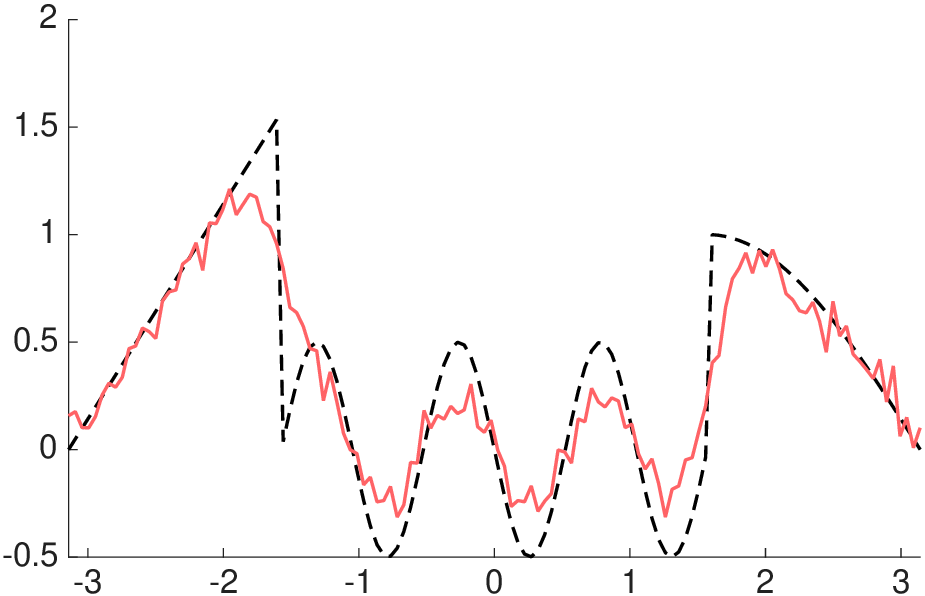}
    \end{subfigure}
    \begin{subfigure}[b]{.30\textwidth}
        \includegraphics[width=\textwidth]{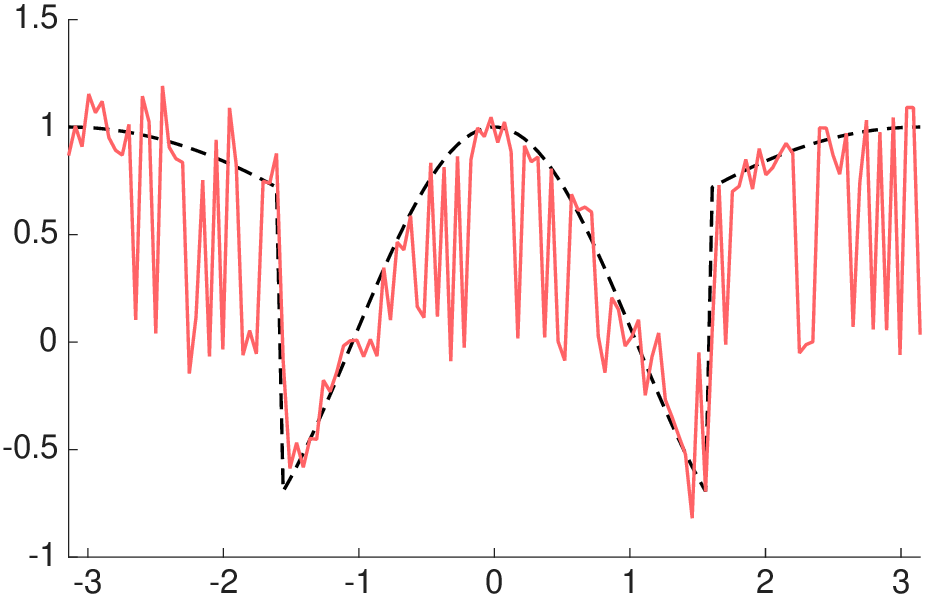}
    \end{subfigure}
        \begin{subfigure}[b]{.30\textwidth}
        \includegraphics[width=\textwidth]{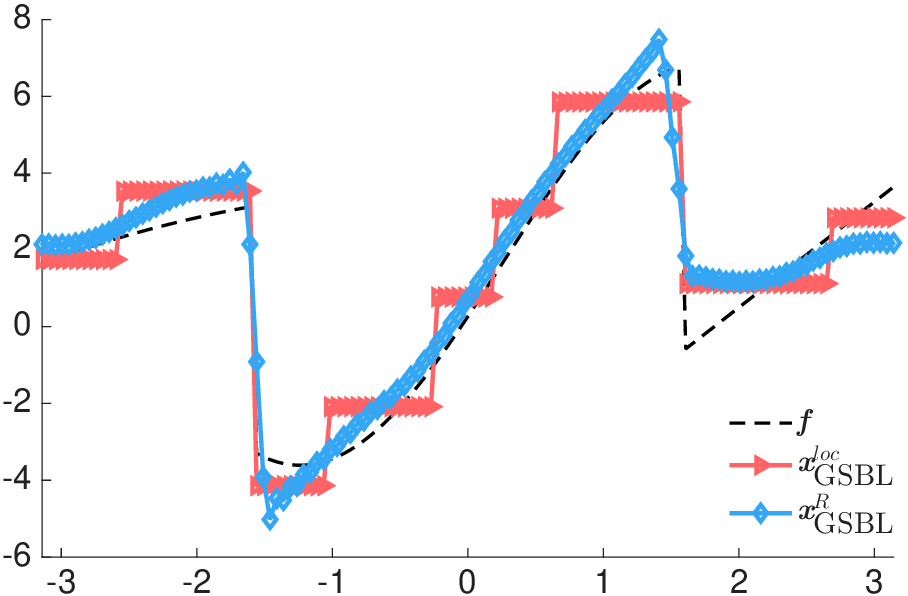}
    \end{subfigure}
    \begin{subfigure}[b]{.30\textwidth}
        \includegraphics[width=\textwidth]{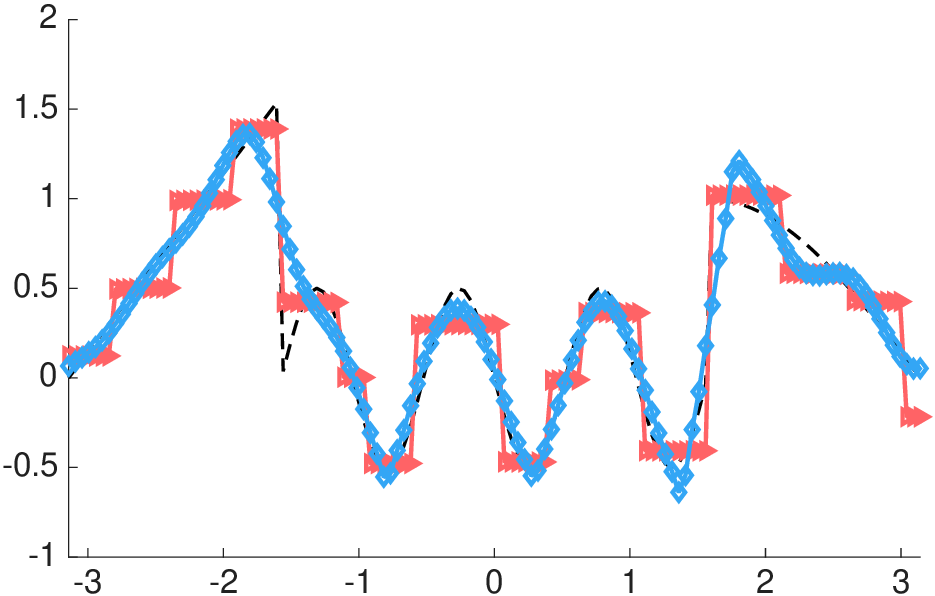}
    \end{subfigure}
    \begin{subfigure}[b]{.30\textwidth}
        \includegraphics[width=\textwidth]{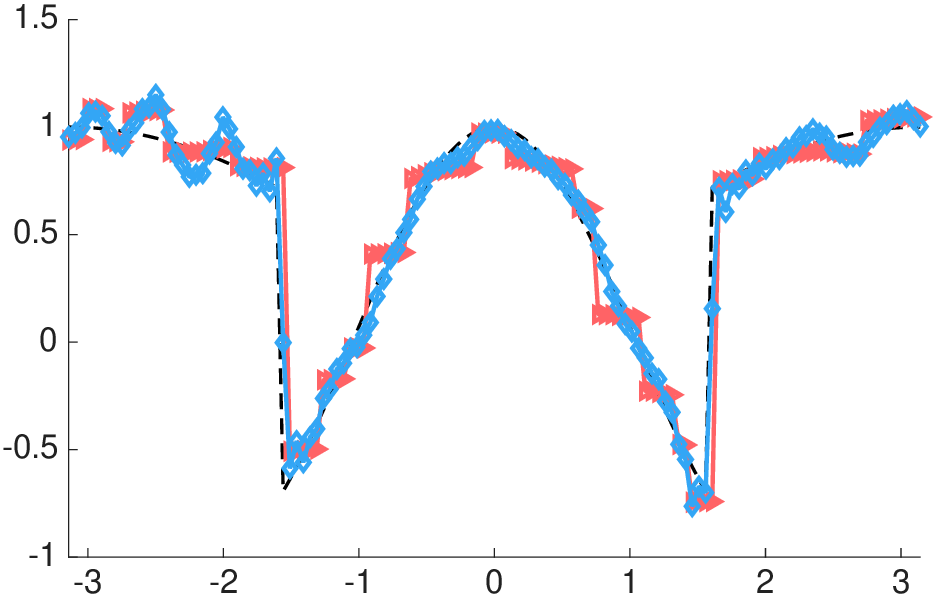}
    \end{subfigure}\\[.5em]
    \begin{subfigure}[b]{.30\textwidth}
        \includegraphics[width=\textwidth]{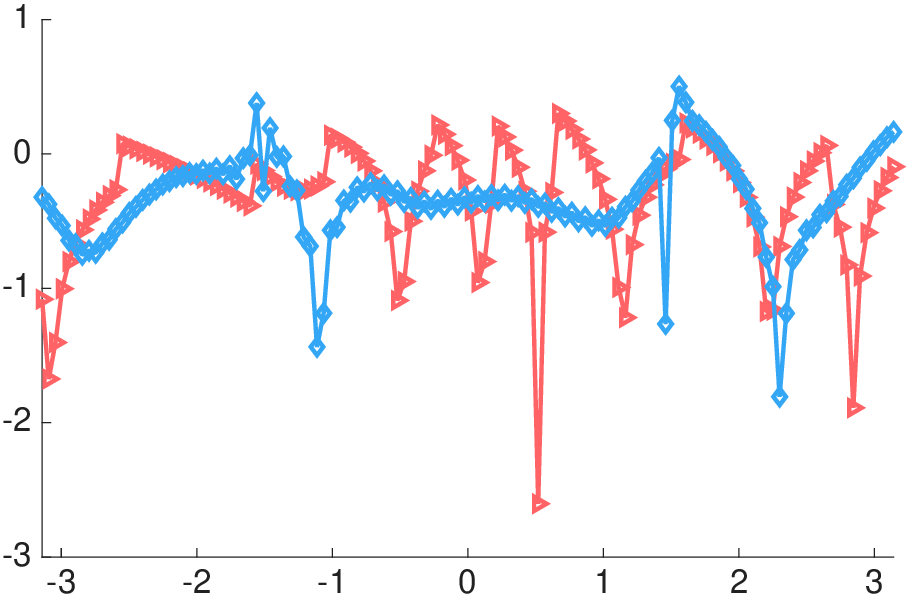}
    \end{subfigure}
    \begin{subfigure}[b]{.30\textwidth}
        \includegraphics[width=\textwidth]{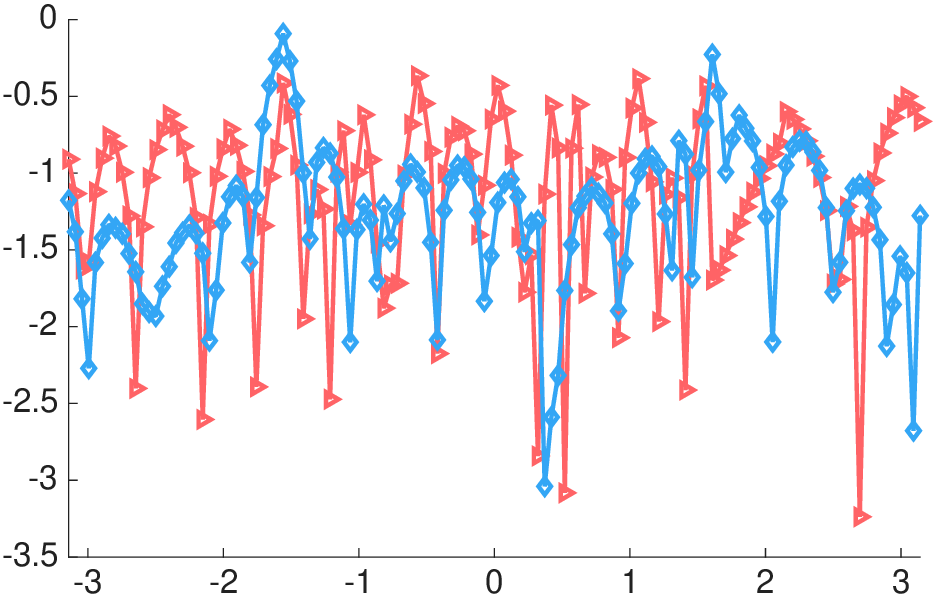}
    \end{subfigure}
    \begin{subfigure}[b]{.30\textwidth}
        \includegraphics[width=\textwidth]{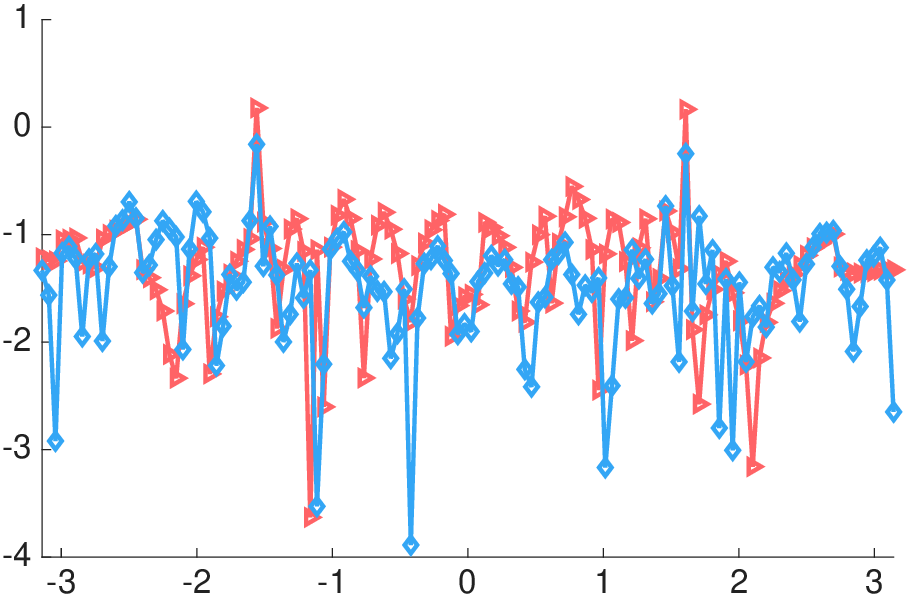}
    \end{subfigure}
    \caption{(top) The red lines show the  measurements, $\bm y_l$, $l = 1,2,3$  as given by \eqref{eq: data_acquisition} for corresponding $f_l(s)$ in \Cref{ex: example_1d}.  (top-left) Noisy data for the piecewise smooth signal $f_1$, measured using the identity matrix $F_1 = I_n$ and contaminated by additive noise at an $\text{SNR}=5$dB. (top-middle) Blurred data for signal $f_2$, using the Gaussian blur operator $F_2=\tilde F$ with psf parameter  $\gamma=0.03$, as defined in \eqref{eq:bluroperator}, and contaminated by additive noise at an $\text{SNR}=20$dB. (top-right) Undersampled data for signal $f_3$, obtained with the operator $F_3=HI_n$ with undersampling ratio  $r=0.3$, and contaminated by additive noise at an $\text{SNR}=20$dB. The dashed black lines represent the corresponding discretized true underlying signals, $\bm f_l$,  $l=1,2,3$. (middle row) MAP estimates obtained by \cref{alg:GSBL}. (bottom) Corresponding $\log_{10}E^{\text{abs}}_j$, $j=1,\dots,n$,  given by \eqref{eq: abs_err}.}
    \label{fig:trueSignals}
\end{figure}

Our first test suite is designed to evaluate the algorithm's performance on recovering multiple functions that share structural properties. Specifically, \Cref{ex: example_1d} considers three functions have common discontinuity locations but exhibit unique variability within their smooth regions. 
\begin{example}
\label{ex: example_1d}
   Let $f_1(s)$, $f_2(s)$, and $f_3(s)$ be $2\pi-$periodic functions defined over $[-\pi,\pi)$ where
    \begin{equation*}
        f_1(s)=\begin{cases}
        \tfrac{11\pi}{4} - 5 - \tfrac{s^2}{5}, & s\in U_1\\
        \tfrac{7}{4} - \tfrac{s}{2} + 6\sin(s-\tfrac{1}{4}), & s\in U_2\\
        \tfrac{11s}{4} - 5, & s\in U_3
        \end{cases},
        f_2(s)=\begin{cases}
        s+\pi, & s\in U_1\\
        \tfrac{-\sin(6s)}{2}, & s\in U_2\\
        \sin(-s+\pi), & s\in U_3
        \end{cases}, 
        f_3(s)=\begin{cases}
        \sin(-\tfrac{s}{2}), & s\in U_1\\
        \cos(\tfrac{3s}{2}), & s\in U_2\\
        \sin(\tfrac{s}{2}), & s\in U_3
        \end{cases},
    \end{equation*}
    with  $U_1=[-\pi,-\tfrac{\pi}{2})$, $U_2=[-\tfrac{\pi}{2},\tfrac{\pi}{2})$, and $U_3=[\tfrac{\pi}{2},\pi)$.  
\end{example}

\Cref{fig:trueSignals} (top) offers a visual introduction to the multimodal data acquisition scenarios in our numerical investigation for \cref{ex: example_1d}. Each panel displays a 1D signal from \cref{ex: example_1d} alongside the corresponding measurement resulting from a specific type of degradation. While the true signals vary, the types of degradation present here, including the forward operators and noise levels, are used for all reconstructions in \cref{sec:1dsignal}. The significant fluctuations in the measurement for $f_3$ are a result of spatial undersampling, and {\em not} a low signal to noise ratio. The multimodal problem forms the foundation for evaluating the performance of our residual prior transform and highlights how leveraging joint information across the datasets leads to enhanced results.

The middle and bottom rows of \Cref{fig:trueSignals} compare recoveries  and their corresponding spatial pointwise $\log_{10}$ absolute errors \eqref{eq: abs_err} across three degradation scenarios using \Cref{alg:GSBL}.
The ``staircasing'' artifact commonly observed when employing the first-order local differencing prior transform $\Phi=T_n^p$ is pronounced 
in each of the three cases. 
Collectively, the recoveries computed by \cref{alg:GSBL} using the residual prior transform $\Phi=R_{n,\zeta}^p$ consistently mitigate the staircasing artifacts and yield greater detail in signal's variable regions compared to the standard first-order differencing prior transform. While this comes at the cost of compromising some local structure, the overall advantage is quantitatively supported by consistently smaller absolute errors across multimodal conditions. 
\begin{figure}[h!]
    \centering
    \begin{subfigure}[b]{.32\textwidth}
        \includegraphics[width=\textwidth]{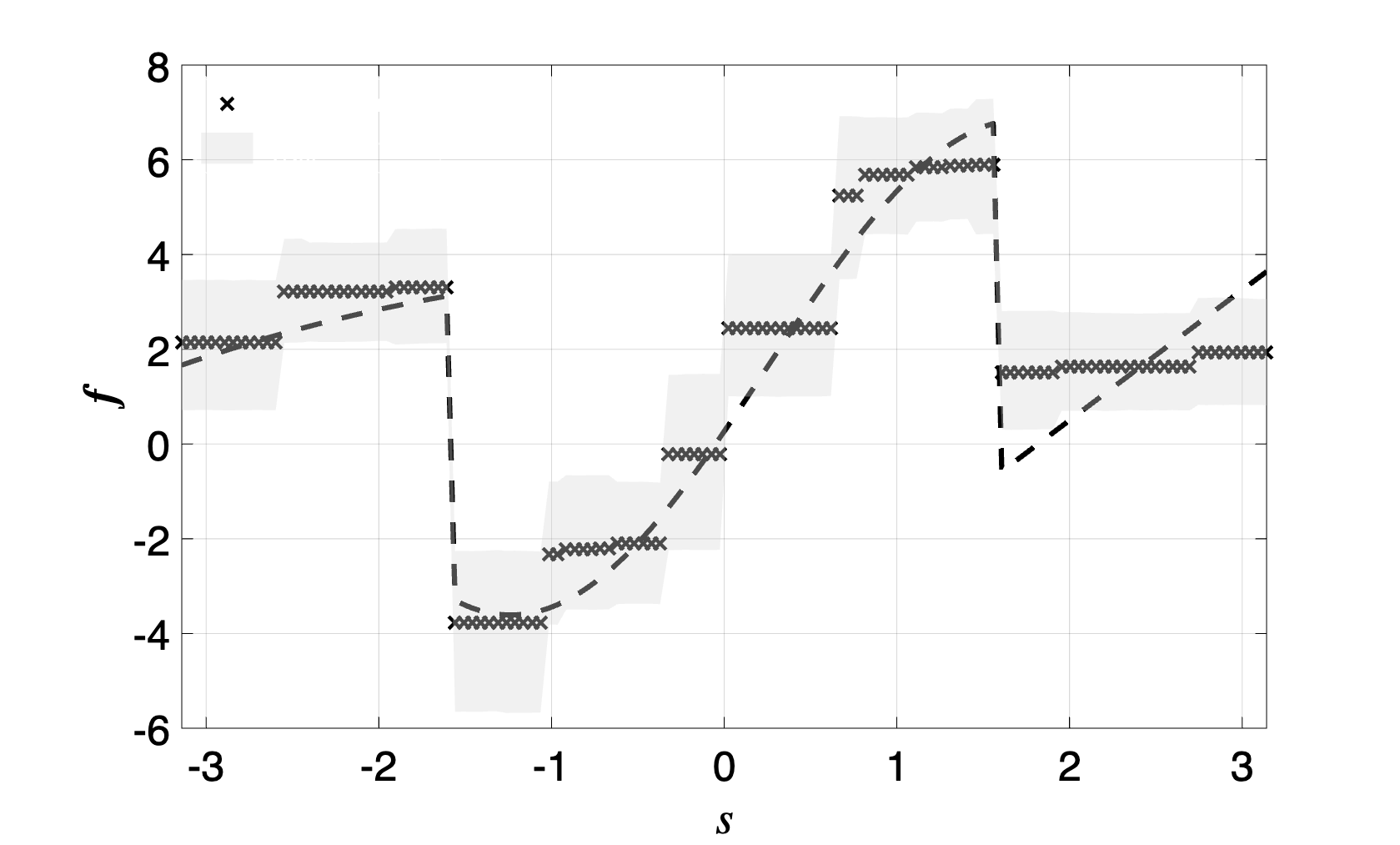}
    \end{subfigure}
    \begin{subfigure}[b]{.32\textwidth}
        \includegraphics[width=\textwidth]{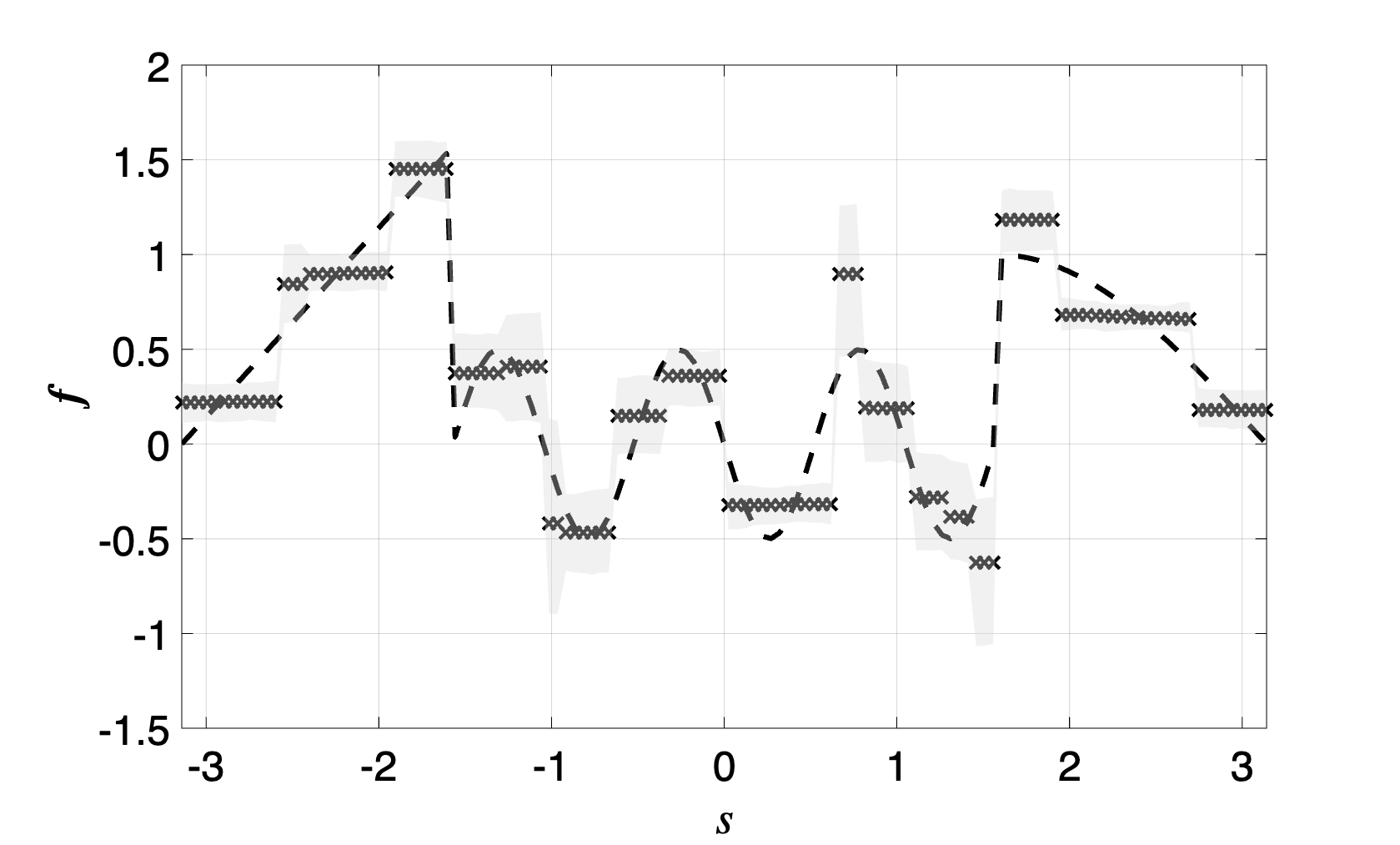}
    \end{subfigure}
    \begin{subfigure}[b]{.32\textwidth}
        \includegraphics[width=\textwidth]{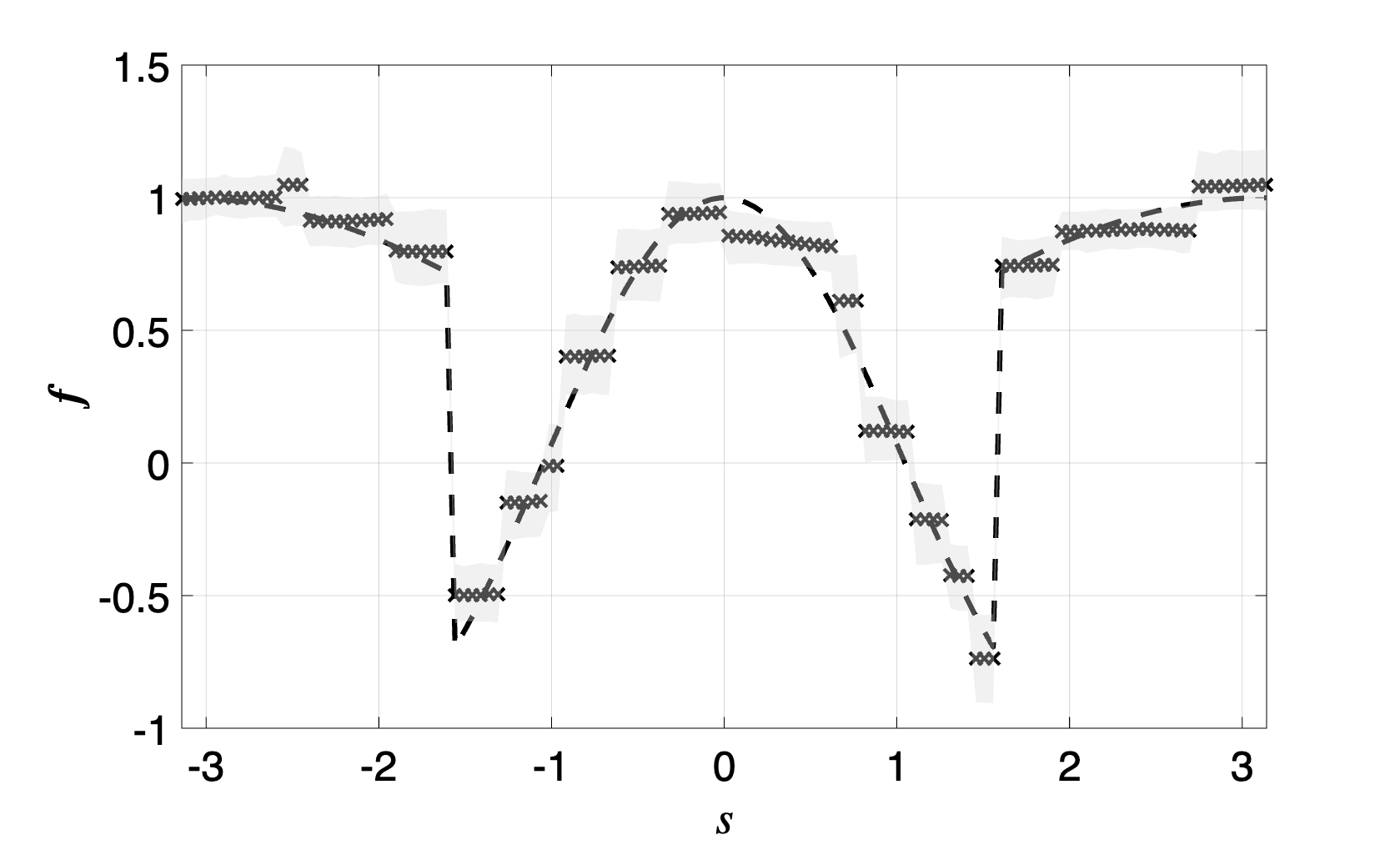}
    \end{subfigure}\\
    \begin{subfigure}[b]{.32\textwidth}
        \includegraphics[width=\textwidth]{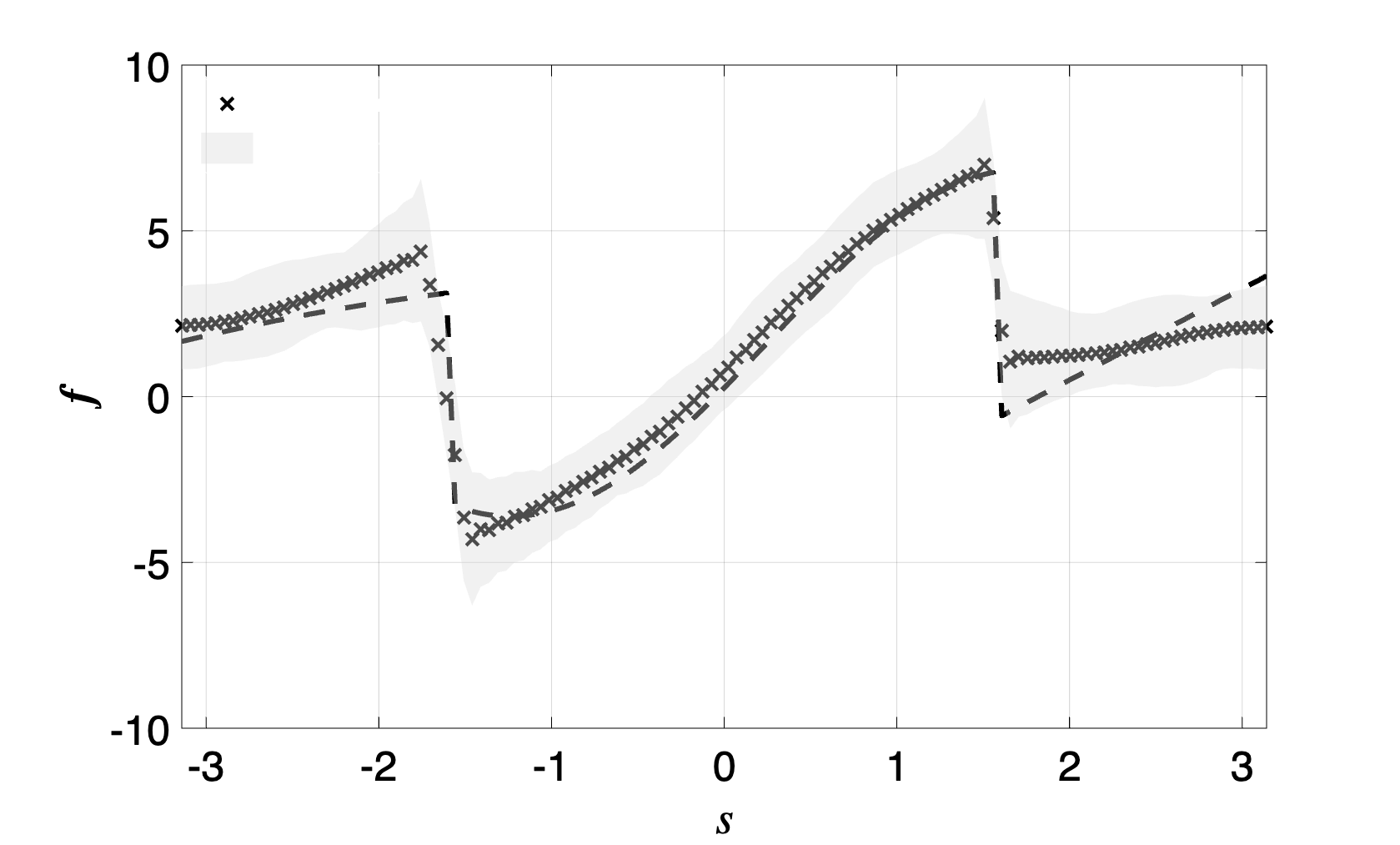}
    \end{subfigure}
    \begin{subfigure}[b]{.32\textwidth}
        \includegraphics[width=\textwidth]{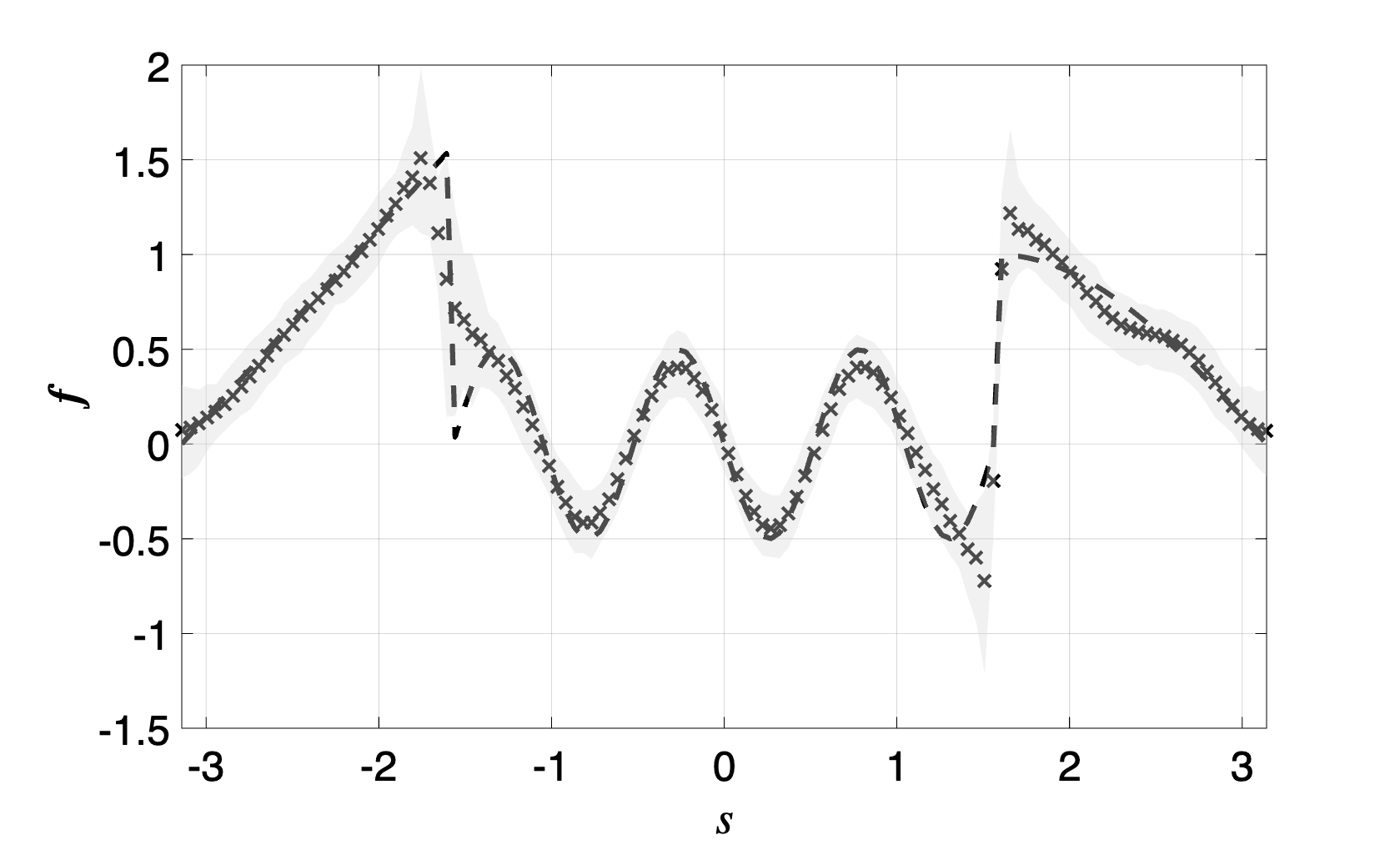}
    \end{subfigure}
    \begin{subfigure}[b]{.32\textwidth}
        \includegraphics[width=\textwidth]{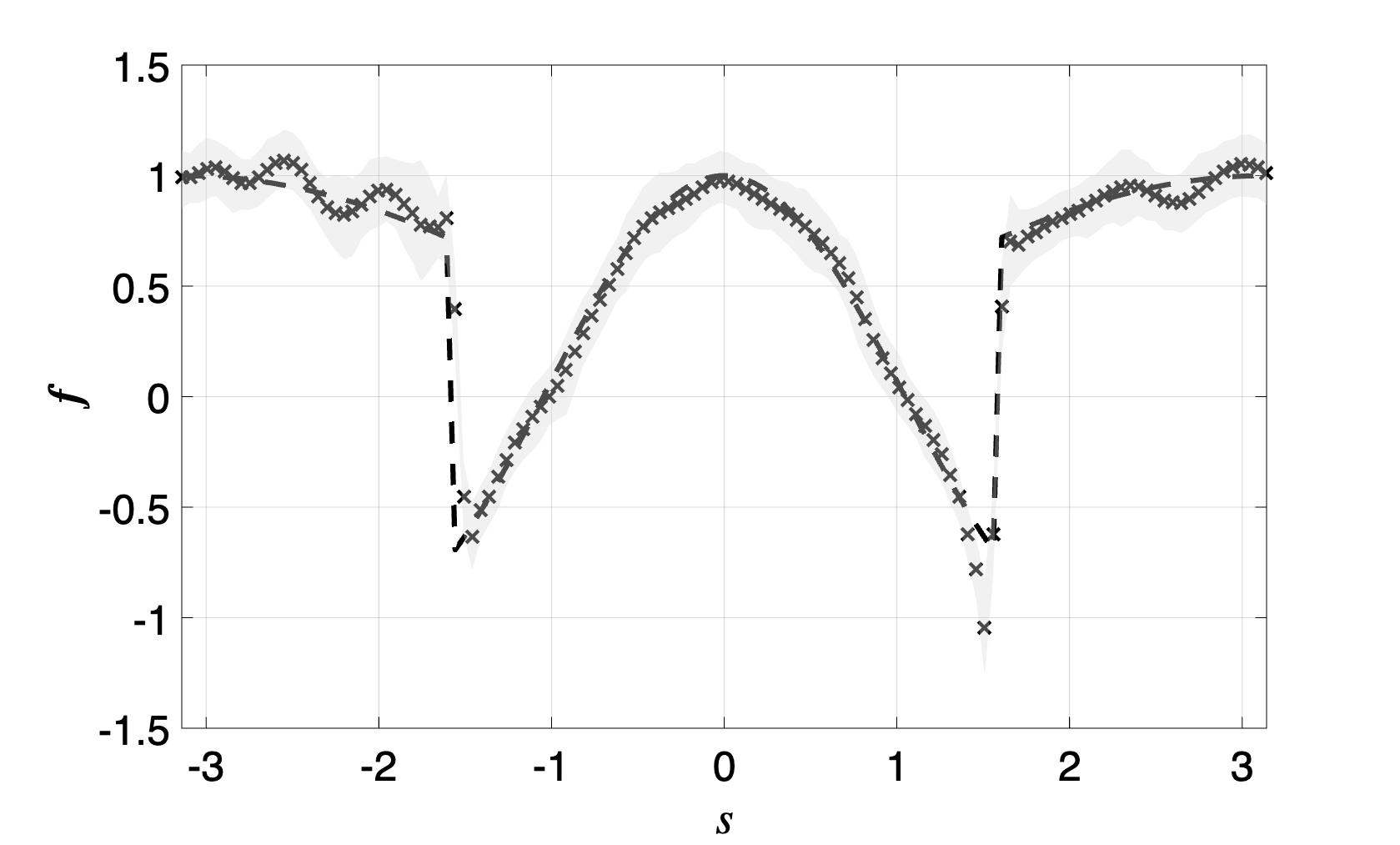}
    \end{subfigure}
    \caption{The MAP estimates for (left) $f_1$, (middle) $f_2$, and (right) $f_3$ and 99\% credible intervals for the recovered signals, using the local differencing prior transform (top) $\Phi = T_n^p$ and the residual prior transform (bottom) $\Phi = R_{n,\zeta}^p$.}
    \label{fig:UQ}
\end{figure}

\Cref{fig:UQ} illustrates uncertainty quantification results for the MAP solutions obtained by \cref{alg:GSBL}. Each subfigure in \cref{fig:UQ} displays both the MAP estimate in \cref{fig:trueSignals} (middle row) and its associated 99\% credible interval, computed from \eqref{eq: pos_mean_var}. Specifically, \cref{fig:UQ} (top) illustrates the UQ for the MAP solution obtained using $\Phi=T_n^p$, while \cref{fig:UQ} (bottom) provides the UQ for the MAP solution utilizing $\Phi=R_{n,\zeta}^p$. These credible intervals, shown as shaded regions around the MAP estimate, offer a comprehensive measure of the estimation uncertainty for each choice of prior transform.

\begin{figure}[h!]
    \centering
    \begin{subfigure}[b]{.32\textwidth}
        \includegraphics[width=\textwidth]{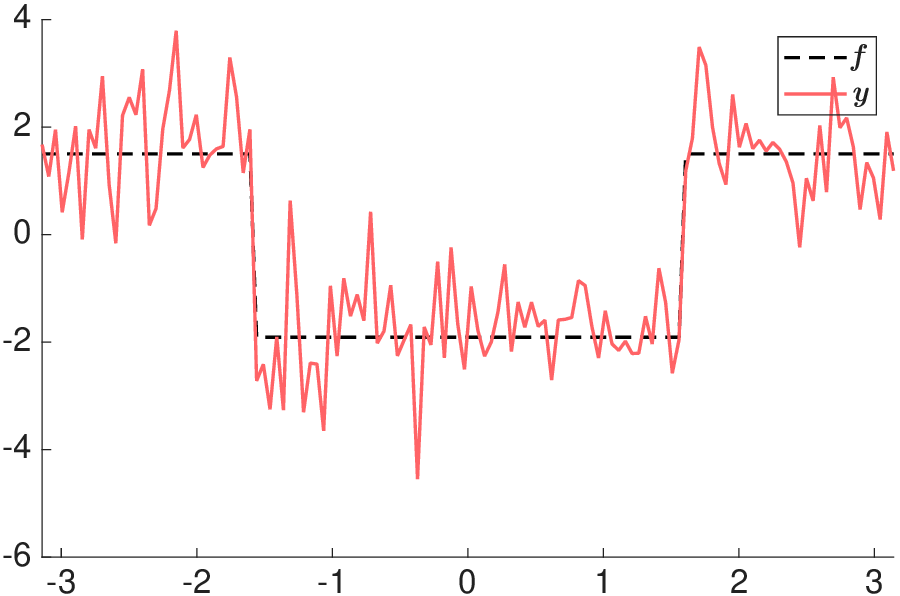}
    \end{subfigure}
    \begin{subfigure}[b]{.32\textwidth}
        \includegraphics[width=\textwidth]{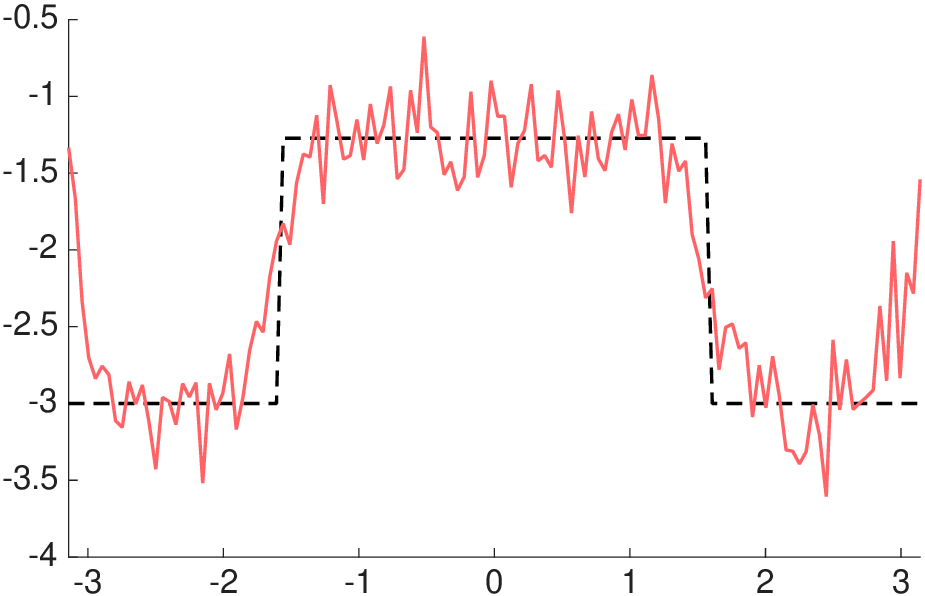}
    \end{subfigure}
    \begin{subfigure}[b]{.32\textwidth}
        \includegraphics[width=\textwidth]{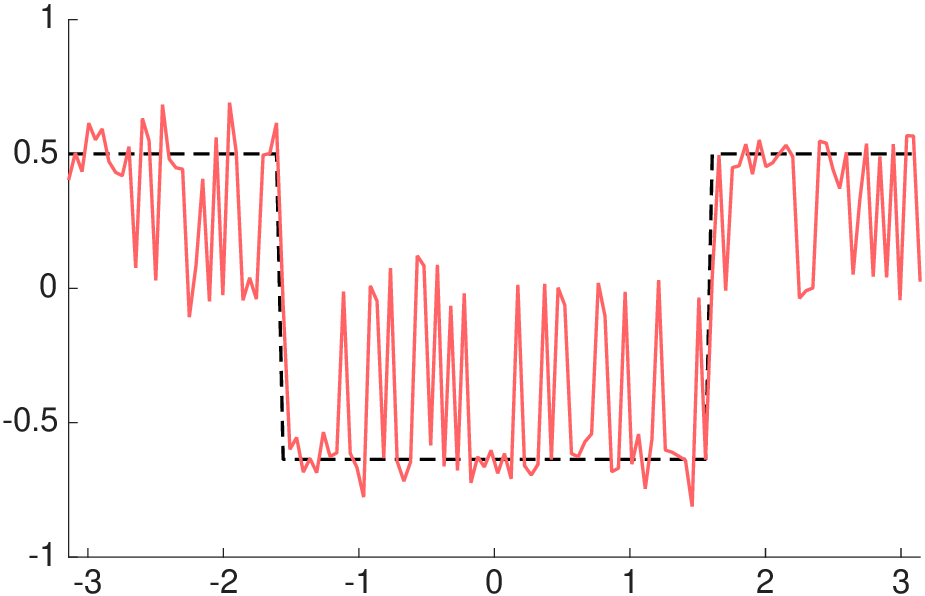}
    \end{subfigure}\\[.5em]
    \begin{subfigure}[b]{.32\textwidth}
        \includegraphics[width=\textwidth]{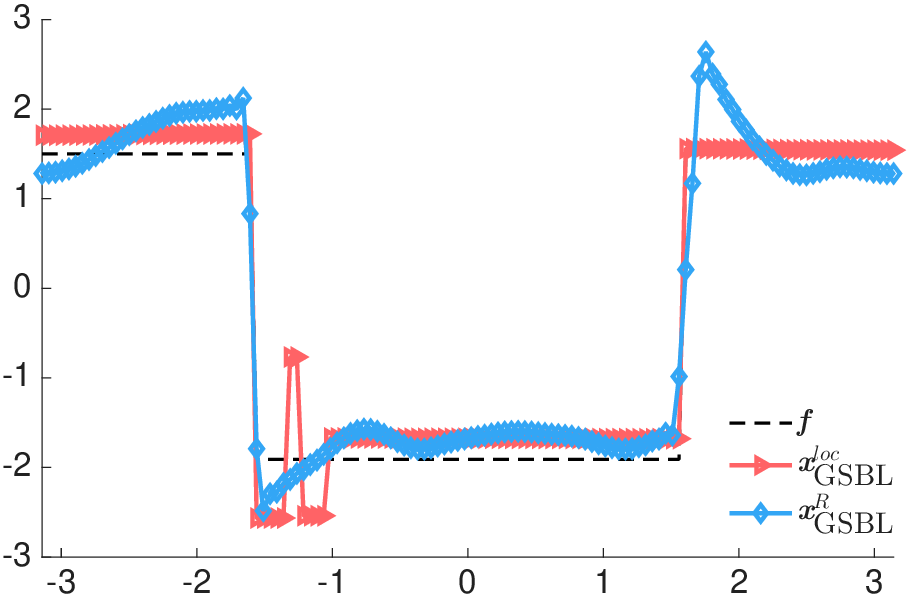}
    \end{subfigure}
    \begin{subfigure}[b]{.32\textwidth}
        \includegraphics[width=\textwidth]{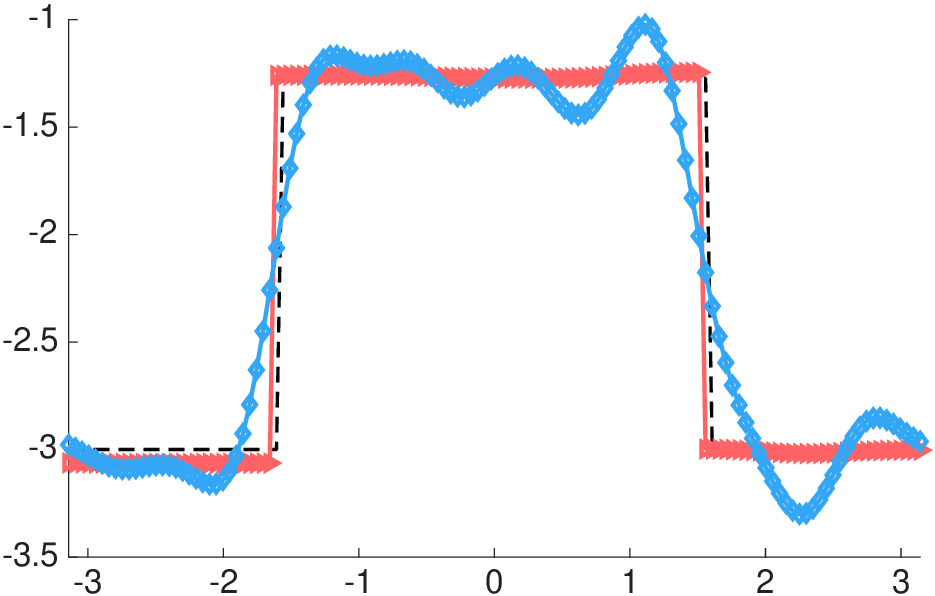}
    \end{subfigure}
    \begin{subfigure}[b]{.32\textwidth}
        \includegraphics[width=\textwidth]{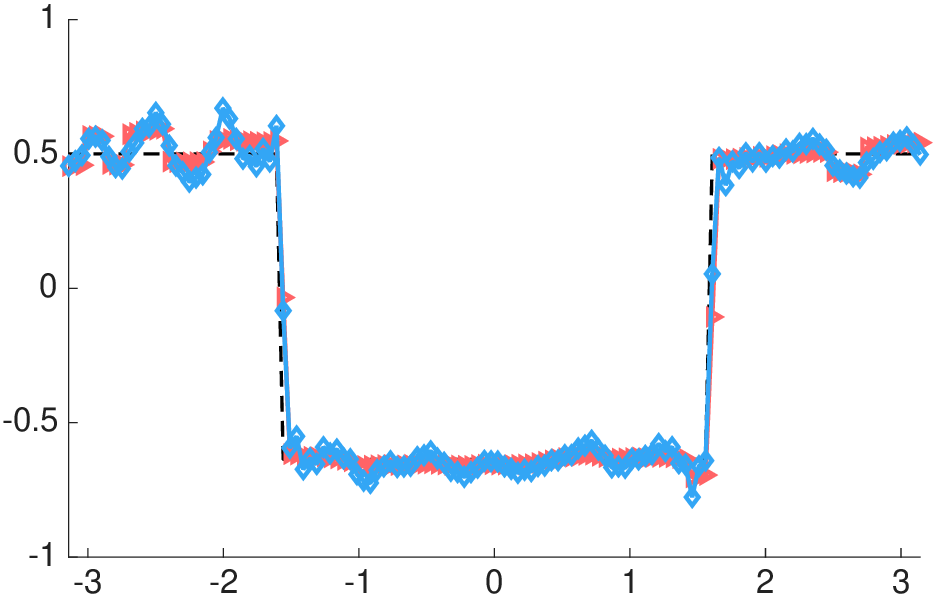}
    \end{subfigure}\\[.5em]
    \begin{subfigure}[b]{.32\textwidth}
        \includegraphics[width=\textwidth]{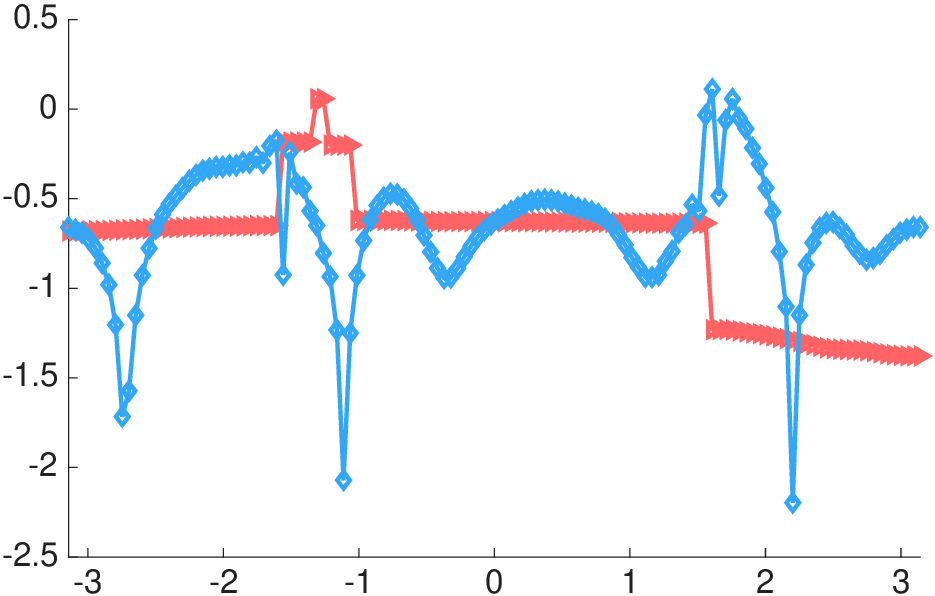}
    \end{subfigure}
    \begin{subfigure}[b]{.32\textwidth}
        \includegraphics[width=\textwidth]{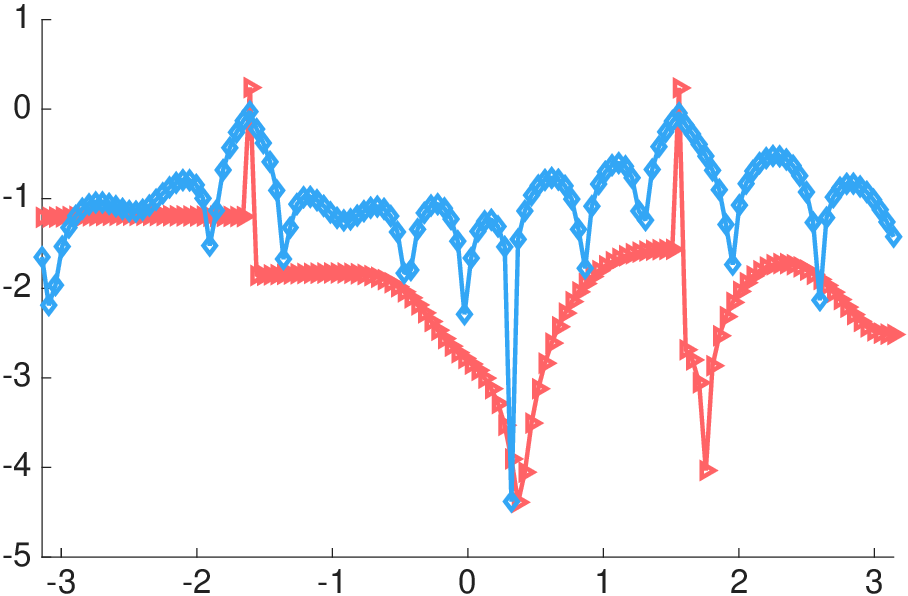}
    \end{subfigure}
    \begin{subfigure}[b]{.32\textwidth}
        \includegraphics[width=\textwidth]{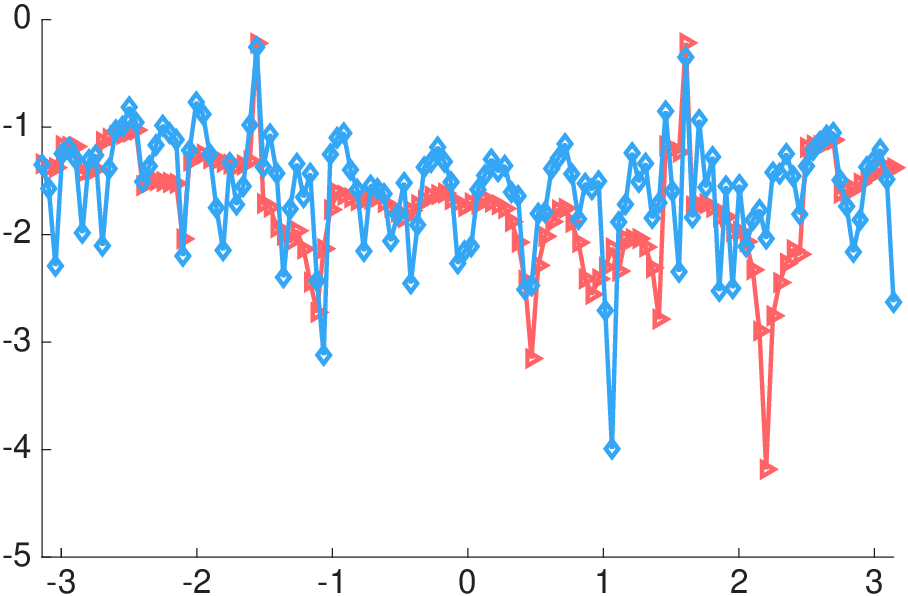}
    \end{subfigure}
    \caption{(top) The red lines show the  measurements in \cref{eq: data_acquisition} that  correspond to $f_{4}(s), f_5(s)$ and $f_6(s)$ in \Cref{ex: example_1d_const}. The multimodal data collections are described in \cref{fig:trueSignals}. (middle row) MAP estimates obtained by \cref{alg:GSBL} from data $\bm y$.  (bottom) Corresponding $\log_{10}E^{\text{abs}}_j$, $j=1,\dots,n$,  given by \eqref{eq: abs_err}.}
    \label{fig:rec_pw_consant}
\end{figure}

\begin{figure}[h!]
    \centering
    \begin{subfigure}[b]{.32\textwidth}
       \includegraphics[width=\textwidth]{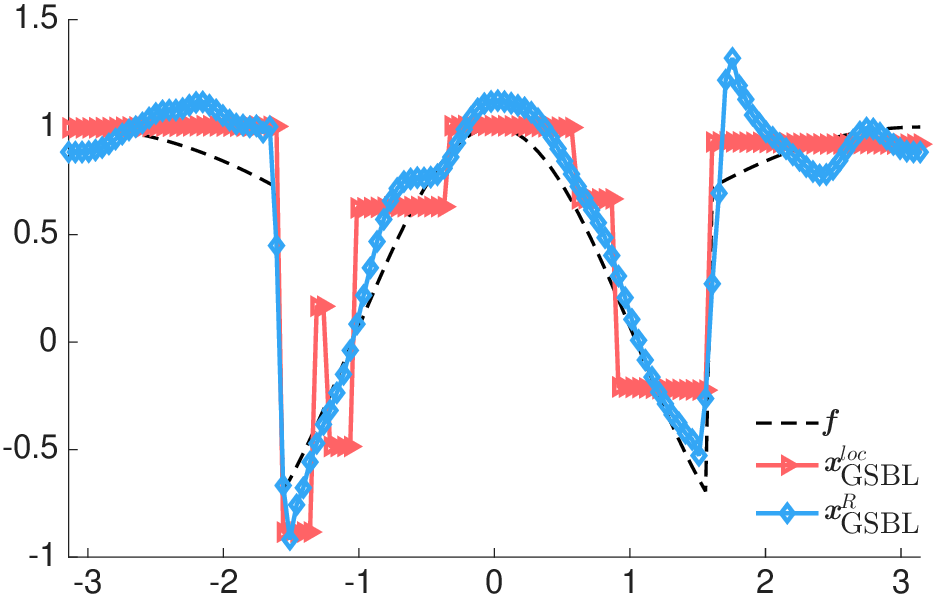}
    \end{subfigure}
    \begin{subfigure}[b]{.32\textwidth}
        \includegraphics[width=\textwidth]{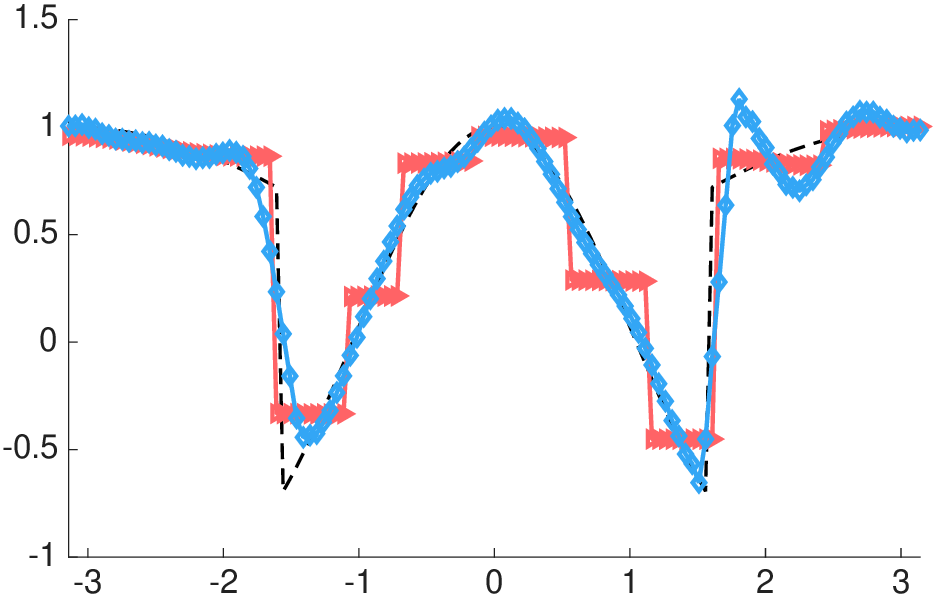}
    \end{subfigure}
    \begin{subfigure}[b]{.32\textwidth}
        \includegraphics[width=\textwidth]{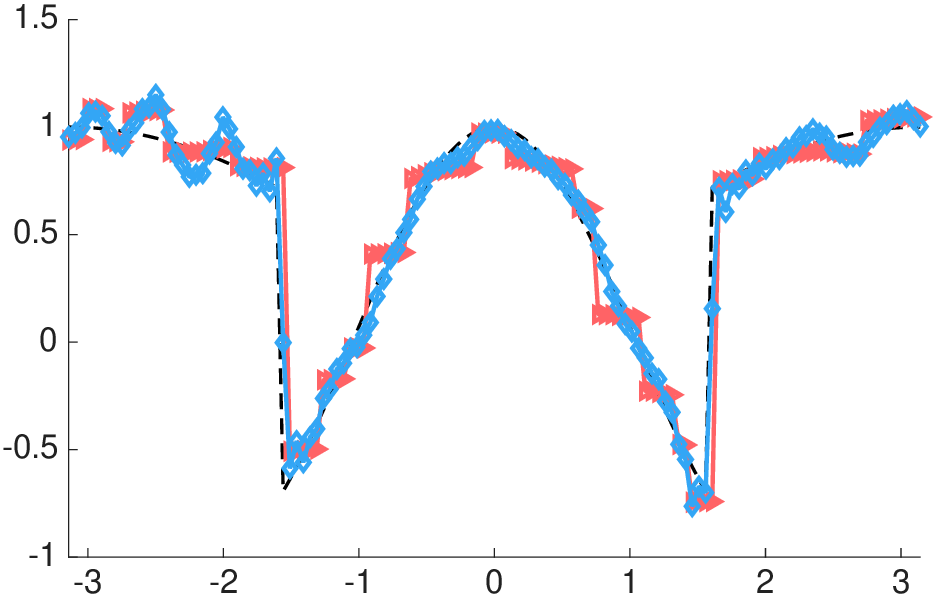}
    \end{subfigure}\\[.5em]
    \begin{subfigure}[b]{.32\textwidth}
        \includegraphics[width=\textwidth]{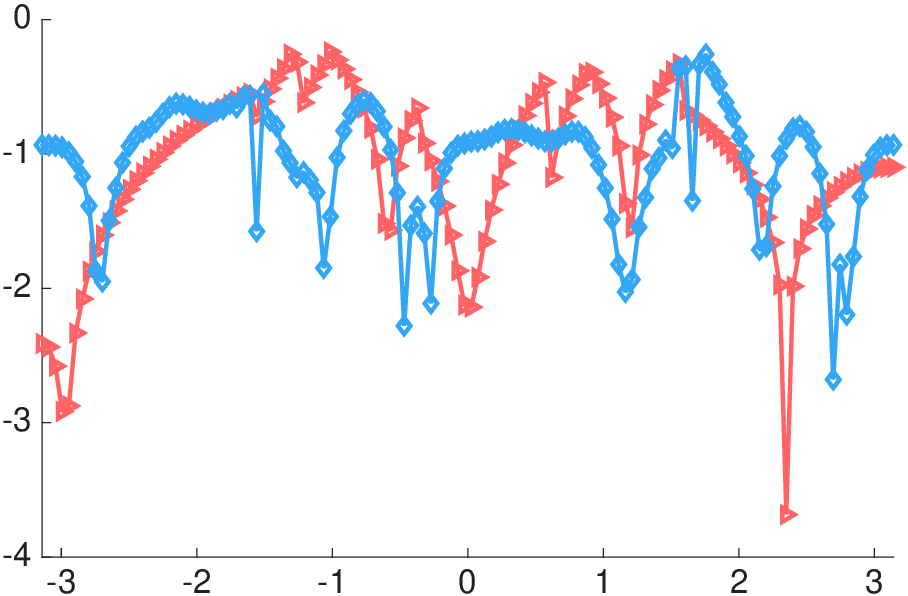}
    \end{subfigure}
    \begin{subfigure}[b]{.32\textwidth}
        \includegraphics[width=\textwidth]{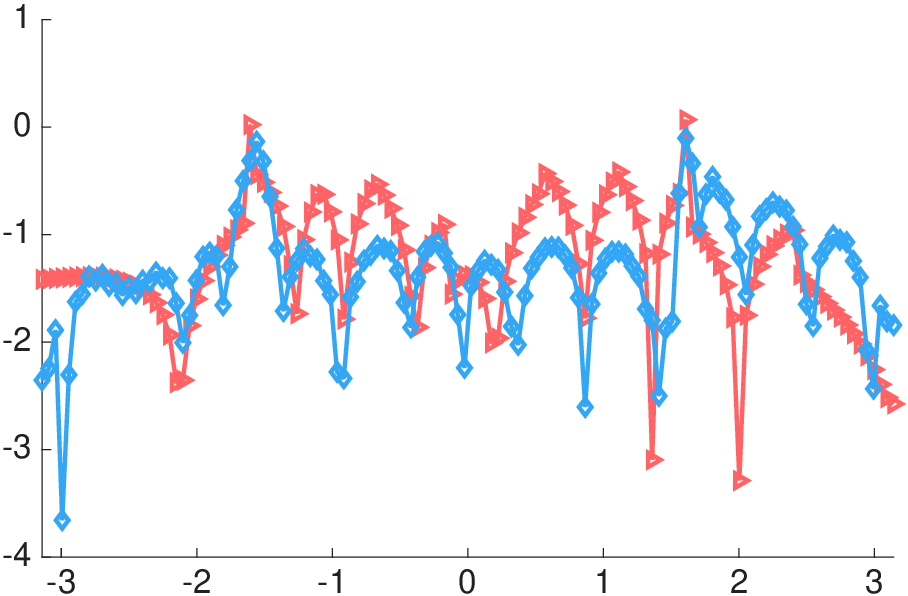}
    \end{subfigure}
    \begin{subfigure}[b]{.32\textwidth}
        \includegraphics[width=\textwidth]{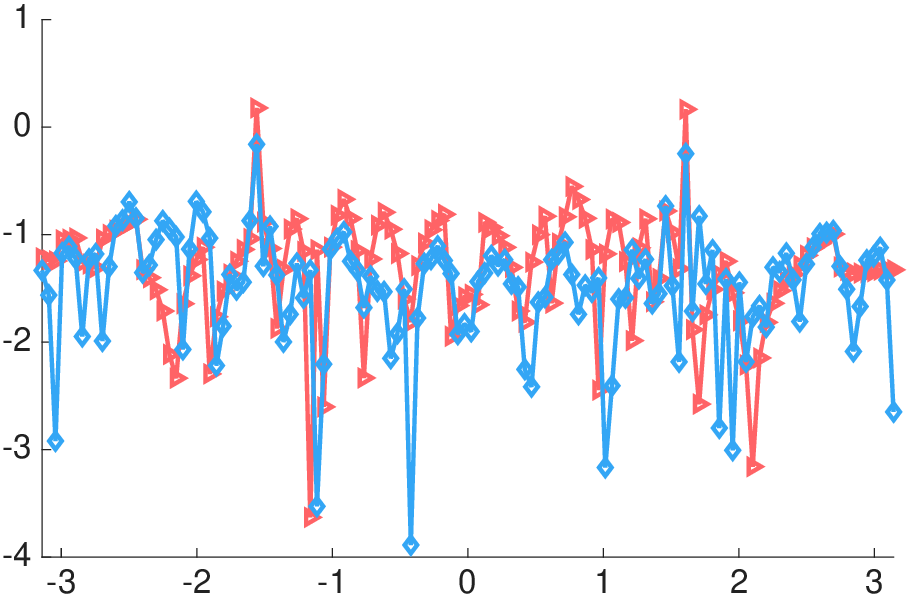}
    \end{subfigure}\\[.5em]
    \begin{subfigure}[b]{.32\textwidth}
        \includegraphics[width=\textwidth]{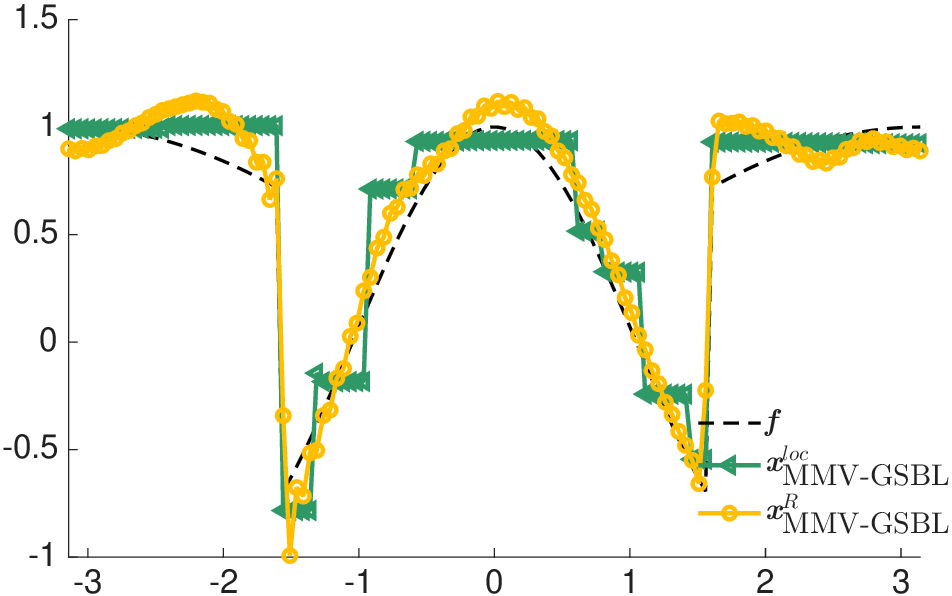}
    \end{subfigure}
    \begin{subfigure}[b]{.32\textwidth}
        \includegraphics[width=\textwidth]{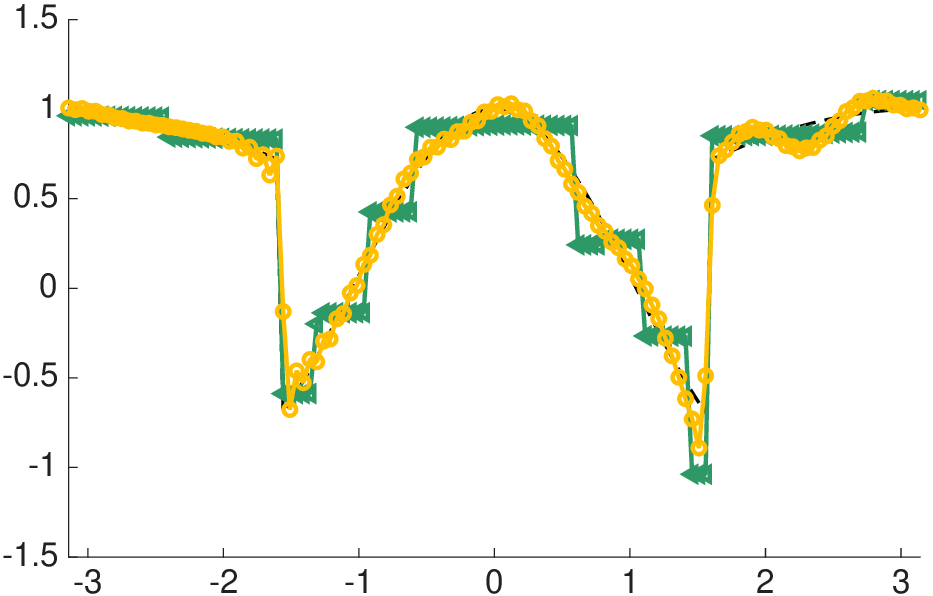}
    \end{subfigure}
    \begin{subfigure}[b]{.32\textwidth}
        \includegraphics[width=\textwidth]{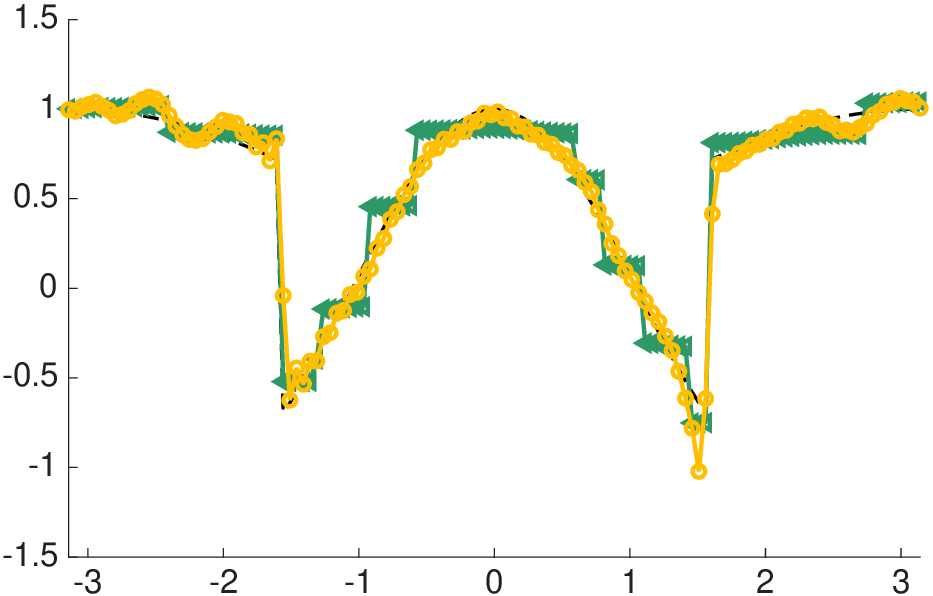}
    \end{subfigure}\\[.5em]
    \begin{subfigure}[b]{.32\textwidth}
        \includegraphics[width=\textwidth]{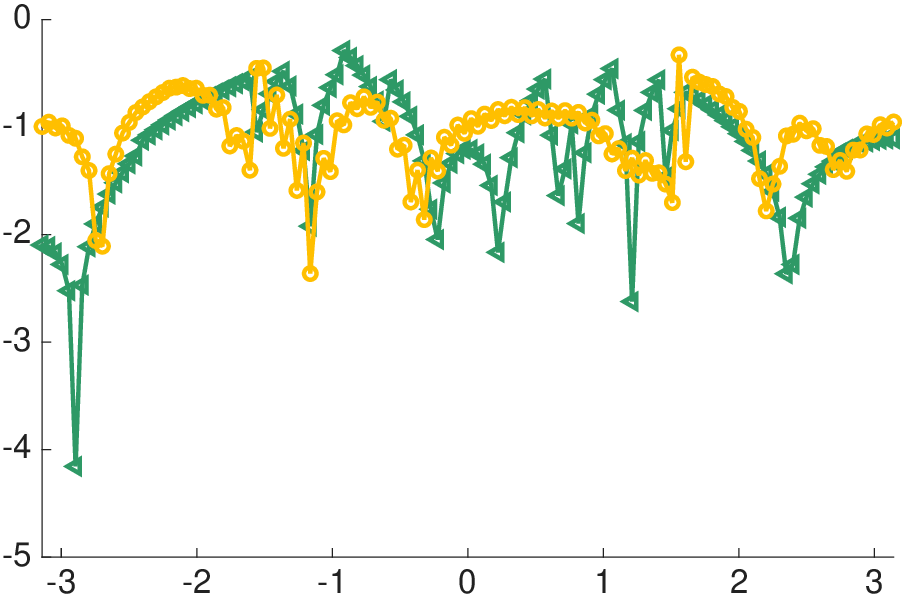}
    \end{subfigure}
    \begin{subfigure}[b]{.32\textwidth}
        \includegraphics[width=\textwidth]{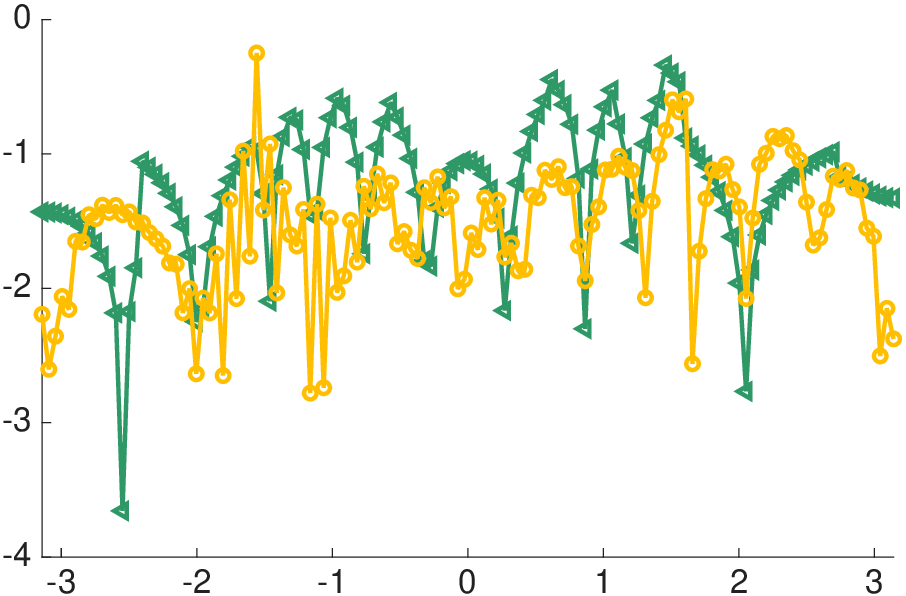}
    \end{subfigure}
    \begin{subfigure}[b]{.32\textwidth}
        \includegraphics[width=\textwidth]{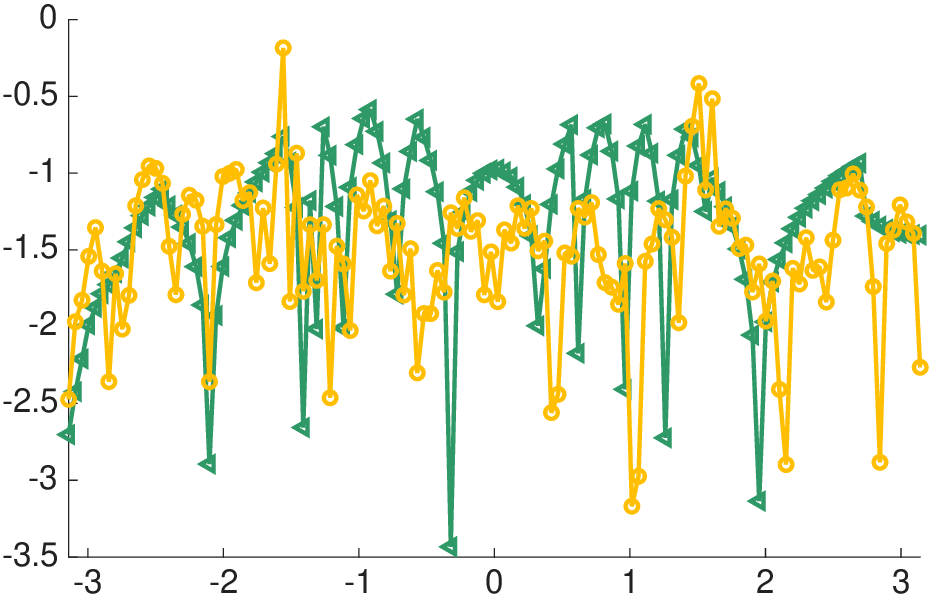}
    \end{subfigure}
    \caption{The multimodal measurement scenarios are consistent with those described in \cref{fig:trueSignals}, but here the recovery focuses exclusively on $f_3$. (top and  third rows) MAP estimates obtained by \cref{alg:GSBL,alg:MMV_GSBL} from data $\bm y$.  (second and bottom rows) Corresponding $\log_{10}E^{\text{abs}}_j$, $j=1,\dots,n$ of the top/third row solutions given by \eqref{eq: abs_err}.}
    \label{fig:rec_pw_multimodal}
\end{figure}

\subsubsection*{Piecewise constant signal recovery}
The first-order differencing prior transform $\Phi=T_n^p$ is specifically designed for reconstructing piecewise constant signals. To establish a benchmark for the residual prior transform $\Phi=R^p_{n,\zeta}$, we apply the same types of data collections used  for \Cref{ex: example_1d} to a collection of piecewise constant signals, defined as  
\begin{example}
\label{ex: example_1d_const}
    $f_4(s)$, $f_5(s)$, and $f_6    (s)$ are $2\pi-$periodic functions defined over $[-\pi,\pi)$ such that
    \begin{equation*}
        f_4(s)=\begin{cases}
        1.5, & s\in U_1\\
        \tfrac{-6}{\pi}, & s\in U_2\\
        1.5, & s\in U_3
        \end{cases},\quad
        f_5(s)=\begin{cases}
        -3, & s\in U_1\\
        \tfrac{-4}{\pi}, & s\in U_2\\
        -3, & s\in U_3
        \end{cases}, \quad
        f_6(s)=\begin{cases}
        .5, & s\in U_1\\
        \tfrac{-2}{\pi}, & s\in U_2\\
        .5, & s\in U_3
        \end{cases},
    \end{equation*}
    where $U_1=[-\pi,-\tfrac{\pi}{2})$, $U_2=[-\tfrac{\pi}{2},\tfrac{\pi}{2})$, and $U_3=[\tfrac{\pi}{2},\pi)$. Note that $f_4$, $f_5$, and $f_6$ share jump discontinuities at $-\tfrac{\pi}{2}$ and $\tfrac{\pi}{2}$. 
\end{example}

\cref{fig:rec_pw_consant} displays recovery results for piecewise constant functions defined in \cref{ex: example_1d_const}. Clearly, $\Phi=T_n^0$ provides the best recovery for each type of data modality, which is to be expected. To underscore the importance of matching the prior to the signal's structure, it is worth noting that a higher-order local differencing prior transform (e.g. $p=1$) would further degrade reconstructions. This detrimental effect of a mismatched prior on piecewise constant signals is demonstrated for Lasso regression in \cite{xiao2025new}.  While the residual prior transform $\Phi=R_{n,\zeta}^p$ is not able to counteract all of the effects of noise, blur, or missing data, it still appears robust, making it a practical choice for real-world applications where signal smoothness is typically unknown and rarely piecewise constant.

\subsubsection*{Employing joint sparsity in signal recovery}
We now consider the case where different modalities are used to measure the same underlying signal and we seek to determine if incorporating a joint prior can improve the recovery of each individual recovery via \cref{alg:MMV_GSBL}.  \Cref{fig:rec_pw_multimodal} displays the recovery of $f_3$ from \cref{ex: example_1d} under the same  set of data acquisitions shown in \Cref{fig:trueSignals}. The top row shows results using \cref{alg:GSBL}, that is, when infomation is considered separately, while the third row shows the joint recovery for each MAP estimate.  We observe that the staircasing artifacts persist when using $\Phi=T_n^p$, suggesting  that exploiting joint information cannot compensate for issues arising from incorrect prior assumptions. 
Conversely, when $\Phi=R_{n,\zeta}^p$  is employed, \cref{alg:MMV_GSBL} effectively leverages joint sparsity to suppress oscillations within smooth regions, yielding a significantly improved reconstruction. The corresponding spatial pointwise $\log_{10}$ errors for the top and third row recoveries are provided in the second and fourth rows, respectively.

\begin{figure}[h!]
    \centering
    \begin{subfigure}[b]{.32\textwidth}
        \includegraphics[width=\textwidth]{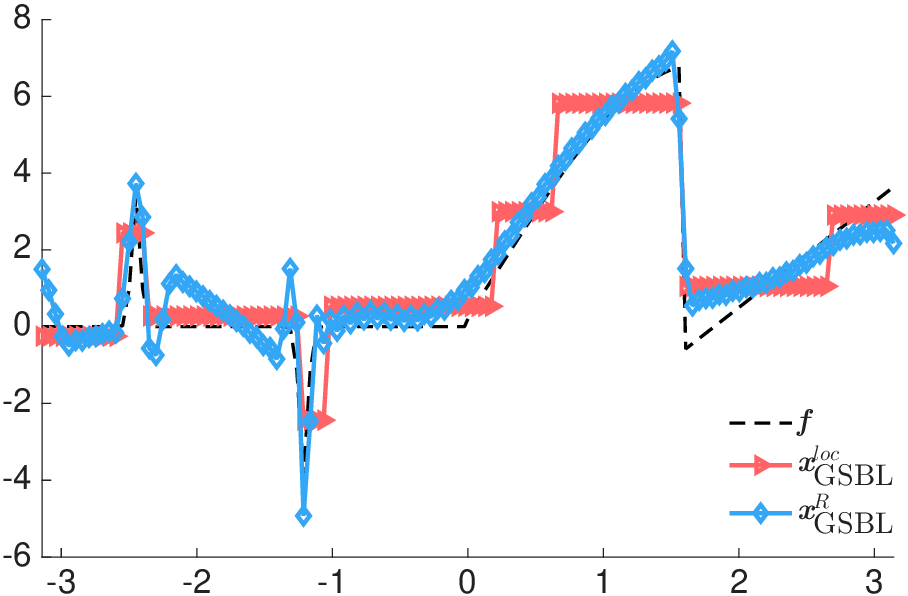}
    \end{subfigure}
    \begin{subfigure}[b]{.32\textwidth}
        \includegraphics[width=\textwidth]{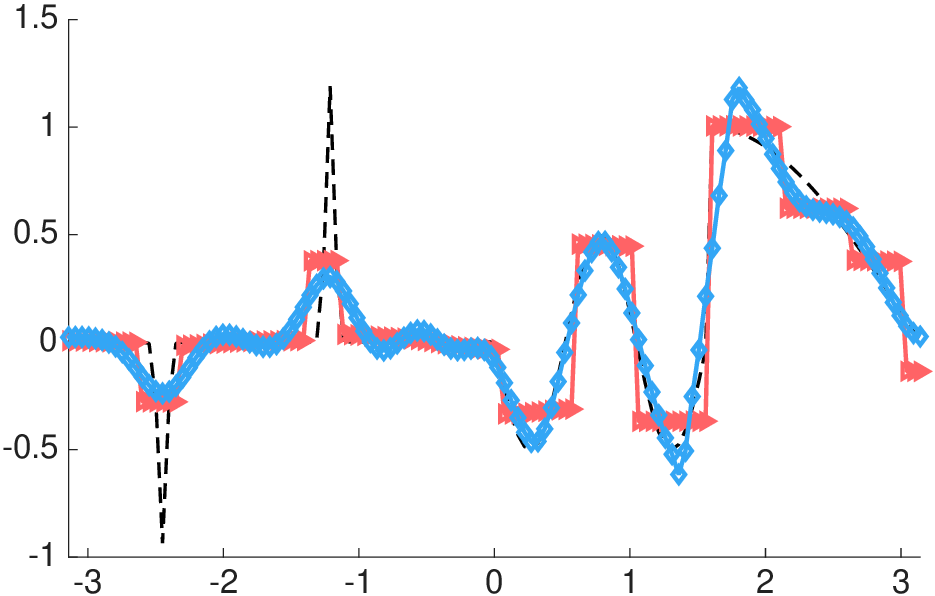}
    \end{subfigure}
    \begin{subfigure}[b]{.32\textwidth}
        \includegraphics[width=\textwidth]{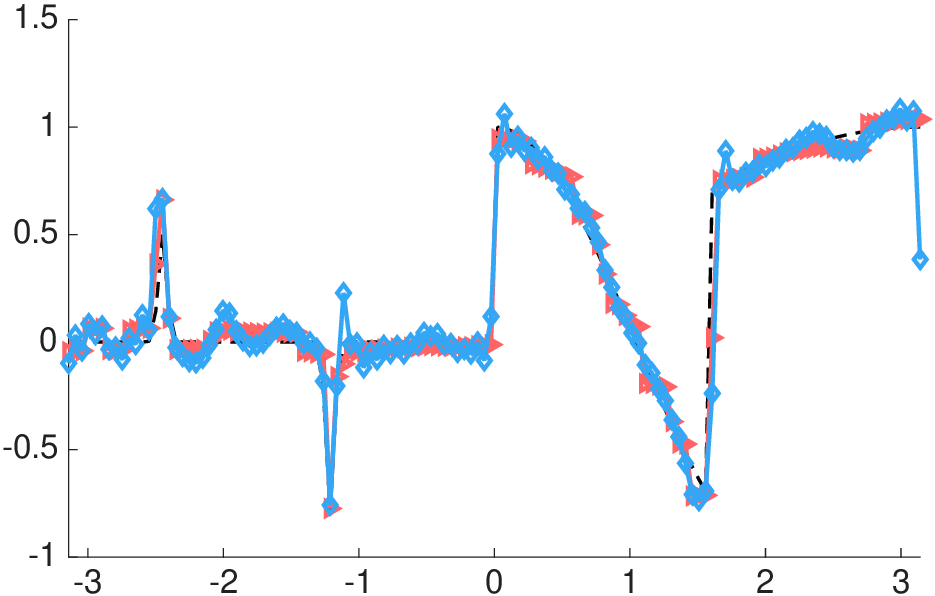}
    \end{subfigure}\\[.5em]
    \begin{subfigure}[b]{.32\textwidth}
        \includegraphics[width=\textwidth]{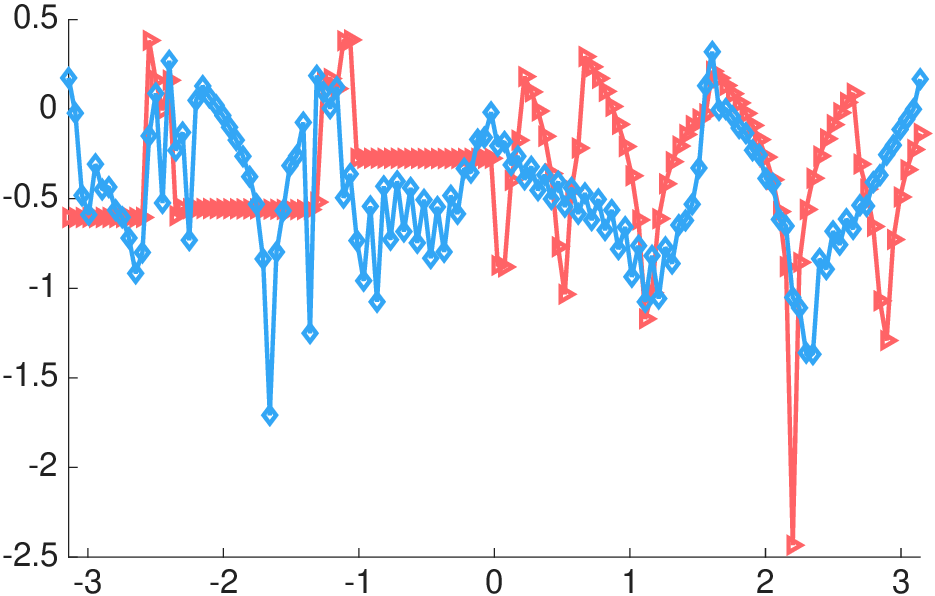}
    \end{subfigure}
    \begin{subfigure}[b]{.32\textwidth}
        \includegraphics[width=\textwidth]{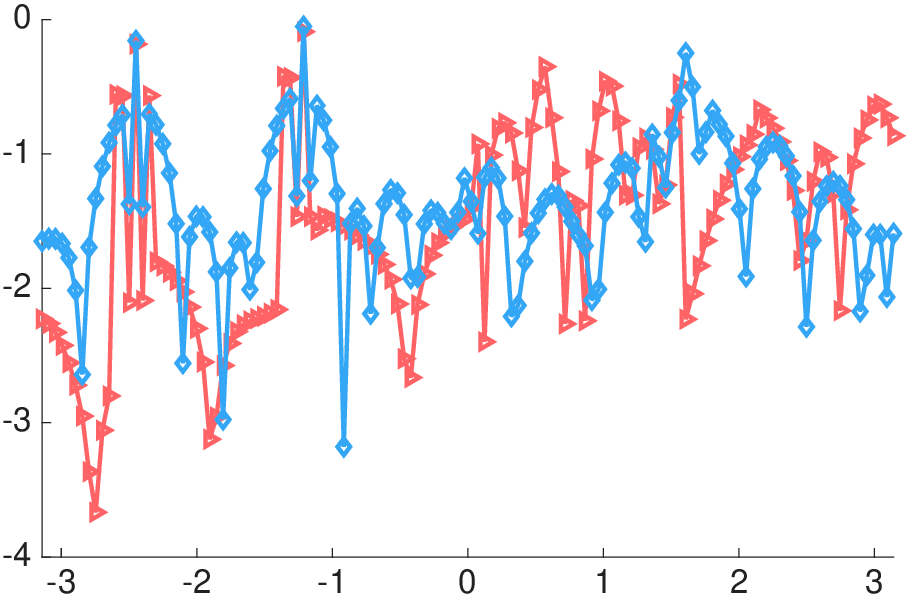}
    \end{subfigure}
    \begin{subfigure}[b]{.32\textwidth}
        \includegraphics[width=\textwidth]{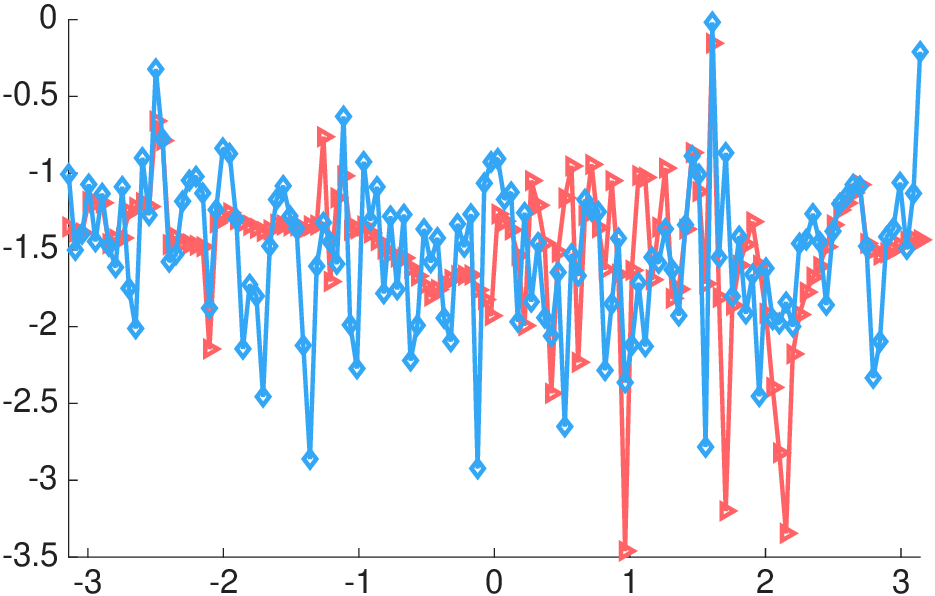}
    \end{subfigure}
    \caption{(top) MAP estimates obtained by \cref{alg:GSBL} from data $\bm y$. (bottom) Corresponding $\log_{10}E^{\text{abs}}_j$, $j=1,\dots,n$ of the top row solutions given by \eqref{eq: abs_err}. }
    \label{fig:rec_pw_combo}
\end{figure}

\subsubsection*{Mixed signal variability}
Our final 1D scenario displayed in \cref{fig:rec_pw_combo} features a collection of functions that are partially sparse and partially piecewise smooth. This structural composition is substantially different in its variability from piecewise constant or polynomial examples defined in \cref{ex: example_1d,ex: example_1d_const}. For each function, the interval $s_j\in[-\pi,0)$ contains non-zero values only at $j=15$ and $j=40$, thus representing a sparse signal component. Moreover, the interval $s_j\in[0,\pi)$ is populated by the piecewise smooth functions (left) $f_1$, (middle) $f_2$, and (right) $f_3$ in \cref{ex: example_1d}. Consistent with observations from \cref{fig:rec_pw_consant}, reconstructions using $\Phi=T_n^p$ tend to flatten over the sparse region but introduce severe staircasing artifacts where the underlying functions exhibit smooth variability. In contrast, recovery with $\Phi=R_{n,\zeta}^p$ successfully captures the nuanced variable behaviors in the smooth regions while maintaining minimal oscillations over the sparse intervals. We note that exploiting joint sparsity via \cref{alg:MMV_GSBL} offers no discernible improvement over the individual recoveries shown in \cref{fig:rec_pw_combo},  regardless of whether $\Phi=T^p_n$ or $\Phi=R^p_{n,\zeta}$ was used.  
\subsection{Joint Recovery of Multimodal Two-dimensional Images}
\label{sec:numerics2d}

We now extend our analysis to 2D images. This allows us to evaluate the performance and effectiveness of the residual prior transform in a higher-dimensional setting.
\begin{example}
    \label{ex: example_2d}
   Consider  2D-periodic pattern images over the domain $[-\pi,\pi]^2$ pixelated for $j,j' = 1, \dots, n$ as 
    \begin{equation*} 
    h_1(j,j') =\begin{dcases} 
    \begin{multlined}[t] 
        3 + (\cos(j) + 1)^2 \\
        + (\cos(j') + 1)^2
    \end{multlined},  \\ 
    4 + \sin(4j) + \sin(4j'), \\ 
    \sin(j),
\end{dcases}\quad
    h_2(j,j') = \begin{cases} \sin(6j), \\ -0.3 \sin(6j'), \\ \sin(-j+\pi),  \end{cases}
    h_3(j,j') = \begin{cases} \cos(2 \rho(j,j')), & \rho(j,j') \in\mathcal{R}_1 \\ \cos(4 \rho(j,j')), & \rho(j,j') \in\mathcal{R}_2 \\ -0.5 \cos(j), & \rho(j,j') \in\mathcal{R}_3 \end{cases}, 
    \end{equation*} where $\rho(j,j') = \sqrt{j^2 + j'^2}$ is the radius from the center of the domain with $\mathcal{R}_1=[0,0.3\pi)$, $\mathcal{R}_2=[0.3\pi, 0.7\pi)$, and $\mathcal{R}_3=[0.7\pi,\pi]$.
\end{example}

\begin{figure}[h!]
    \centering
    \begin{subfigure}[b]{.32\textwidth}
        \includegraphics[width=\textwidth]{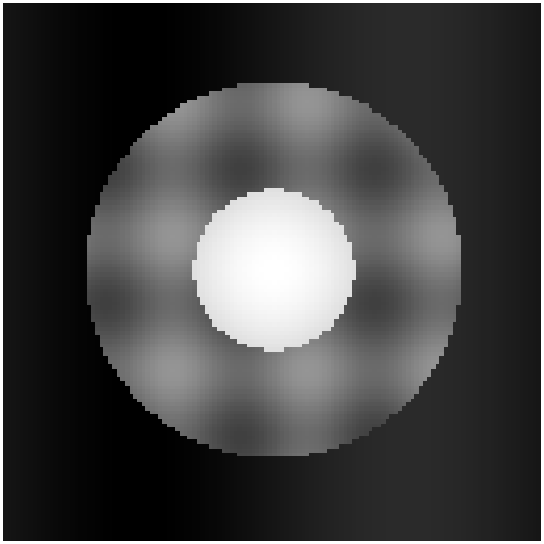}
    \end{subfigure}
    \begin{subfigure}[b]{.32\textwidth}
        \includegraphics[width=\textwidth]{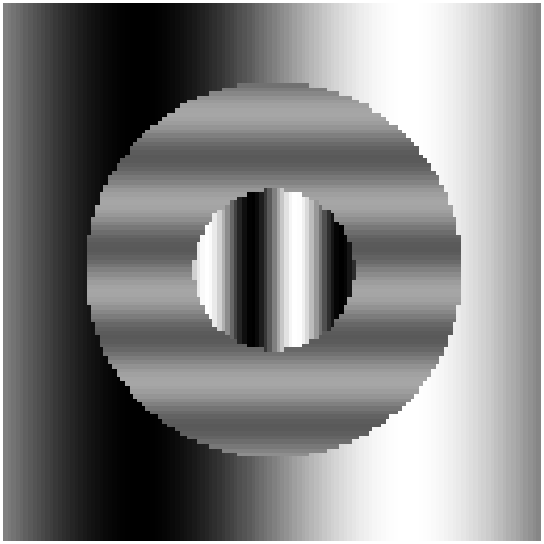}
    \end{subfigure}
    \begin{subfigure}[b]{.32\textwidth}
        \includegraphics[width=\textwidth]{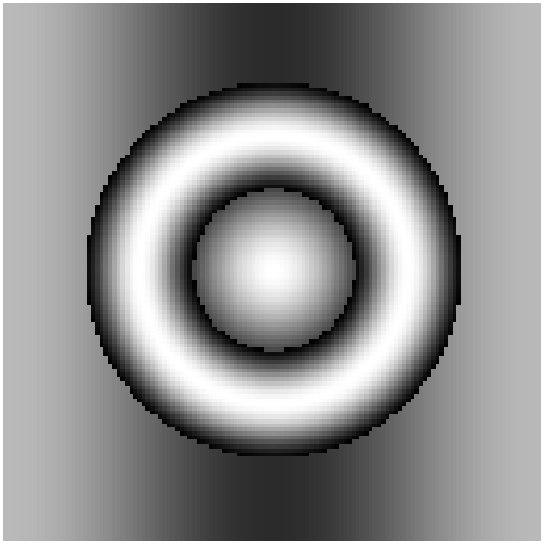}
    \end{subfigure}
    \caption{Discretized 2D images in \Cref{ex: example_2d}: (left) $\bm h_1$, (middle) $\bm h_2$, and (right) $\bm h_3$. }
    \label{fig:trueImages}
\end{figure}

\begin{figure}[h!]
    \centering
    \begin{subfigure}[b]{.32\textwidth}
        \includegraphics[width=.8\textwidth]{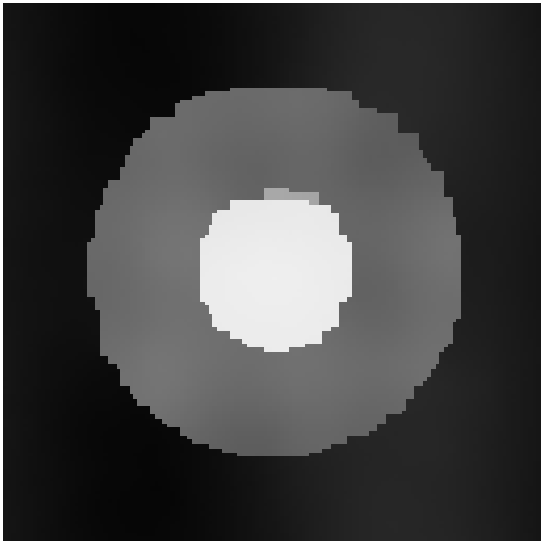}
    \end{subfigure}
    \begin{subfigure}[b]{.32\textwidth}
        \includegraphics[width=.8\textwidth]{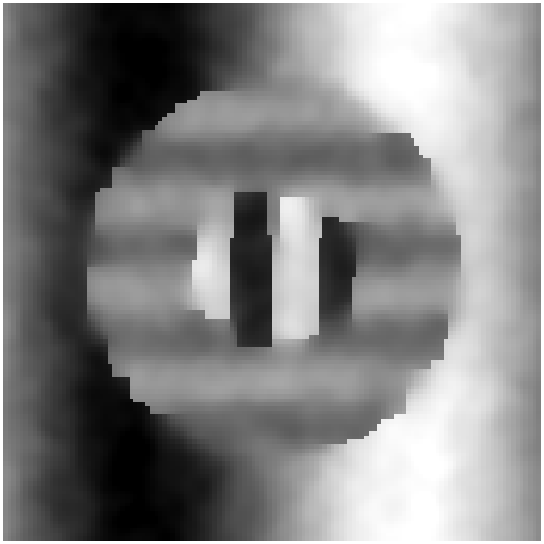}
    \end{subfigure}
    \begin{subfigure}[b]{.32\textwidth}
        \includegraphics[width=.8\textwidth]{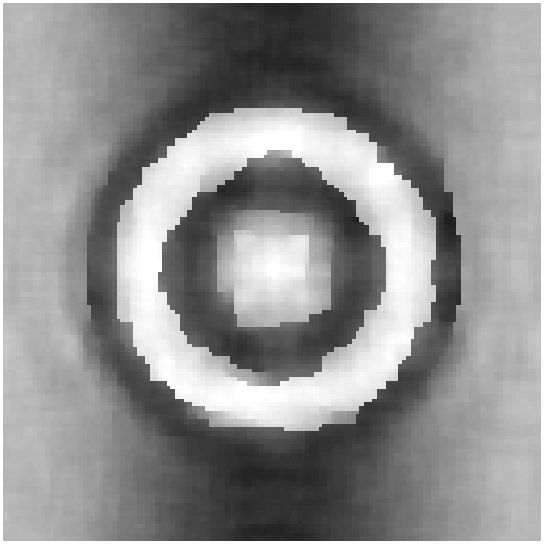}
    \end{subfigure}\\[.5em]
    \begin{subfigure}[b]{.32\textwidth}
        \includegraphics[width=.8\textwidth]{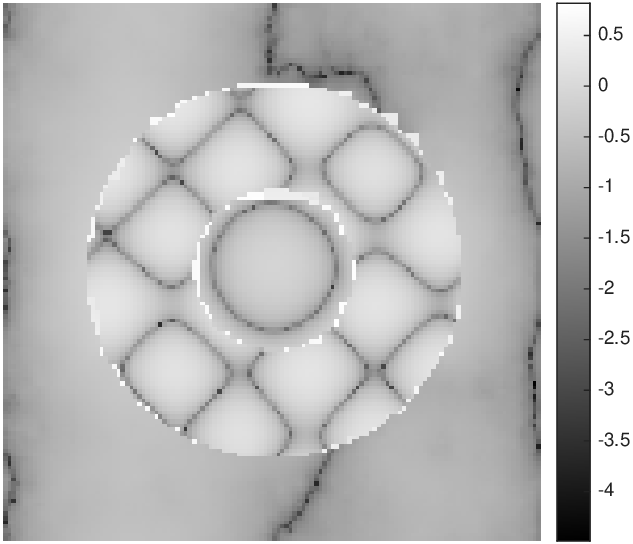}
    \end{subfigure}
    \begin{subfigure}[b]{.32\textwidth}
        \includegraphics[width=.8\textwidth]{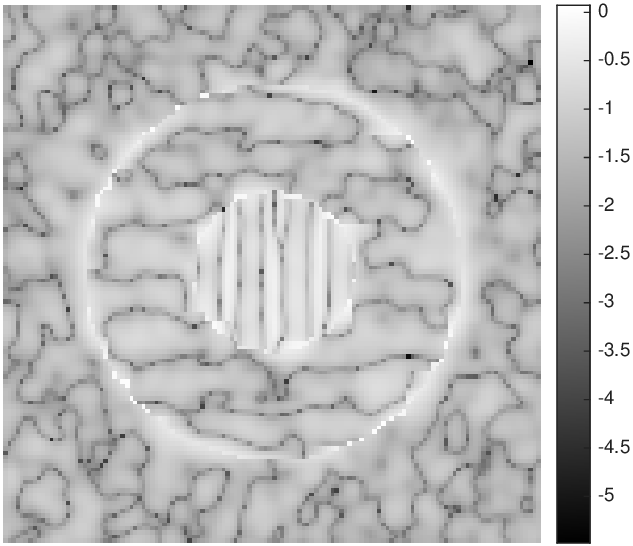}
    \end{subfigure}
    \begin{subfigure}[b]{.32\textwidth}
        \includegraphics[width=.8\textwidth]{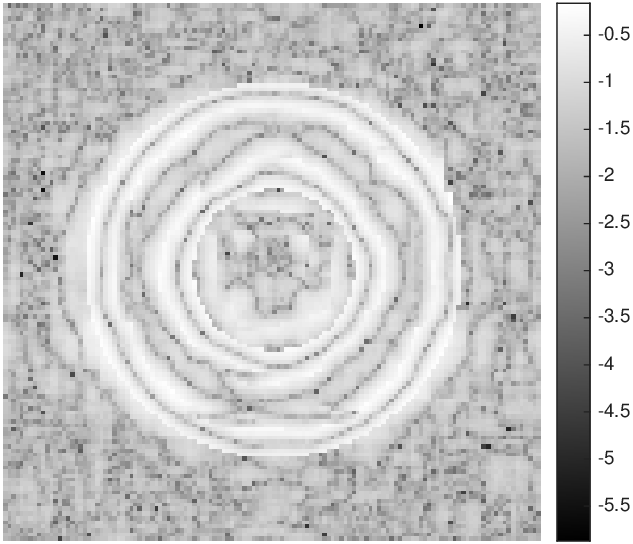}
    \end{subfigure}\\[.5em]
    \begin{subfigure}[b]{.32\textwidth}
        \includegraphics[width=.8\textwidth]{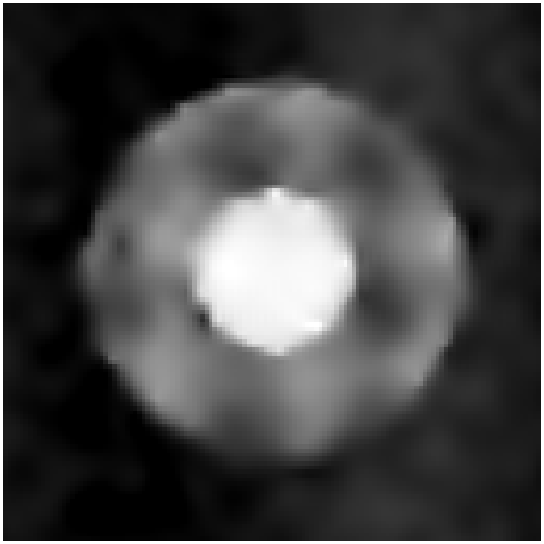}
    \end{subfigure}
    \begin{subfigure}[b]{.32\textwidth}
        \includegraphics[width=.8\textwidth]{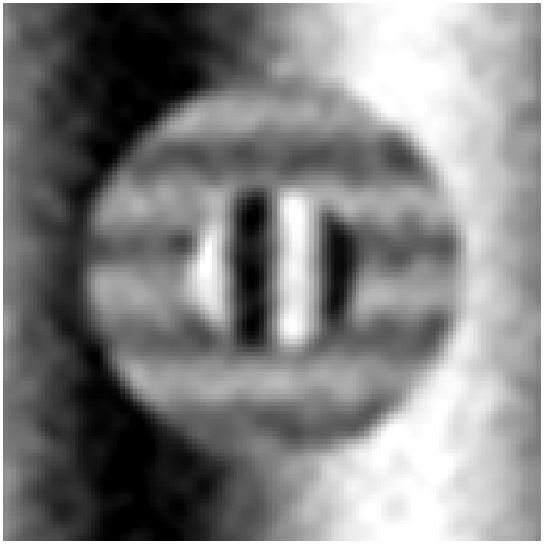}
    \end{subfigure}
    \begin{subfigure}[b]{.32\textwidth}
        \includegraphics[width=.8\textwidth]{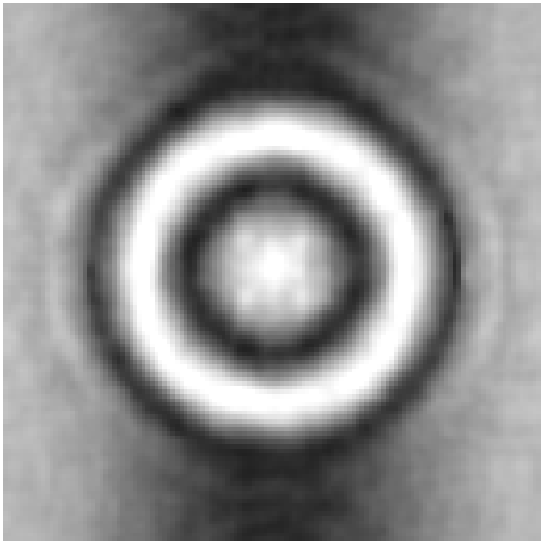}
    \end{subfigure}\\[.5em]
    \begin{subfigure}[b]{.32\textwidth}
        \includegraphics[width=.8\textwidth]{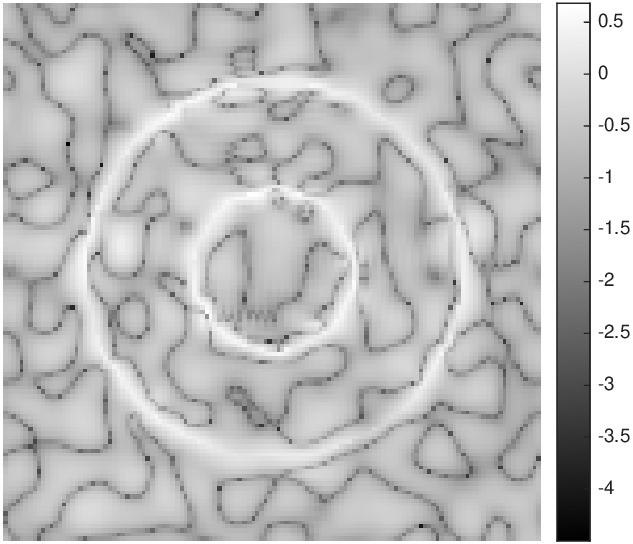}
    \end{subfigure}
    \begin{subfigure}[b]{.32\textwidth}
        \includegraphics[width=.8\textwidth]{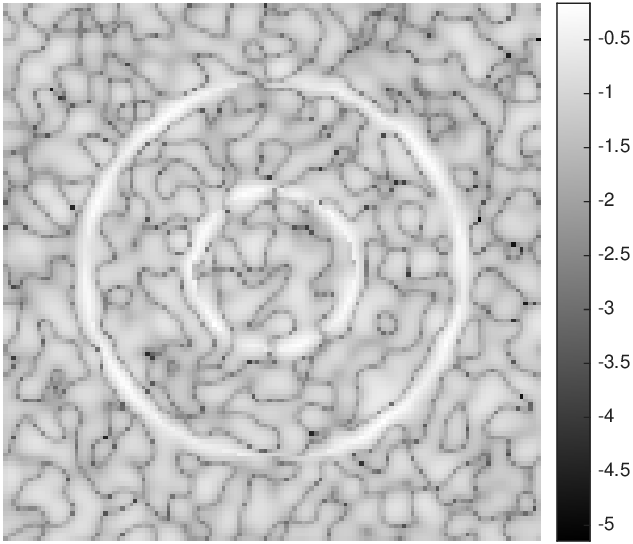}
    \end{subfigure}
    \begin{subfigure}[b]{.32\textwidth}
        \includegraphics[width=.8\textwidth]{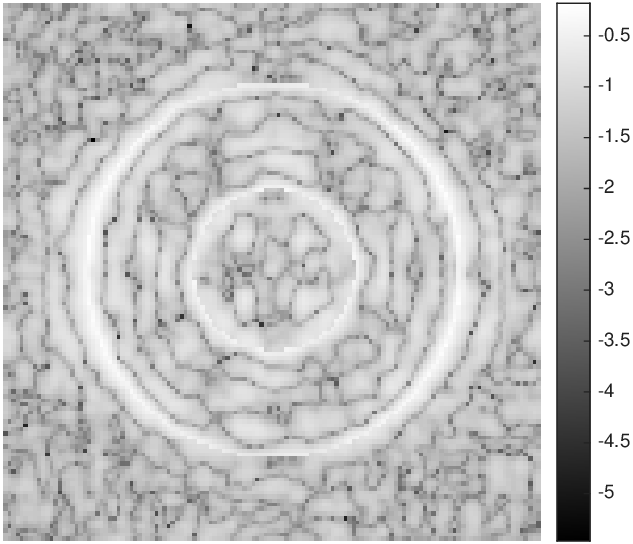}
    \end{subfigure}
    \caption{Recoveries of images in \cref{ex: example_2d} using \cref{alg:GSBL}: (left) $\bm h_1$ measured with $F_1=H_1I_n$ when $r_1=0.3$. (middle column) $\bm h_2$ measured with $F_2$ using a psf parameter $\gamma=0.01$. (right) $\bm h_3$ measured with $F_3=H_3\hat F$ when $r_3=0.7$. The SNR of additive noise is 5dB. (top) MAP estimates and (second row) corresponding absolute log errors ($\log_{10}E^\text{abs}$) for the standard first-order local differencing prior transform ($T_n^p$). (third row) MAP estimates and (bottom) corresponding absolute log errors ($\log_{10}E^\text{abs}$) for the residual prior transform ($R_{n,\zeta}^p$).}
    \label{fig:sep_rec_img}
\end{figure}

\begin{figure}[h!]
    \centering
    \begin{subfigure}[b]{.32\textwidth}
        \includegraphics[width=\textwidth]{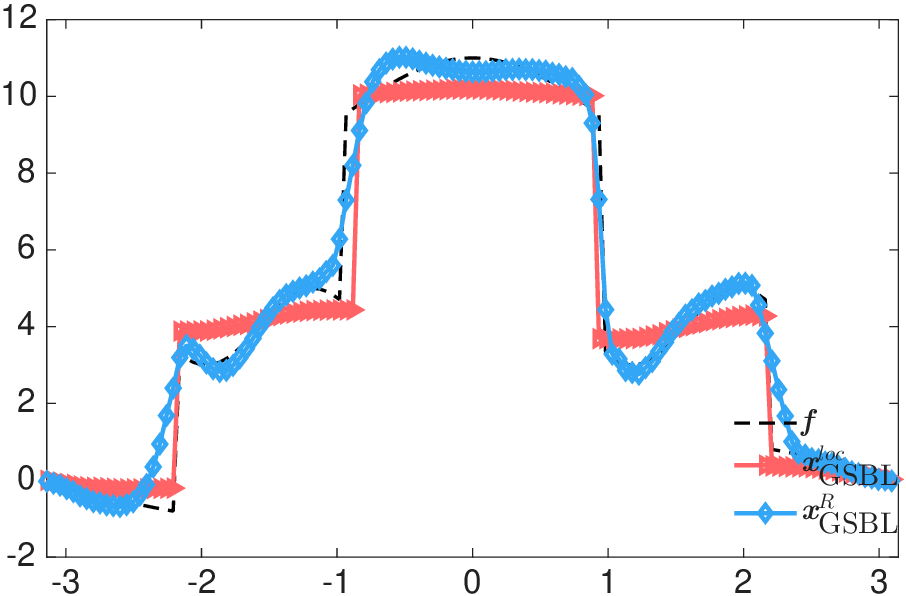}
    \end{subfigure}
    \begin{subfigure}[b]{.32\textwidth}
        \includegraphics[width=\textwidth]{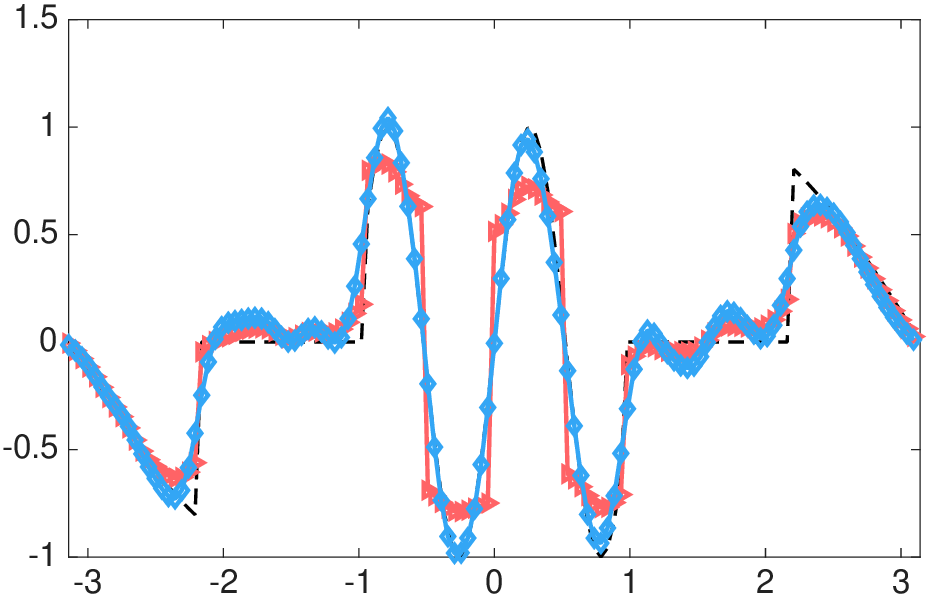}
    \end{subfigure}
    \begin{subfigure}[b]{.32\textwidth}
        \includegraphics[width=\textwidth]{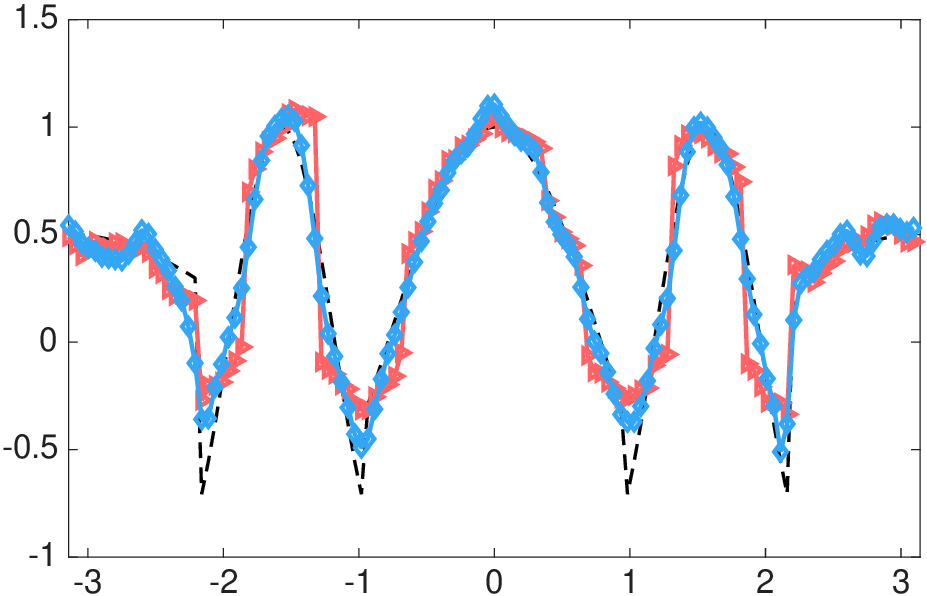}
    \end{subfigure}
    \caption{Horizontal slices of the MAP estimates  for (left) $\bm h_1(64,j')$, (middle) $\bm h_2(64,j')$, and (right) $\bm h_3(64,j')$, $j' = 1,\dots,n$. These slices correspond to the reconstructions detailed in \cref{fig:sep_rec_img}.}
    \label{fig:slice_2d}
\end{figure}

\begin{remark}
    Due to its separable nature, the residual prior transform extension to 2D  is straightforward.  Specifically, $R_{n,\zeta}^p$ can be applied  independently along each of the image's axes, without having to consider directional information. This separable application robustly isolates discontinuities, such as edges, in the 2D data by effectively suppressing the smooth variations in each direction independently. Moreover, although not discussed in any applications here, the core principles of the 1D transform matrix can be applied efficiently to any higher-dimensional problem, and over any range of pixels.
\end{remark}

As observed in \cref{fig:trueImages},  while each image in \Cref{ex: example_2d} has a unique pattern, they share common edges at radii of $0.3\pi$ and $0.7\pi$. \Cref{fig:sep_rec_img}  shows the corresponding MAP estimates obtained by \cref{alg:GSBL} with (top) $\Phi=T_n^p$ and (third row) $\Phi=R_{n,\zeta}^p$. Analogous to our observations in 1D case, the recoveries using $\Phi=T_n^p$ suffer from a pronounced staircasing effect, resulting in visibly ``chunky'' artifacts in the image.  Solutions obtained using $\Phi=R_{n,\zeta}^p$ successfully capture the intricate variability, achieving a significantly higher degree of structural fidelity. The second and bottom rows present the spatial pointwise $\log_{10}$ errors for the reconstructions with $\Phi=T^p_n$  and $\Phi=R_{n,\zeta}^p$ respectively. A horizontal slice $\bm h_l(64,j')$, $j'  = 1,\dots,n$,  in \cref{fig:slice_2d} further corroborates the observed performance differences.

\begin{figure}[h!]
    \centering
    \begin{subfigure}[b]{.45\textwidth}
        \includegraphics[width=\textwidth]{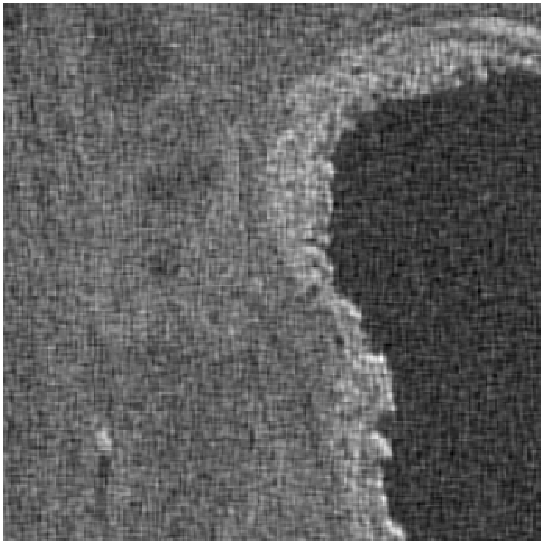}
    \end{subfigure}\hspace{2em}
    \begin{subfigure}[b]{.45\textwidth}
        \includegraphics[width=\textwidth]{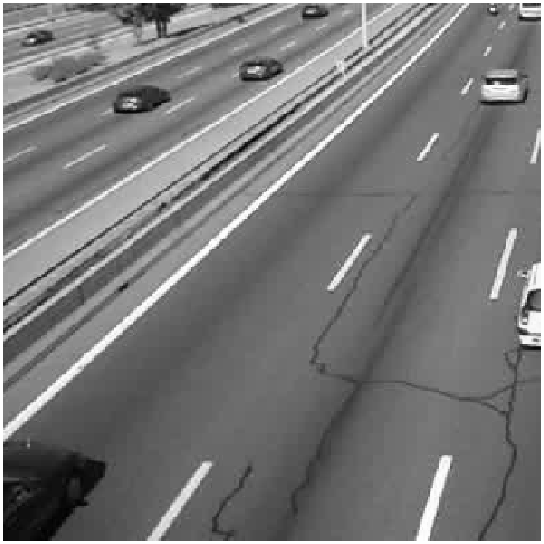}
    \end{subfigure}
    \caption{(left) SAR image of a golf course \cite{wang2007fast}.
    (right) A sample image from the GRAM road-traffic monitoring dataset \cite{guerrero2013vehicle}. }
    \label{fig:SAR_DRIVE_data}
\end{figure}

\begin{figure}[h!]
    \centering
    \begin{subfigure}[b]{.32\textwidth}
        \includegraphics[width=\textwidth]{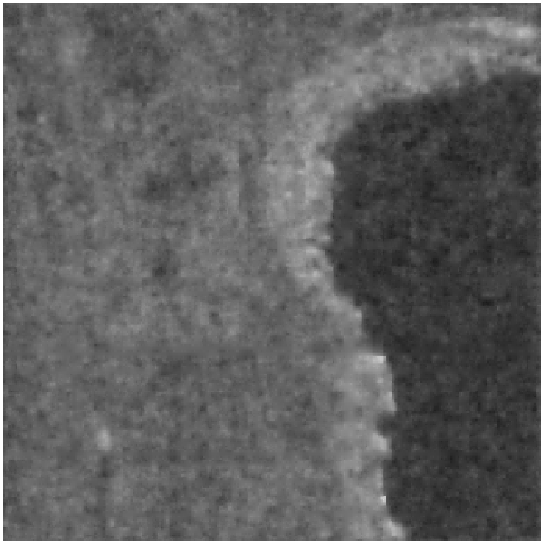}
    \end{subfigure}
    \begin{subfigure}[b]{.32\textwidth}
        \includegraphics[width=\textwidth]{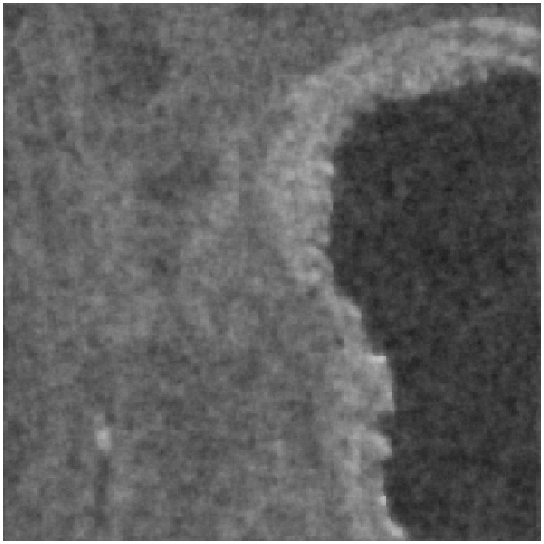}
    \end{subfigure}
    \begin{subfigure}[b]{.32\textwidth}
        \includegraphics[width=\textwidth]{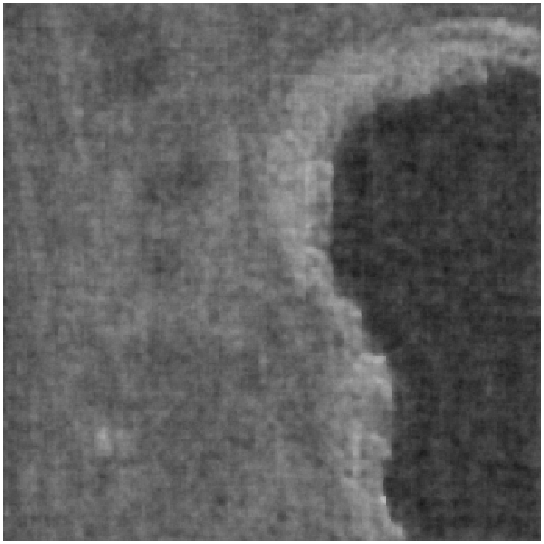}
    \end{subfigure}\\
    \begin{subfigure}[b]{.32\textwidth}
        \includegraphics[width=\textwidth]{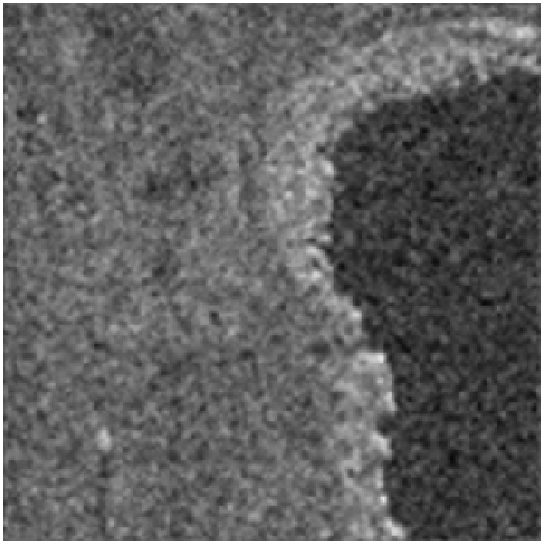}
    \end{subfigure}
    \begin{subfigure}[b]{.32\textwidth}
        \includegraphics[width=\textwidth]{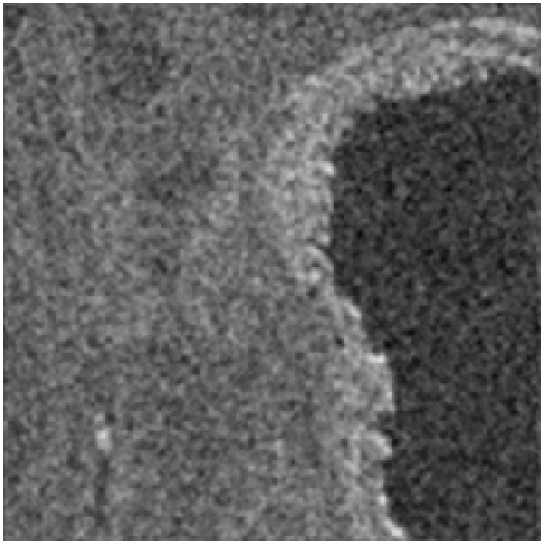}
    \end{subfigure}
    \begin{subfigure}[b]{.32\textwidth}
        \includegraphics[width=\textwidth]{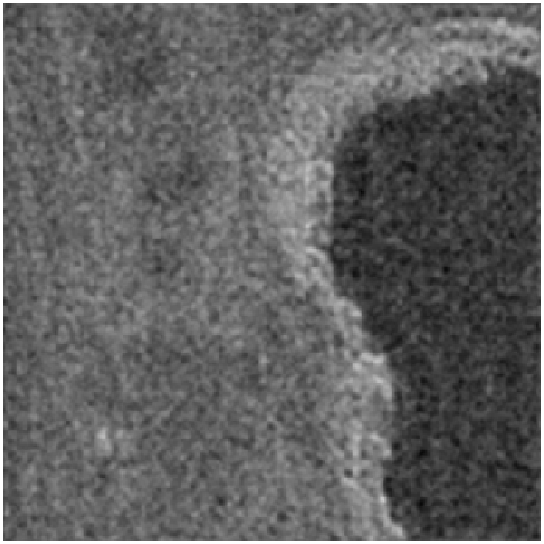}
    \end{subfigure}
    \caption{The resolution is $n=400$. MAP estimates of \cref{fig:SAR_DRIVE_data} (left) SAR image obtained by \cref{alg:MMV_GSBL} using prior transforms 
    (top) $\Phi=T_n^p$ and (bottom) $\Phi=R_{n,\zeta}^p$.}
    \label{fig:SAR}
\end{figure}

\begin{figure}[h!]
    \centering
    \begin{subfigure}[b]{.32\textwidth}
        \includegraphics[width=\textwidth]{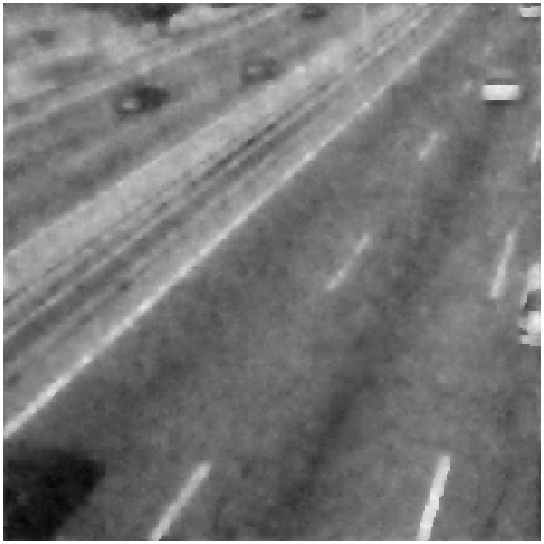}
    \end{subfigure}
    \begin{subfigure}[b]{.32\textwidth}
        \includegraphics[width=\textwidth]{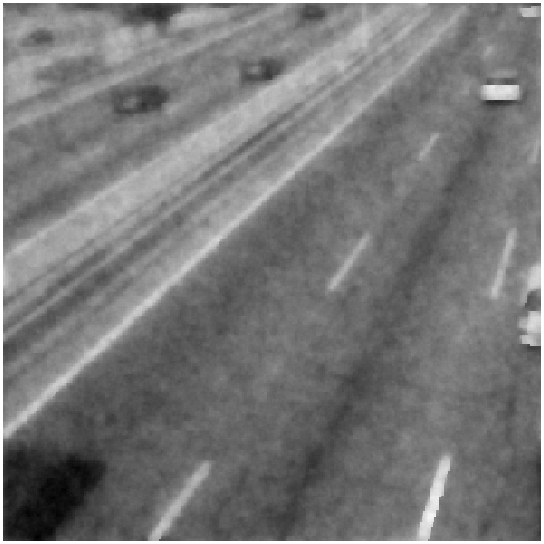}
    \end{subfigure}
    \begin{subfigure}[b]{.32\textwidth}
        \includegraphics[width=\textwidth]{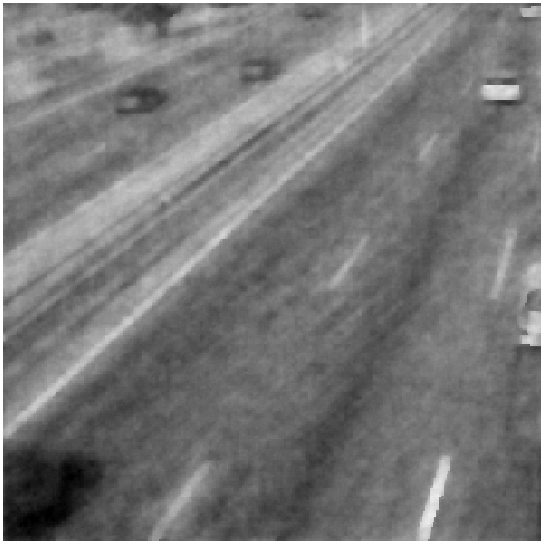}
    \end{subfigure}\\
    \begin{subfigure}[b]{.32\textwidth}
        \includegraphics[width=\textwidth]{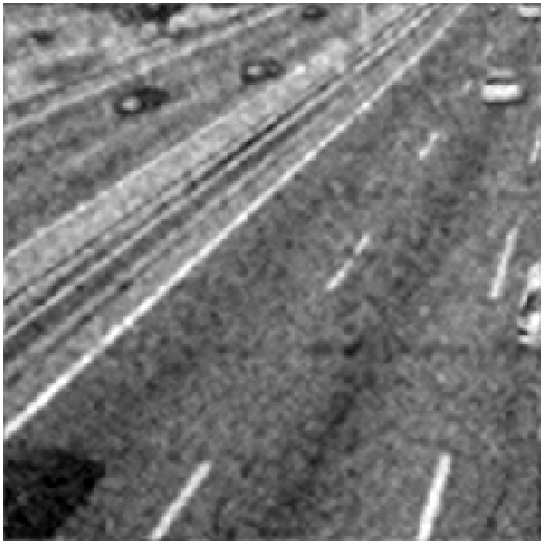}
    \end{subfigure}
    \begin{subfigure}[b]{.32\textwidth}
        \includegraphics[width=\textwidth]{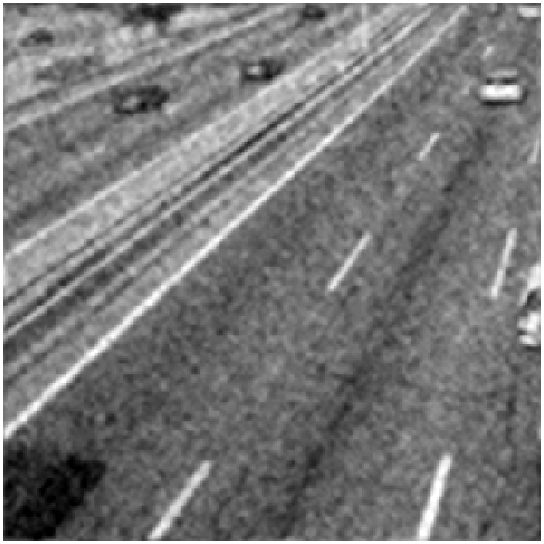}
    \end{subfigure}
    \begin{subfigure}[b]{.32\textwidth}
        \includegraphics[width=\textwidth]{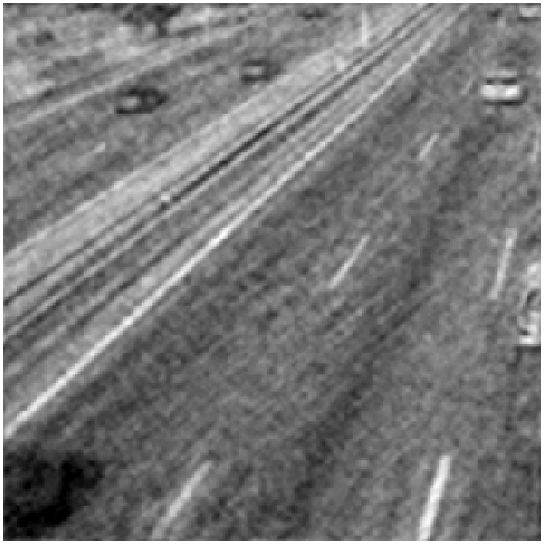}
    \end{subfigure}
    \caption{The resolution is $n=400$. MAP estimates of \cref{fig:SAR_DRIVE_data} (right) road-traffic footage image obtained by \cref{alg:MMV_GSBL} using prior transforms 
    (top) $\Phi=T_n^p$ and (bottom) $\Phi=R_{n,\zeta}^p$.}
    \label{fig:DRIVE}
\end{figure}

To demonstrate the practical applicability of the residual prior transform, we conclude our analysis with experiments on two real-world examples, presented in \Cref{fig:SAR_DRIVE_data}, namely (left) a complex SAR image of a golf course \cite{wang2007fast} and (right) a representative frame from the GRAM road-traffic monitoring dataset \cite{guerrero2013vehicle}.
Each source image is then subjected to a simulated multimodal data acquisition process, comprising three distinct degradation processes. The first modality consists of direct spatial undersampling by $F_1=H_1I_n$ with $r_1 = 0.3$, where 30\% of the image pixels are randomly measured. The second is a fully sampled, noisy observation of the entire image $F_2 = I_n$, posing a standard denoising task.  The third involves measuring a random 50\% subset of the image's Fourier coefficients using $F_3=H_3\hat F$ with $r_3 = 0.5$. For each scenario, the SNR of additive noise is 5dB.

\Cref{fig:SAR} presents a visual comparison of the MAP reconstructions for \Cref{fig:SAR_DRIVE_data} (left) SAR imagery, obtained by \cref{alg:MMV_GSBL}. The results contrast the performance of (top) $\Phi=T_n^p$ with (bottom) $\Phi=R_{n,\zeta}^p$. Consistent with previous findings, the reconstruction with $\Phi=R_{n,\zeta}^p$ introduces minor oscillatory artifacts into the constant dark area of the lake. However, this drawback is significantly outweighed by its ability to capture fine details of the landscapes in the scene, whereas $\Phi=T_n^p$ tends to produce blocky and over-smoothed results that obscure these structural features.

A visual comparison of the MAP estimates of \cref{fig:SAR_DRIVE_data} (right) traffic image obtained by \cref{alg:MMV_GSBL} is presented in \Cref{fig:DRIVE}, contrasting the performance of (top) $\Phi=T_n^p$ with (bottom) $\Phi=R_{n,\zeta}^p$. The reconstruction utilizing $\Phi=T_n^p$ suffers from pronounced blocky artifacts, which degrade image quality and result in a significant feature loss of vehicles in the scene. In contrast, the solution using $\Phi=R_{n,\zeta}^p$ demonstrates superior performance, successfully avoiding such artifacts and preserving the finer details and structural integrity of key features, such as the road and individual vehicles.

\section{Concluding remarks} 
\label{sec:summary} 

We have demonstrated the efficacy and robustness of the residual prior transform for solving ill-posed inverse problems within a hierarchical sparse Bayesian learning framework, particularly for signals and images with complex variability. The residual prior transform is so named because it is constructed from the residual between two distinct transform matrices, in our case, a first-order local differencing $\Phi_1=T_n^p$ and a global edge approximant $\Phi_2=S_{n,\zeta}^{\sigma_{2p+1}}$. While either prior transform used individually can introduce modeling errors if mismatched with true variability of the signal, the resulting residual prior transform $R_{n,\zeta}^p=T_n^p-S_{n,\zeta}^{\sigma_{2p+1}}$ is designed by construction to be small everywhere except at the discontinuities, effectively avoiding the errors caused by false smoothness assumptions. This makes the recovery process far less sensitive to the specific choice of smoothness order in the prior than traditional methods. 
A key contribution of this paper is the application of the residual prior transform to challenging multimodal problems, moving beyond its original single measurement design, particularly the joint recovery of multiple signals with shared structural information and the fusion of MMVs of one scene from diverse sources. The inherent flexibility of the residual prior transform allows for uniform application across multimodal datasets, obviating the need for a priori assumptions about each signal's specific smoothness. This greatly simplifies prior transform selection and enhances the robustness of reconstructions in complex, real-world applications.

Our numerical results provide compelling evidence of the practical advantages of the residual prior transform. In a wide array of tests, from highly variable synthetic signals and images to real-world SAR and traffic surveillance data, the residual prior transform consistently outperformed the standard first-order local differencing prior transform by successfully capturing complex structural features.  By contrast,   the staircasing artifacts inevitably introduced by the first-order local differencing prior transform can lead to misidentification of important structure.  Furthermore, the strong performance was maintained across several challenging multimodal recovery tasks from noisy, blurred, and spatially or spectrally undersampled data. Overall, the residual prior transform is demonstrated to be a practical and highly flexible tool for addressing realistic inverse problems.

While this paper illustrates the efficiency of the residual prior transform within a hierarchical sparse Bayesian learning framework for producing robust point estimates and uncertainty quantification, it also opens several promising avenues for future research. The separable nature of the 2D residual prior transform, which operates on one dimension at a time, enables its computationally efficient extension to higher dimensions, e.g.~3D volumetric multispectral data, while avoiding the need for complex non-separable operators.  Moreover, the idea of leveraging the residual of two operators is not limited to our specific choices of $\Phi_1$ and $\Phi_2$ in this paper. Exploring alternative operator pairings, such as wavelets, could yield priors that are better suited to other data types. To this end, while we considered only uniform discretizations, the residual prior transform is not inherently limited as such.  That is, the recovery could be obtained on points suitable for the given application. Finally, the residual prior transform concept can be extended from the spatial to the temporal domain, where a residual of two temporal operators could effectively identify key features in time series data, such as regions of temporal coherence or anomalies. This could prove valuable in applications involving physiological signals, such as analyzing heart rate impulses, or in video processing.


\section*{Acknowledgements}
This work was supported by  the DOD (ONR MURI) \#N00014-20-1-2595 and  DOE ASCR \#DEDE-SC0025555.

\bibliographystyle{siamplain}
\bibliography{literature}

\end{document}